%% file: min.tex
\input pmacro2.tex

\overfullrule=0pt
\vcorrection{-1cm}
\topmatter
\title
The $CR$ structure of minimal\\
orbits in complex flag manifolds
\endtitle
\author
Andrea Altomani \\
Costantino Medori \\
Mauro Nacinovich
\endauthor
\address
A. Altomani:
Dipartimento di Matematica, II Universit\`a di Roma ``Tor Vergata'', 
Via della Ricerca Scientifica, 00185 Roma (Italy)
\endaddress
\email altomani\@mat.uniroma2.it\endemail
\thanks During the preparation of this paper
the first Author has been partially supported by a grant of the
University of Basilicata.  
This research has been partially supported by M.I.U.R. \endthanks
\address
C.Medori,
Dip. di Matematica, Universit\`a di Parma,
Parco Area delle Scienze, 53/A
- 43100 Parma - Italy
\endaddress
\email medori\@prmat.math.unipr.it\endemail
\address
M. Nacinovich:
Dipartimento di Matematica, II Universit\`a di Roma ``Tor Vergata'', 
Via della Ricerca Scientifica, 00185 Roma (Italy)
\endaddress
\email nacinovi\@mat.uniroma2.it\endemail
\rightheadtext{$CR$ structure of minimal orbits}
\leftheadtext{A.Altomani\quad C.Medori\quad M.Nacinovich}
\date{July 14, 2005}\enddate
\abstract 
Let $\Hat{\bold{G}}$ be a complex semisimple Lie group, $\bold{Q}$ a parabolic
subgroup and $\bold{G}$ a real form of $\Hat{\bold{G}}$.
The flag manifold $\Hat{\bold{G}}/\bold{Q}$ decomposes into finitely many
$\bold{G}$-orbits; among them there is exactly one orbit of minimal dimension,
which is compact.
We study these minimal orbits from the point of view of $CR$ geometry.
In particular we characterize those minimal orbits that are of finite type and
satisfy various nondegeneracy conditions, compute their fundamental group and
describe the space of their global $CR$ functions.
Our main tool are parabolic $CR$ algebras, which give an infinitesimal 
description of the $CR$ structure of minimal orbits.
\endabstract
\toc\widestnumber\head{\S 14}
\head{\S 1.} Introduction\endhead
\head{\S 2.} Preliminaries on $CR$ manifolds\endhead
\head{\S 3.} The minimal orbit in a complex flag manifold\endhead
\head{\S 4.} Homogeneous $CR$ manifolds and $CR$ algebras\endhead
\head{\S 5.} Parabolic $CR$ algebras\endhead
\head{\S 6.} Parabolic minimal $CR$ algebras and cross-marked Satake
\linebreak 
{diagrams}\endhead
\head{\S 7.} $\frak{g}$-equivariant fibrations\endhead
\head{\S 8.} A rigidity theorem for minimal orbits\endhead
\head{\S 9.} Fundamental $CR$ algebras\endhead
\head{\S 10.} Totally real and totally complex $CR$ algebras\endhead
\head{\S 11.} Weak nondegeneracy\endhead
\head{\S 12.} Strict nondegeneracy\endhead
\head{\S 13.} Essential pseudoconcavity for minimal orbits\endhead
\head{\S 14.} $CR$ functions on minimal orbits\endhead
\head{} Appendix: Table of noncompact real forms and Satake diagrams\endhead
\endtoc
\endtopmatter

\document
\input 01min.tex

\input 02min.tex

\input 03min.tex

\input 04min.tex

\input 05min.tex
\input 06min.tex
\input 07min.tex
\input 08min.tex

\input 09min.tex

\input 10min.tex
\input 11min.tex

\input 12min.tex
\input 13min.tex

\input 14min.tex
\input 15min.tex

\input bibl.tex
\vfill\eject
\enddocument

\end

%% file: pmacro2.tex
\input amstex.tex
\documentstyle{amsppt}
\expandafter\xdef\csname makeatother\endcsname{\catcode`\noexpand\@=
\the\catcode`\@}
\catcode`\@=11

\def\proclaimheadfont@{\sc}
\makeatother
\magnification=\magstep1
\refstyle{A}
\nologo
\PSAMSFonts
\loadbold

\input xy
\xyoption{all}
\CompileMatrices

\def\line{\hbox to\hsize}

\font\sc=cmcsc10

\def\C{\Bbb C}
\def\Z{\Bbb Z}
\newcount\t
\newcount\q
\newcount\x
\t=0 \q=0 \x=0

\def\epsilon{\varepsilon}
\long\def\se#1{\advance\q by1
\t=0 \x=0
\csname head\endcsname {\S\number\q.\ }#1\endhead}

\long\def\sse#1{\global\advance\x by1
\csname subhead\endcsname {\number\q.\number\x.\ }{#1} \endsubhead}

\long\def\form#1{\global\advance\t by1%
$$ #1 \tag\number\q.\number\t $$}

\long\def\thm#1{\advance\x by1
\csname proclaim\endcsname {Theorem \number\q.\number\x} {#1}\endproclaim}

\long\def\lem#1{\advance\x by1
\csname proclaim\endcsname {Lemma \number\q.\number\x} {#1}\endproclaim}

\long\def\prop#1{\advance\x by1
\csname proclaim\endcsname {Proposition \number\q.\number\x} {#1}\endproclaim}

\long\def\cor#1{\advance\x by1
\csname proclaim\endcsname {Corollary \number\q.\number\x} {#1}\endproclaim}

\def\dimo{\demo{Proof}}
\def\enddimo{\qed\enddemo}

\long\def\rem#1{\advance\x by1
\csname definition\endcsname {Remark \number\q.\number\x} {#1}\enddefinition}

\long\def\exam#1{\advance\x by1
\csname definition\endcsname {Example \number\q.\number\x} {#1}\enddefinition}

\def\supp{\roman{supp}}

\def\cfrecc{\ar@{=}|*{\SelectTips{cm}{}\dir{<}}[r]}
\def\bfrecc{\ar@{=}|*{\SelectTips{cm}{}\dir{>}}[r]}

\long\def\comment#1\endcomment{}

\def\Ad{\mathop\roman{Ad}}

%% file: 01min.tex
\se{Introduction}
In the last decades the study of
$CR$ manifolds grew to
become an increasingly
important theme of research
(see e.g. \cite{AnF}, \cite{AnH}, \cite{BaERo}, \cite{HN1}, \cite{Tr}).
Particularly important was the work of 
N.Tanaka (\cite{T1}, \cite{T2}). 
He  considered the 
$CR$ manifolds as generalized contact manifolds
carrying a partial complex structure, and showed that under
some regularity and
strict nondegeneracy assumptions 
the study of their differential geometrical invariants
fits into the scheme
of Cartan geometry. Later Chern and Moser (\cite{CM})
deeply investigated the $CR$ invariants of strictly 
Levi nondegenerate
hypersurfaces.\par
More recently, there has been an increasing interest in the
study of {\it homogeneous} $CR$ manifolds, both  
of the hypersurface type and of arbitrary $CR$ 
codimension (\cite{AS1}, \cite{AS2}, \cite{AS3}, 
\cite{AzHuR}, \cite{K}, \cite{KZ}, \cite{St}). 
They provide the  natural examples that suggest and
motivate also the directions in which it is reasonable
to pursue the analysis on the general abstract $CR$ manifolds.
In fact the generalization in \cite{HN2} of the classical
notion of (local) pseudoconcavity of \cite{HN1} 
was largely motivated by the work on homogeneous models of
\cite{AlN}, \cite{LN}, \cite{MeN1},
\cite{MeN2}, \cite{MeN3}, \cite{MeN4}. 
\par
In \cite{MeN1} we
associated to the graded Lie algebras introduced by N.Tanaka
in \cite{T1} and \cite{T2},
and that we called in \cite{MeN1}
{\it Levi-Tanaka algebras},  some {\it standard} $CR$ manifolds, 
showing
in \cite{MeN4} that they are characterized by special rigidity properties.
These standard $CR$ manifolds are compact if and only if
they are minimal orbits for the action of a real form in
a complex flag manifold. 
In turn, all the minimal orbits 
of the action of real forms in
complex flag manifolds are compact
homogeneous $CR$ manifolds. However, not all of them are standard.
Some, which  are strictly 
Levi nondegenerate, correspond to non compact
standard models. Beside, there are others
(see \S 12 for the complete classification) 
which have an {\it irreducible}
$CR$ structure, 
but are not strictly Levi nondegenerate.
These considerations
lead us to investigate in \cite{MeN5} more general objects, that we
called 
{\it $CR$ algebras}.
They are canonically associated to homogeneous $CR$ manifolds,
as will be explained in \S 4. Unlike the Levi-Tanaka algebras, they
are not required to be
graded. Indeed, the
existence 
of a $CR$ compatible $\Bbb{Z}$ or $\Bbb{Z}_2$-gradation
of the $CR$ algebra
was shown in \cite{LN} to be related to that of a
Riemannian $CR$-symmetric structure (in the sense of \cite{KZ}) of 
the associated $CR$ manifold.
In particular in \cite{MeN5} we discussed
weaker nondegeneracy assumptions to ensure
finite dimensionality of the Lie algebra of infinitesimal
$CR$ isomorphism, and hence the possibility
of utilizing the homogeneous model as the $0$-curvature
object of a Cartan geometry. Among the different nondegeneracy conditions
for the partial complex structure of a $CR$ manifold $M$, we discussed
in \cite{MeN5} 
the concept of {\it weak nondegeneracy}. 
A homogeneous $CR$ manifold $M$ which is not
weakly nondegenerate is locally the product of
a $CR$ manifold and of a complex manifold of positive dimension.
A complete classification of the minimal orbits that are
weakly nondegenerate is obtained in \S 11.
\par
In this paper we concentrate on the $CR$ structure of the minimal orbit
$M$ of the action of a real form $\bold{G}$ 
in a complex flag manifold $\Hat{\bold{G}}/\bold{Q}$. These orbits,
especially the open ones,
have been studied in connection 
with representation theory 
(cf. \cite{GiMa}, \cite{Hu}, \cite{HuW}, \cite{W1}, \cite{Zi}).
Our point of view here is strictly that of $CR$ geometry.
\par\smallskip
The arguments are organized as follows. In \S 2, 3, 4
we rehearse 
the essential
definitions and notions to prepare
the general setting for the study of the $CR$ geometry of the
minimal orbits.
\par
In \S 5, 6 we give the notion of {\it parabolic $CR$ algebras} 
and we associate to each
minimal orbit $M$ a 
special parabolic $CR$ algebras that we call
{\it minimal}. Minimal orbits and 
parabolic minimal $CR$ algebras
are in a one to one correspondence. 
We classify parabolic minimal $CR$ algebras,
and thus the minimal orbits,  
by attaching to each of them a
{\it cross-marked Satake diagram}.
\par
In \S 7 we study some special morphisms of $CR$ algebras, which are 
infinitesimal analogues of smooth $\bold G$ equivariant fibrations. 
They will be an essential tool in the following sections.
\par
In \S 8 we compute the fundamental group of $M$ and
we show that, under a condition (F) that
is shared by all minimal orbits that are {\it fundamental}, (i.e.
those in which the Cauchy-Riemann distribution generates the
full tangent space), 
all $\bold{G}$-homogeneous $CR$ manifolds that are locally
$CR$-diffeomorphic to $M$ are simply connected and globally 
$CR$-diffeomorphic to $M$.
\par
In \S 9 we read off the cross-marked
Satake diagrams the property of being
fundamental and prove that every $\bold{G}$-homogeneous $CR$ manifolds 
that is locally $CR$-diffeomorphic to $M$ admits a 
{\it fundamental reduction} which is a $CR$ fibration on a totally
real basis with a connected and simply connected fiber.
\par
In \S 10 we characterize totally real and totally complex minimal orbits, 
and in \S 11, 12 we read off the cross-marked Satake diagram
the property of weak and strong nondegeneracy.
\par
In \S 13 we classify all minimal orbits that are essentially
pseudoconcave, justifying the experimental claim made in
\cite{HN2} that 
``{\sl the vast majority of them are essentially pseudoconcave}''.
\par
Finally, in \S 14, we study the space of global smooth $CR$ functions
on $M$.
\par\smallskip
A Satake diagram gives a 
graphic representation of 
the conjugation
defined by 
a real form 
on the Dynkin diagram of
the corresponding complex simple Lie algebra.
We largely utilize Satake diagrams 
in the presentation of our
results. Thus we found expedient,
to fix the notation and for identifying 
specific set of simple roots used to define special parabolic
subalgebras,
to add, at the
end of the paper, the table of 
the Satake diagrams of all the non compact real forms
of simple real Lie algebras of the real type. In all statements
concerning cross-marked Satake diagrams, we understand that the
notation refers to that table.

%% file: 02min.tex
\se{Preliminaries on $CR$ manifolds}\edef\sectpreliminaries{\number\q}
\edef\sectionb{\number\q}
We briefly rehearse some basic notions for $CR$ manifolds
(see e.g. \cite{AnF}, \cite{HN1}, \cite{MeN1}, \cite{BaERo}).\par
An (abstract) almost $CR$ manifold of type $(n,k)$ is
a triple $(M,HM,J)$,
consisting of a paracompact smooth manifold $M$
of real dimension $(2n+k)$, of a smooth subbundle $HM$ of $TM$ of 
even rank $2n$, its {\it holomorphic tangent space},
and of a smooth partial complex structure 
$J:HM@>>>HM$, $J^2=-\pmb{1}$, on the fibers of $HM$.
The integer $n\geq 0$ is the {\it $CR$ dimension} and
$k$ the {\it $CR$ codimension} of $(M,HM,J)$.
\par Let $T^{1,0}M$ and 
$T^{0,1}M$ be the complex subbundles of the complexification $\Bbb{C}HM$
of $HM$, which correspond to the $i$- and 
$(-i)$-eigenspaces of $J$:
\form{
T^{1,0}M\,=\,\{X\,-\, i\,J\,X\, \big{|}\, X\in HM\}\,
\; ,\;
 T^{0,1}M\,=\,\{X\,+\, i\,J\,X\, \big{|}\, X\in HM\}\,.
} 
We say
that $(M,HM,J)$ is a $CR$ manifold if the formal integrability
condition 
\form{\left[ \Cal C^\infty(M,T^{0,1}M),\,
\Cal C^\infty(M,T^{0,1}M)\right]\,
\subset\, \Cal C^\infty(M,T^{0,1}M)
\, }\edef\intcondit{\number\q.\number\t}
\noindent
\!\!{holds} [we get an equivalent condition by
substituting $T^{1,0}$ for $T^{0,1}$ in (\intcondit)].
When $k=0$, we have $HM=TM$ and,
via the Newlander-Nirenberg theorem, we recover the
definition of a complex manifold. A smooth real manifold 
of real dimension $k$ can
always be considered as a {\it totally real} $CR$ manifold, i.e.
a $CR$ manifold of 
$CR$ dimension $0$ and $CR$ codimension $k$.
\par
\smallskip
Let $(M_1,HM_1,J_1)$, $(M_2,HM_2,J_2)$ be two abstract 
smooth $CR$ manifolds.
A smooth map $f:M_1@>>>M_2$,
with differential $f_*:TM_1@>>>TM_2$,
is a {\it $CR$ map} if $f_*(HM_1)\subset HM_2$,
and $f_*(J_1v)=J_2f_*(v)$ for every $v\in HM_1$.
We say that $f$ is a {\it $CR$ diffeomorphism} if
$f:M_1@>>>M_2$ is a smooth diffeomorphism and
both $f$ and $f^{-1}$ are $CR$ maps.\par
A {\it $CR$ function} is a 
$CR$ map $f:M@>>>\Bbb{C}$
of a $CR$ manifold $(M,HM,J)$ in $\Bbb{C}$, endowed with the standard
complex structure. 
\par
Let $M$ be a $CR$ manifold.
Denote by ${T^*}^{1,0}M$ the
annihilator of $T^{0,1}M$ in the complexified cotangent bundle 
$\Bbb{C}T^*M$ and by $Q^{0,1}M$ 
the quotient bundle 
\hfill\linebreak
$\Bbb{C} T^* M /{T^*}^{1,0}M$, 
with projection $\pi_Q$.  It is a rank $n$ complex vector
bundle on $M$, dual to $T^{0,1}M$.  The $\bar\partial_M$-operator
acts on smooth complex valued 
functions by\,: $\bar\partial_M= \pi_Q\circ d$.
The $CR$ functions on $M$ are the smooth
solutions $u$ of $\bar\partial_Mu=0$, i.e. of
${L}u=0$ for all ${L}\in T^{0,1}M$.
We shall denote by $\Cal{O}_M(M)$ the space of smooth
$CR$ functions on $M$.\par
A $CR$ {\it embedding}
$\phi$ of an abstract $CR$ manifold $(M,HM,J)$ into
a complex manifold $\frak{M}$, with complex structure
$J_{\frak{M}}$, is a $CR$ map which is a smooth embedding
and satisfies
$\phi_*(H_pM)=\phi_*(T_pM)\cap J_{\frak{M}}(\phi_*(T_pM))$
for every $p\in M$. We say
that the embedding is {\it generic} if the complex dimension of 
$\frak{M}$
is $(n+k)$, where $(n,k)$ is the type of $M$. A real analytic
$CR$ manifold $(M,HM,J)$ always admits an embedding into a
complex manifold $\frak{M}$ (see \cite{AnF}).
\par
If $\phi:M@>>>\frak{M}$ is a smooth embedding of a paracompact smooth
manifold $M$ into a complex manifold $\frak{M}$, for each point $p\in M$
we can define $H_pM$ to be the set of tangent vectors $v\in T_pM$
such that $J_{\frak{M}}\phi_*(v)\in\phi_*(T_pM)$. For $v\in H_pM$, let
$J_Mv$ be the unique tangent vector $w\in H_pM$ satisfying
$\phi_*(w)=J_{\frak{M}}\phi_*(v)$. If
the dimension of $H_pM$ is constant for $p\in M$, then 
$HM=\overset{.}\to{\bigcup}_{p\in M}{H_pM}$ and $J_M$
are smooth and define 
the unique $CR$ manifold structure on $M$ for which
$(M,HM,J_M)$ is a $CR$ manifold and $\phi:M@>>>\frak{M}$ 
a $CR$ embedding.
\smallskip
The {\it characteristic bundle} $H^0M$ 
of $(M,HM,J)$ is defined to be the
annihilator of $HM$ in $T^*M$. It parametrizes
the Levi form: recall that the {\it Levi form}
of $M$ at $p$ is defined for $\xi\in H^0_pM$ and $X\in H_pM$ by
\form{\Cal L(\xi;X)\, = \, d\tilde\xi(X,JX)=\langle\xi,
[J\tilde X,\tilde X]\rangle\, ,}
where $\tilde\xi\in\Cal C^\infty(M,H^0M)$ and $\tilde
X\in\Cal C^\infty(M,HM)$
are smooth extensions of $\xi$ and $X$.
For each fixed $\xi$ it is a Hermitian quadratic form for the
complex structure $J_p$ on $H_pM$. \par
The map
$HM\ni X@>>>{\ssize\frac{1}{2}}(X-iJX)\in T^{1,0}M$ yields
for each $p\in M$ an $\Bbb{R}$-linear isomorphism
of $H_pM$ with 
the complex linear space $T^{1,0}_pM$, 
in such a way that the antiinvolution $J_p$ 
on $H_pM$ becomes in $T^{1,0}_pM$
the multiplication by the imaginary unit $i$. In this way we 
associate to the Levi form $\Cal{L}(\xi; \,\cdot\,)$ a unique
Hermitian symmetric form 
$T^{1,0}_p\times T^{1,0}_p\ni (Z_1,Z_2) @>>> \Cal{L}_{\xi}(Z_1,Z_2)
\in\Bbb{C}$
such that $\Cal{L}(\xi;X)=\frac{1}{4}\Cal{L}_{\xi}(X-iJX,X-iJX)$
for all $X\in H_pM$.\smallskip
In the next sections, to shorten notation, we shall write simply 
$M$, or $M^{n,k}$, for a $CR$ manifold $(M,HM,J)$ of type
$(n,k)$, as the $CR$ structure will in general be clear from the context.
\medskip
We recall that a $CR$ manifold $M$ is:
\roster
\item"$\bullet$" of {\it finite kind\/} (or {\it finite type\/})
at $p\in M$ if the higher order 
commutators of $C^\infty(M,HM)$, evaluated at $p$, 
span $T_pM$;
\item"$\bullet$" {\it strictly nondegenerate\/} (or {\it Levi nondegenerate\/})
at $p\in M$ if for each 
$\bar Z\in C^\infty(M,T^{0,1}M)$,
with $\bar Z(p)\neq 0$, there exists
$Z'\in C^\infty(M,T^{1,0}M)$ such that:
$$
[Z',\bar Z](p)\not\in T^{1,0}_pM+T^{0,1}_pM;
$$
This is equivalent to the following\,:
for every $Z\in T^{1,0}_pM\setminus\{0\}$ there exist $\xi\in H^0_pM$
and $Z'\in T^{1,0}_pM$ such that $\Cal{L}_{\xi}(Z,Z')\neq 0$.
\item"$\bullet$" {\it weakly nondegenerate\/} at $p\in M$ if for each 
$\bar Z\in C^\infty(M,T^{0,1}M)$,
with $\bar Z(p)\neq 0$, there exist $m\in\Bbb N$ and
$Z_1,\dots,Z_m\in C^\infty(M,T^{1,0}M)$ such that:
$$
[Z_1,\dots,Z_m,\bar Z](p)=[Z_1,[Z_2,\dots,[Z_m,\bar Z]\dots]](p)
\not\in T^{1,0}_pM+T^{0,1}_pM.
$$\endroster

%% file: 03min.tex
\se{The minimal orbit in a complex flag manifold}\edef\paragflag{\number\q}
Throughout this paper, we shall consistently use the symbol
$\Hat{V}$ to indicate the {\it complexification} of a real
vector space $V$. 
\medskip
A {\it complex flag manifold} is a coset space
$\frak{M}=\Hat{\bold{G}}/\bold{Q}$, where 
$\Hat{\bold{G}}$ is a connected complex semisimple Lie group 
and $\bold{Q}$ is parabolic in $\Hat{\bold{G}}$.
The manifold $\frak{M}$ 
is a closed complex projective variety which only depends on the
Lie algebras $\Hat{\frak{g}}$ of $\Hat{\bold{G}}$ and
$\frak{q}$ of $\bold{Q}$\,: 
this is a consequence of the fact that the center of a
connected and simply connected complex Lie group 
is contained in each of its parabolic subgroups. 
\par\smallskip
A {\it real form} of $\Hat{\bold{G}}$ is a real subgroup 
$\bold G$ of $\Hat{\bold{G}}$ whose 
Lie algebra $\frak g$ is a real form of $\hat\frak g$
(i.e.\ $\hat\frak g=\C\otimes_{\Bbb R} \frak g$). 
The
real form $\frak{g}$ is the set of fixed points of
an anti-involution $\sigma$ in $\Hat{\frak g}$\,:  
$\frak g=
\roman{Fix}_{\hat\frak g}(\sigma)=\{X\in\hat\frak g\,|\,\sigma(X)=X\}$. 
\par
A 
real form $\bold G$ acts on the complex flag manifold $\frak{M}$ 
by left multiplication, and $\frak{M}$ 
decomposes into
a disjoint union of $\bold G$-orbits. 
In \cite{W1} it is shown that there are  finitely many orbits, and
a unique one which is closed (hence compact).
This orbit 
$M$ has
minimal dimension and is connected. In particular, 
the connected component of the identity $\bold{G}^{\circ}$
of $\bold{G}$ is transitive on $M$. 
Thus, while studying $M$, we can as well assume
that $\bold{G}=\bold{G}^{\circ}$ is connected.
\par
Moreover, up to conjugation, we can arrange  that the 
closed orbit is
$M=\bold G\cdot o$, where $o=e\bold Q$. 
We shall denote by $\bold G_+=\bold G\cap\bold Q$ 
the isotropy
subgroup of $\bold G$ at $o$ and 
by $\frak g_+=\frak{g}\cap\frak q$ its Lie 
algebra.
\par
\smallskip
The closed orbit $M$ has a $CR$ structure 
induced by its
embedding in the complex manifold $\frak{M}$.
This can also be described by using the canonical identifications\,:
$T_o\frak{M}\simeq\hat\frak g/\frak q$ 
and $T_oM\simeq\frak g/\frak g_+\subset T_o\frak{M}$.
We have then\,:
\form{
H_oM\simeq \big(\frak g\cap(\frak q+\bar\frak q)\big)/\frak g_+
}
and
\form{\cases
J_o:H_oM @>>> H_oM\quad\text{is defined by\,:}\\
J_o(X+\bar X+\frak g_+)=iX-i\bar X+\frak g_+
\quad \text{for $X\in\frak q$\,.}
\endcases}
The bundle  
$HM$ and 
the partial complex
structure $J$ 
are defined, at all points of $M$, 
in such a way that $\bold G$ acts on $M$ as a group of 
$CR$ automorphisms. By using the identification
$T_{g\bold{Q}}\frak{M}=\Hat{\frak{g}}/(\Ad(g)\frak{q})$, 
we obtain\,: 
\form{\cases
T_{g\bold G_+}M\simeq\frak g/(\Ad(g)\frak g_+)\\
H_{g\bold G_+}M=\big(\frak g\cap(\Ad(g)\frak q
+\Ad(g)\bar\frak q)\big)/\Ad(g)\frak g_+\\
J_{g\bold G_+}(X+\bar X+\Ad(g)\frak g_+)=iX-i\bar X+\Ad(g)\frak g_+\quad
\text{for $X\in\Ad(g)\frak q$.}
\endcases}
Note that the embedding of $M$ into $\frak{M}$ is always generic,
because $T_p\frak{M}=T_pM+\overline{T_pM}$ at every $p\in M$.

%% file: 04min.tex
\se{Homogeneous $CR$ manifolds and $CR$ algebras}
\edef\secta{{\number\q}}
To a $CR$ manifold $M$, which is homogeneous for the action of a
real Lie group $\bold G$ of $CR$ transformations, we associate a
{\it $CR$ algebra} $(\frak g,\frak q)$.
This is a pair consisting of the real Lie algebra $\frak g$ of 
the group $\bold G$ and of a complex subalgebra $\frak q$ of its
complexification $\hat\frak g$.
This subalgebra $\frak q$ is the inverse image of $T_o^{0,1}M$
by the complexification 
$\Hat{\pi}_*$ of the differential $\pi_*:\frak g\simeq 
T_e\bold G @>>>T_oM$ of the group action at $e$.
Note that the fact that $\frak q$ is a complex Lie subalgebra of 
$\hat\frak g$ is a consequence of the formal integrability condition
(\intcondit) for $CR$ manifolds.
\par
Let $(\frak{g},\frak{q})$ be a $CR$ algebra.
The real Lie subalgebra $\frak g_+=\frak g\cap\frak q$ is 
called its {\it 
isotropy\/}.
Let $\bold{G}$ be a connected and simply connected Lie group with
Lie algebra $\frak{g}$. 
Assume that the analytic subgroup $\bold{G}_+$
of $\bold{G}$ corresponding to the Lie subalgebra $\frak{g}_+$
is closed in $\bold{G}$.
Then the homogeneous space $M=\bold{G}/\bold{G}_+$ is
a smooth paracompact manifold and
has a unique $CR$ structure such that\,: 
\roster
\item"$\bullet$" $T^{0,1}_oM=\Hat{\pi}_*(\frak{q})$
\item"$\bullet$"
$\bold{G}$ acts on $M$ by $CR$ diffeomorphisms.
\endroster
We denote the $CR$ manifold
$\bold{G}/\bold{G}_+$
by $\tilde{M}(\frak{g},\frak{q})$.  
\par
For general definitions and basic properties of $CR$ algebras we refer 
the reader to \cite{MeN5}. In particular, we recall that a
{\it morphism of $CR$ algebras} 
$\phi:(\frak{g},\frak{q})@>>>(\frak{g}',\frak{q}')$ is
a homomorphism $\phi:\frak{g}@>>>\frak{g}'$ of real Lie algebras
whose complexification $\Hat{\phi}$ satisfies
$\Hat{\phi}(\frak{q})\subset\frak{q}'$; it is a {\it $CR$ submersion}
if $\phi(\frak{g})+\frak{g}_+'=\frak{g}'$ and
$\Hat{\phi}(\frak q)+\frak{q}'\cap\bar{\frak{q}}'=\frak{q}'$. 
\par\smallskip
We say that the $CR$ algebra $(\frak g,\frak q)$ is:
\roster
\item"$\bullet$" {\it effective\/} if there are
no ideals of $\frak g$ contained in $\frak g_+$;
\item"$\bullet$"  {\it fundamental\/} if $\frak q+\bar\frak q$ 
generates $\Hat{\frak g}$\,;
\item"$\bullet$"  {\it ideal nondegenerate\/} if there is no ideal 
$\frak a$ of $\frak g$ with 
$\frak a\subset\frak q+\bar\frak q$ 
and $\frak a\not\subset\frak g_+$;
\item"$\bullet$"  {\it weakly nondegenerate\/} if there are no complex 
subalgebras $\frak q'\subset\hat\frak g$ \; with\; \par\centerline{
$\frak q\subsetneq\frak q'\subset\frak q+\bar\frak q$\,;}\par 
\item"$\bullet$"  {\it strictly nondegenerate\/} if for every 
$Z\in\frak q\setminus\bar\frak q$ 
there exists $Z'\in\bar\frak q$ such that
$[Z,Z']\notin\frak q+\bar\frak q$\,.
\endroster
If $(\frak g,\frak q)$ is the $CR$ algebra associated to a 
$\bold G$-homogeneous $CR$ manifold $M$, these notions 
express geometric 
properties of 
$M$ (see \cite{MeN5})\,:
effectiveness is equivalent to almost effectiveness
(i.e.\ discreteness of the isotropy subgroup) of the $\bold G$ action;
fundamental to finite kind; ideal nondegeneracy to
holomorphic nondegeneracy (see \cite{BaERo}); weak and strict nondegeneracy
to weak and strict nondegeneracy as defined at the end of 
\S{\sectpreliminaries}.
\par

When $(\frak{g},\frak{q})$ is weakly degenerate, it was proved
in \cite{MeN5} that there is a $CR$-fibration 
$\tilde{M}(\frak{g},\frak{q})@>>>M'$ of the corresponding homogeneous
$CR$ manifold $\tilde{M}(\frak{g},\frak{q})$ 
on a $CR$ manifold $M'$ with
the same $CR$ codimension, having a non trivial complex fiber.
For homogeneous 
simply connected $CR$ manifolds, the condition of
weak degeneracy of \S{\paragflag} 
is in fact necessary and sufficient
for the existence of $CR$ fibrations with non trivial complex
fibers. Indeed, for general $CR$ manifolds, the existence
of a $CR$ fibration with non trivial complex fibers 
implies weak degeneracy, as we have\,:
\prop{Let $M$ and $M'$ be $CR$ manifolds. We assume that
$M'$ is locally embeddable and that
there exists a $CR$ fibration $M@>\pi>> M'$ with totally complex
fibers of positive dimension. 
Then $M$ is weakly degenerate.}\edef\proposizioneda{\number\q.\number\x}
\dimo
Let $f$ be any smooth $CR$ function defined on a neighborhood $U'$ of
$p'\in M'$. Then $\pi^*f$ is a $CR$ function in $U=\pi^{-1}(U')$,
that is constant along the fibers of $\pi$. Then, if 
$L\in\Cal C^{\infty}(M,T^{1,0}M)$ is tangent to the fibers of $\pi$
in $U$, we obtain that 
$[\bar Z_1\,,\,\,\hdots\, ,\,\,\bar Z_m,\,\,L]\left( \pi^*f\right)=0$ 
for every choice of 
$\bar Z_1,\,\hdots ,\,\, \bar Z_m\in\Cal C^\infty(M,T^{0,1}M)$.
Assume by contradiction that $M$ is weakly nondegenerate
at some $p$ with $\pi(p)=p'$. Then for
some choice of  
$\bar Z_1,\,\hdots ,\,\, \bar Z_m\in\Cal C^\infty(M,T^{0,1}M)$
we would have 
$v_{p}=[\bar Z_1\,,\,\,\hdots\, ,\,\,\bar Z_m,\,\,L]
\notin T_p^{1,0}M\oplus
T_p^{0,1}M$. 
Since the fibers of $\pi$ are totally complex,
$\pi_*(v_{p})\neq 0$. By the assumption that $M'$ is locally embeddable
at $p$, the real parts of the 
(locally defined) $CR$ functions give local coordinates
in $M'$ and therefore
there is a $CR$ function $f$ defined on a neighborhood $U'$ 
of $p'$ with $v_p(\pi^*f)=\pi_*(v_{p})(f)\neq 0$. This gives a contradiction,
proving our statement.
\enddimo
We can always reduce to the case of an almost
effective action of $\bold{G}$\,: at the level of $CR$ algebras,
this corresponds to substituting to $(\frak g,\frak q)$ its
{\it effective quotient}, which is the 
$CR$ algebra $(\frak g/\frak a,\frak q/\hat\frak a)$, where $\frak a$ is 
the maximal ideal of $\frak g$ that is contained in
$\frak{g}_+$, and $\Hat{\frak{a}}$ its
complexification in $\Hat{\frak{g}}$ (see \cite{MeN5, Lemma 4.7}).

%% file: 05min.tex
\se{Parabolic $CR$ algebras}
In the following, we shall restrict our consideration to the the case of
{\it parabolic\/} $CR$ algebras, i.e.\ those $CR$ algebras 
$(\frak g,\frak q)$
where $\frak g$ is finite dimensional and $\frak q$ is a parabolic 
subalgebra of $\hat\frak g$. In this section we explain some of their
simplest properties.
\par
\prop{A parabolic $CR$ algebra $(\frak g,\frak q)$ is effective if
and only if the following two conditions are satisfied:
\roster
\item"$(i)$" {$\frak g$ is semisimple},
\item"$(ii)$" {no simple ideal of $\hat\frak g$ is 
contained in $\frak q\cap\overline{\frak q}$.}
\endroster}
\dimo The statement follows by observing that: (a) for a parabolic
$(\frak g,\frak q)$ the radical $\frak r$ of $\frak g$ is contained
in $\frak g_+$; (b) if an ideal $\frak a$ of $\hat\frak g$ is
contained in $\frak q\cap\overline{\frak q}$, then $\frak a+\overline{
\frak a}$ is the complexification of an ideal $\frak b$ of
$\frak g$ contained in $\frak g_+$.
\enddimo
\par
To an effective parabolic $CR$ algebra $(\frak g,\frak q)$ we
associate a $CR$ manifold $M=M(\frak g,\frak q)$, unique
modulo isomorphisms,
defined as the orbit
$\bold G\cdot o$ in $\hat\bold G/\bold Q$, where:
\roster
\item"$\bullet$" $\hat\bold G$ is a connected and simply connected Lie 
group with Lie algebra $\hat\frak g$\,;
\item"$\bullet$" $\bold Q=\bold N_{\hat\bold G}(\frak q)$ is the 
parabolic subgroup of $\hat\bold G$ with Lie algebra $\frak q$\,;
\item"$\bullet$" $\bold G$ is the analytic real subgroup of $\hat\bold G$
with Lie algebra $\frak g$.
\endroster
Note that for a parabolic $(\frak{g},\frak{q})$ also
 $\tilde{M}(\frak{g},\frak{q})$ is well defined and 
is the universal cover of ${M}(\frak{g},\frak{q})$.
\par
The following proposition reduces the study of effective parabolic
$CR$ algebras $(\frak g,\frak q)$ to the case where $\frak g$ is a
{\it simple} real Lie algebra.
\prop{Let $(\frak g,\frak q)$ be an effective parabolic $CR$ algebra
and let $\frak g =\frak g_1\oplus \cdots \oplus\frak g_\ell$ be the
decomposition of $\frak g$ into the direct sum of its simple ideals.
Then:
\roster
\item"($i$)" $\frak q=\frak q_1\oplus \cdots \oplus\frak q_\ell$
where $\frak q_j=\frak q\cap\hat\frak g_j$ for $j=1,\hdots,\ell$;
\item"($ii$)" for each $j=1,\hdots,\ell$, $(\frak g_j,\frak q_j)$
is an effective parabolic $CR$ algebra;
\item"($iii$)" $(\frak g,\frak q)$ is ideal (resp. weakly, strictly)
nondegenerate if and only if for each $j=1,\hdots,\ell$, the
$CR$ algebra $ (\frak g_j,\frak q_j)$ is  ideal (resp. weakly, strictly)
nondegenerate;
\item"($iv$)" $(\frak g,\frak q)$ is fundamental
if and only if for each $j=1,\hdots,\ell$, the
$CR$ algebra $ (\frak g_j,\frak q_j)$ is fundamental.
\endroster}\edef\teoremaeb{\number\q.\number\x}
\dimo In fact $\Hat{\frak{g}}=\bigoplus_{j=1}^\ell{\Hat{\frak{g}}_j}$
is a decomposition of $\Hat{\frak{g}}$ into a direct sum of ideals.
The decomposition ($i$)
of $\frak q$ follows then from the decomposition
$\Hat{\frak{h}}=\bigoplus_{j=1}^{\ell}{\left(
\Hat{\frak{h}}\cap\Hat{\frak{g}}_j\right)}$ of any Cartan subalgebra
of $\Hat{\frak{g}}$ contained in $\frak q$ (see \cite{B2, Ch.VII,\S 2,Prop.2}).
\par
The proof of the
other statements is straightforward.
\enddimo
We note that, with the notation of Proposition \number\q.\number\x,\par
\smallskip
\centerline{
$M(\frak g,\frak q)\simeq M(\frak g_1,\frak q_1) \times\cdots\times
M(\frak g_\ell,\frak q_\ell)$}\par
\smallskip
\centerline{
$\tilde{M}(\frak g,\frak q)\simeq \tilde{M}(\frak g_1,\frak q_1) 
\times\cdots\times
\tilde{M}(\frak g_\ell,\frak q_\ell)$}\par\smallskip\noindent
where "$\simeq$" means isomorphism of $CR$ manifolds.\par\smallskip
When $\frak q$ is parabolic in $\hat\frak g$,
 its conjugate $\overline{\frak q}$
with respect to the real form $\frak g$ is also parabolic in
$\hat\frak g$. Therefore the intersection $\frak q\cap\overline{\frak q}$
contains a Cartan subalgebra $\hat{\frak h}$ that is invariant under
conjugation (see e.g.\ \cite{B2, Ch.VIII, \S 3, Prop.10}. We observe
that $\frak{g}_+$ contains an element $A$ 
with $\roman{ad}_{\frak q\cap\bar\frak q}(A)$ 
semisimple and of maximal rank. The centralizer $\Hat{\frak{h}}$
of $A$ in $\Hat{\frak{g}}$ is a Cartan subalgebra of $\Hat{\frak{g}}$
that is contained in $\frak q\cap\bar\frak q$ and is invariant under
conjugation). 
The intersection $\frak h=\hat{\frak h}\cap\frak g$ is
a Cartan subalgebra of $\frak g$. Thus we have\,:
\lem{Let $(\frak g,\frak q)$ be a parabolic $CR$ algebra.
Then $\frak g_+=\frak q\cap\frak g$ contains a Cartan subalgebra of
$\frak g$.\qed}
Moreover we have:
\prop{Let $(\frak g,\frak q)$ be an effective parabolic $CR$ algebra.
The set $\frak n_+$ 
of the elements $A$ of the radical
$\frak r(\frak g_+)$ of $\frak g_+$, for which
$\roman{ad}_{\frak g}(A):\frak g@>>>\frak g$ is nilpotent, 
is a nilpotent ideal of $\frak g_+$ and there exists
a reductive subalgebra $\frak w$ of $\frak g_+$ such that
\form{\frak g_+=\frak n_+\oplus\frak w\, .}
The reductive subalgebra $\frak w$ is uniquely determined modulo
the subgroup of inner automorphisms of $\frak g_+$ generated
by those of the form $\exp\left(\roman{ad}_{\frak g_+}(X)\right)$
with $X\in\frak n_+$.  }
\dimo
Indeed $\frak q$, being parabolic, contains the semisimple and
nilpotent parts of its elements. If $X\in\frak q$ belongs to the
real form $\frak g$, then also its semisimple and nilpotent parts
belong to $\frak g$. Therefore $\frak g_+$ is splittable, i.e. contains
the semisimple and nilpotent part of its elements and we can apply
\cite{B2, Prop.7, \S 5, Ch.VII} to obtain our statement.\enddimo
\medskip
Let $\frak z$ denote the center of $\frak w$ and let 
$\frak s=[\frak w,\frak w]$ be its semisimple ideal. Then
\form{\frak w\, = \, \frak z\oplus\frak s\, .}
Thus, a Cartan subalgebra $\frak h\subset\frak g_+$ of $\frak g$ 
can be taken as the direct sum 
\form{\frak h \, = \, \frak z\oplus\frak h'}
\noindent
of the center\edef\cartdec{\number\q.\number\t} 
$\frak z$ of $\frak w$ and a Cartan subalgebra $\frak h'$
of $\frak s$. Vice versa, every Cartan subalgebra $\frak h$ of $\frak g$
contained in $\frak g_+$ has the form (\cartdec) for some reductive
subalgebra $\frak w$ of $\frak g_+$.
\medskip
Choose a Cartan subalgebra $\frak h$ of $\frak g$, and denote by
$\Hat{\frak{h}}$ its complexification, which 
is a Cartan subalgebra of 
$\Hat{\frak{g}}$. Fix a Cartan decomposition 
\form{\frak{g}\, = \, \frak{k}
\oplus\frak{p}}\edef\formulacartan{\number\q.\number\t}
of $\frak{g}$, where $\frak{k}$ is a maximal compact subalgebra
of $\frak{g}$ and $\frak{p}$ its orthogonal in $\frak{g}$
with respect to the Killing form, such that\,:
\form{\frak{h}\, = \, \frak{h}^+\oplus\frak{h}^-\qquad
\roman{where}\quad \frak{h}^+=\frak{h}\cap\frak{k}\quad\roman{and}\quad
\frak{h}^-=\frak{h}\cap\frak{p}\,;}\edef\formuladopocartan{\number\q.\number\t}
$\frak{h}^+$ is the {\it toroidal} and 
$\frak{h}^-$ the {\it vector part} of $\frak{h}$. \par
Set $\frak{h}_{\Bbb R}=\,(i\,\frak{h}^+)\,\oplus\frak{h}^-$.
This is a real subalgebra of $\Hat{\frak{h}}$ and the Cartan subalgebra
of a real split form of $\Hat{\frak{g}}$\,;
we denote by $\frak{h}_{\Bbb R}^*$ its dual.
Let
$\Cal R=\Cal R(\Hat{\frak g},\Hat{\frak h})
\subset\frak{h}_{\Bbb R}^*$ 
be
the associated root system 
and denote by
$\Hat{\frak{g}}^{\alpha}\subset\Hat{\frak{g}}$ the eigenspace
corresponding to
the root $\alpha\in\Cal R$. \par
We shall denote by $\bold{W}_{\Hat{\frak h}}$ the Weyl group
of $\Cal{R}$\,: it is the group of isometries of 
$\frak{h}_{\Bbb{R}}^*$ generated by the reflections $s_{\alpha}$
with respect to the elements $\alpha\in\Cal{R}$.
We recall that it is canonically 
identified to the quotient 
$\bold{N}_{\bold{Int}(\Hat{\frak{g}})}(\Hat{\frak{h}})/
\bold{Z}_{\bold{Int}(\Hat{\frak{g}})}(\Hat{\frak{h}})$. 
This will be also called the {\it algebraic} Weyl group,
to distinguish it from the {\it analytic} Weyl group 
$\bold{W}_{{\frak h}}$, which is the image in $\bold{W}_{\Hat{\frak h}}$
of the
composed homomorphism:\par\smallskip\centerline{
$\bold{N}_{\bold{Int}({\frak{g}})}({\frak{h}})/
\bold{Z}_{\bold{Int}({\frak{g}})}({\frak{h}})\longrightarrow
\bold{N}_{\bold{Int}(\Hat{\frak{g}})}(\Hat{\frak{h}})/
\bold{Z}_{\bold{Int}(\Hat{\frak{g}})}(\Hat{\frak{h}})
\longrightarrow \bold{W}_{\Hat{\frak h}}$.}\par\smallskip
We also consider the group $\bold{A}_{\Hat{\frak{h}}}$ of
all the isometries of $\frak{h}_{\Bbb{R}}^*$ that transform
$\Cal{R}$ into $\Cal{R}$. The (algebraic) Weyl group 
$\bold{W}_{\Hat{\frak h}}$ is a normal subgroup of 
$\bold{A}_{\Hat{\frak{h}}}$. We have a natural isomorphism
$\bold{N}_{\bold{Aut}(\Hat{\frak{g}})}(\Hat{\frak{h}})/
\bold{Z}_{\bold{Aut}(\Hat{\frak{g}})}(\Hat{\frak{h}})
@>{\sim}>>\bold{A}_{\Hat{\frak{h}}} $, yielding a commutative diagram:
\smallskip
$${\CD
@. @. \pmb{1} @. \pmb{1}\\
@. @. @VVV @VVV \\
\pmb{1} @>>>
\bold{Z}_{\bold{Int}(\Hat{\frak{g}})}(\Hat{\frak{h}})
@>>> \bold{N}_{\bold{Int}(\Hat{\frak{g}})}(\Hat{\frak{h}})
@>>> \bold{W}_{\Hat{\frak{h}}} @>>>\pmb{1}\\
@. @| @VVV  @VVV \\
\pmb{1} @>>>
\bold{Z}_{\bold{Aut}(\Hat{\frak{g}})}(\Hat{\frak{h}})
@>>> \bold{N}_{\bold{Aut}(\Hat{\frak{g}})}(\Hat{\frak{h}})
@>>> \bold{A}_{\Hat{\frak{h}}} @>>>\pmb{1}\\
@. @. @. @VVV \\
@. @. @. \bold{Aut}(\Delta)\\
@. @. @. @VVV \\
@. @. @. \pmb{1}
\endCD}$$
\smallskip\noindent
with exact rows and columns, in which we denoted by
$ \bold{Aut}(\Delta)$ the group of automorphisms of the Dynkin
diagram $\Delta$ associated to $\Cal{R}$. 
\par
We finally define the group $\bold{A}_{\frak{h}}$ 
as the image in
$\bold{A}_{\Hat{\frak{h}}}$ of 
$\bold{N}_{\bold{Aut}({\frak{g}})}({\frak{h}})$,
identified to a subgroup of
$\bold{N}_{\bold{Aut}(\Hat{\frak{g}})}(\Hat{\frak{h}})$, 
 by the
homomorphism described above.
\bigskip
Denote by $\frak{C}(\Cal{R})$ the set of Weyl chambers associated
to the root system $\Cal R$. 
Choose a Weyl chamber $C\subset\frak{h}_{\Bbb R}$ and let
$\prec$ be the corresponding
partial order in $\frak{h}_{\Bbb R}^*$, defined by
\form{
\alpha\prec\beta\quad\Longleftrightarrow
\quad\alpha(H)<\beta(H)\qquad \forall H\in C\, .}
Let 
$\Cal R^+=\Cal R^+(C)=\{\alpha\in\Cal R\, | \, \alpha\succ 0\}$
and 
{$\Cal R^-=\Cal R^-(C)=\{\alpha\in\Cal R\, | \, \alpha\prec 0\}$}
be the set of positive 
and the set of negative roots with respect to $C$,
respectively,
and denote by 
$\Cal{B}=\Cal{B}(C)$ the set of simple roots in $\Cal{R}^+$.
 If\par\smallskip\centerline{
$
\alpha=\sum_i n_i\alpha_i,\quad \alpha_i\in\Cal B,\ n_i\in\Z,
$}\par\smallskip\noindent
we set\par\smallskip\centerline{ 
$\supp(\alpha)=\supp_C(\alpha)=\{\alpha_i\in\Cal B\,|\,n_i\neq 0\}$.}
\par\medskip
Let 
$\Phi$ be a subset of $\Cal B$. The set $\check\Phi^r$ of those
$\beta\in\Cal R^-$  for which 
$\supp(\beta)\cap\Phi=\emptyset$
is a closed system of roots 
($\beta_1,\beta_2\in\check\Phi^r$ and $\beta_1+\beta_2\in\Cal R$
$\Longrightarrow$ $\beta_1+\beta_2\in\check\Phi^r$).
Then
\form{\frak q_{\Phi}=
\hat\frak h\oplus\sum_{\beta\in{\Cal R}^+}{\Hat{\frak{g}}^{\beta}}\oplus
\sum_{\beta
\in\check\Phi^r}{\Hat{\frak{g}}^{\beta}}}\edef\formparab{\number\q.\number\t}
\noindent 
\!\!{is} a parabolic subalgebra of $\Hat{\frak g}$, 
and every parabolic
subalgebra of $\Hat{\frak g}$ that contains $\frak h$ can be
described in this way, by a suitable choice of the Weyl chamber $C$ and 
of the set of simple roots $\Phi\subset\Cal{B}$.\par
Each set\,: 
\form{\left\{\aligned
\Cal Q\;&=\Cal Q_\Phi=\Cal R^+\cup\check\Phi^r\\
\Cal Q^r&=\Cal Q^r_\Phi=\{\alpha\in\Cal Q|-\alpha\in\Cal Q\}
=\check\Phi^r\cup(-\check\Phi^r)\\
\Cal Q^n&=\Cal Q^n_\Phi=\{\alpha\in\Cal Q|-\alpha\not\in\Cal Q\}
=\Cal R^+\setminus(-\check\Phi^r)
\endaligned\right.}
is a closed system of roots. We set\,:
\form{\cases
\frak q^r=\frak q^r_\Phi=\hat\frak h\oplus\sum_{\alpha\in\Cal Q^r}
\Hat{\frak{g}}^\alpha\\\vspace{2\jot}
\frak q^n=\frak q^n_\Phi=\sum_{\alpha\in\Cal Q^n}
\hat\frak g^\alpha\, ,
\endcases}
to obtain the decomposition\,: 
\form{\frak q=\frak q^r\oplus\frak q^n\,,}
of $\frak{q}$ into the direct sum of its
nil radical $\frak q^n$
and a reductive subalgebra
$\frak q^r$. 
\smallskip
We say that a Cartan subalgebra $\frak h$ of $\frak g$ is {\it adapted}
to the parabolic effective $CR$ algebra $(\frak g,\frak q)$ if, in
the decomposition (\cartdec), the Cartan subalgebra $\frak h'$ of
$\frak s$ has {\it maximal vector part}.
\footnote{
A Cartan subalgebra of a semisimple real Lie algebra
$\frak s$ with a maximal vector part
is obtained in the following way:
Assume that $\frak s=\frak k\oplus\frak p$ is a Cartan decomposition of
$\frak s$. Take any maximal Abelian Lie subalgebra $\frak h'_a$ of $\frak s$
contained in $\frak p$ and let $\frak h'_t$ be the centralizer of
$\frak h'_a$ in $\frak k$. Then $\frak h'=\frak h'_a\oplus\frak h'_t$
is a Cartan subalgebra of $\frak s$ with maximal vector part.
All Cartan subalgebras of $\frak s$ with maximal vector part
are conjugate and can be obtained in this way from a suitable
Cartan decomposition of $\frak s$.}
\medskip
A real Lie subalgebra 
$\frak t$ of $\frak g$
is {\it triangular} if all linear maps 
$\roman{ad}_{\frak g}(X)\in\frak{gl}_{\Bbb R}(\frak g)$ with $X\in\frak t$
can be simultaneously represented by triangular matrices in a suitable
basis of $\frak g$. All maximal triangular subalgebras of $\frak g$ are
conjugate by an inner automorphism (\cite{Mo}, \S 5.4, or \cite{V}). 
A real Lie subalgebra of $\frak g$ containing a maximal triangular
subalgebra of $\frak g$ is called a $t$-{\it subalgebra}.\par
An effective parabolic $CR$ algebra $(\frak g,\frak q)$ will be
called {\it minimal} if $\frak g_+=\frak q\cap\frak g$ is
a $t$-subalgebra of $\frak g$. \par
We observe that a maximal triangular subalgebra of $\frak g$ contains
a maximal Abelian subalgebra of semisimple elements having real eigenvalues.
Hence:
\prop{An adapted Cartan subalgebra of an effective parabolic minimal
$CR$ algebra $(\frak g,\frak q)$ has maximal vector part as a
Cartan subalgebra of $\frak g$.\qed}\edef\propaf{\number\q.\number\x}
\thm{Let $\frak g$ be a semisimple real Lie algebra and $\frak q$ a parabolic
subalgebra of its complexification $\Hat{\frak g}$. Then, up to
$CR$ isomorphisms, there is a unique parabolic minimal $CR$ algebra
$(\frak g',\frak q')$ with $\frak g'$ isomorphic to $\frak g$ and
$\frak q'$ isomorphic to $\frak q$.
}\edef\teoremaag{\number\q.\number\x}
\dimo
Fix a maximal triangular subalgebra $\frak t$ of $\frak g$. Its
complexification $\hat\frak t$ is solvable and therefore 
is contained
in a maximal solvable subalgebra,
i.e.\ a Borel subalgebra, $\frak b$ of $\Hat{\frak g}$.
Modulo an inner automorphism of $\Hat{\frak g}$, we can assume that
$\frak b\subset\frak q$. The $CR$ algebra $(\frak g,\frak q)$ is
parabolic minimal. \par
Let $\frak q$, $\frak q'$ be parabolic subalgebras of $\Hat{\frak g}$
such that $\frak g_+=\frak q\cap\frak g$ and $\frak g'_+=\frak q'\cap\frak g$
are $t$-subalgebras of $\frak g$. By an inner automorphism of $\frak g$,
we can assume that $\frak g_+$ and  $\frak g'_+$ contain the same
maximal triangular subalgebra $\frak t$ of $\frak g$
and hence a same maximal Abelian subalgebra of $\frak g$ of
semisimple elements having real eigenvalues. Hence,
using another inner automorphism of $\frak g$, 
we can assume that $\frak q$ and $\frak q'$ contain the same
maximal vectorial Cartan subalgebra $\frak h$ of $\frak g$. \par
The inner automorphism of $\Hat{\frak g}$ transforming $\frak q$ into
$\frak q'$ can now be taken to be an element of the analytic Weyl group,
leaving the Cartan subalgebra $\frak h$ and hence $\frak g$ invariant.
It defines a $CR$ isomorphism between $(\frak g,\frak q)$
and $(\frak g,\frak q')$.
\enddimo
We recall that a $CR$ algebra $(\frak g,\frak q)$ is {\it totally real}
if $\frak q=\overline{\frak q}$, or, equivalently, if
$\frak g_+=\frak q\cap\frak g=\Cal H_+=\frak g\cap\left(\frak
q+\overline{\frak q}
\right)$. This is equivalent to the fact that $M(\frak{g},\frak{q})$
is totally real, i.e. a $CR$ manifold with $CR$ dimension $0$.
For a totally real effective parabolic $CR$ algebra
$(\frak g,\frak q)$ the real subalgebra $\frak g_+$ of $\frak g $
is parabolic, hence a $t$-subalgebra of $\frak g$.
Thus we have\,:
\prop{A totally real effective parabolic $CR$ algebra is minimal.\qed}
Effective parabolic minimal $CR$ algebras correspond to minimal orbits.
In fact we have\,:
\thm{The $CR$ manifold 
$M(\frak{g},\frak{q})$, 
associated to an effective parabolic subalgebra
$(\frak g,\frak q)$, 
is compact if and only if $(\frak g, \frak q)$ is minimal.}
\dimo Indeed, since $\bold{G}$ is a linear group, 
a $\bold G$-homogeneous space  $\bold G/\bold G_+$ is compact if and only if
$\bold G_+$ contains a maximal connected triangular subgroup (see 
\cite{O, II, Ch.5, \S1.1}), i.e.\ if $\frak g_+$ is a $t$-subalgebra of 
$\frak g$.
\enddimo

%% file: 06min.tex
\se{Parabolic minimal $CR$ algebras and cross-marked Satake diagrams}
Denote by $\sigma:\frak{h}_{\Bbb{R}}^*@>>>\frak{h}_{\Bbb{R}}^*$
the involution induced by the conjugation defined by the real form
$\frak{g}$ of $\Hat{\frak{g}}$.
If $\vartheta$ is the complexification of the
Cartan involution associated to
the decomposition (\formulacartan), then the conjugation
equals $(-\vartheta)$ on $\frak{h}_{\Bbb{R}}$, so that
$\sigma=\,-\,\left(\vartheta\left|_{\frak{h}_{\Bbb{R}}}\right.\right)^*$.
\par
We note that $\sigma(\Cal R)=\Cal R$. A root $\alpha\in\Cal R$ is
called {\it real} if $\bar\alpha=\sigma(\alpha)=\alpha$,
{\it imaginary} if $\bar\alpha=\sigma(\alpha)=-\alpha$.
We shall denote by $\Cal R_{\bullet}$ the set of {\it imaginary} roots
in $\Cal R$. 
We recall from \cite{Ar}\,:
\prop{Let $\sigma:\frak{h}_{\Bbb{R}}^*@>>>\frak{h}_{\Bbb{R}}^*$
be the involution associated to the conjugation induced by the real
form $\frak g$ of $\Hat{\frak g}$. 
The real Cartan subalgebra
$\frak h$ of $\frak g$ has maximal vector part if and only if
$\Hat{\frak{g}}^{\alpha}\subset\hat\frak{k}
=\Bbb{C}\otimes_{\Bbb{R}}\frak{k}$
for all $\alpha\in\Cal{R}_\bullet$.\par
Let $\frak{h}$ be a Cartan subalgebra of $\frak{g}$ with maximal
vector part and $\Cal{R}=\Cal{R}(\Hat{\frak{g}},\Hat{\frak{h}})$.
Then there exists a Weyl chamber
$C\in\frak{C}(\Cal R)$ such that:
\roster
\item"($i$)" $\bar\alpha=\sigma(\alpha)\succ 0$ for all 
$\alpha\in\Cal{R}^+(C)\setminus\Cal{R}_\bullet$;
\item"($ii$)" there are pairwise orthogonal roots
$\beta_1,\,\hdots,\,\beta_m\in\Cal R_{\bullet}$
such that $s_{\beta_1}\circ \cdots \circ s_{\beta_m}$ 
is the element $w_{(C,\bar C)}$ 
of the Weyl group that transforms
$C$ into $\bar{C}$; in particular $w_{(C,\bar C)}$ is an involution\,:
$w_{(C,\bar C)}^2=\pmb{1}$\,;
\item"($iii$)" there is an involution 
$\varepsilon_C\in\bold{A}_{\Hat{\frak{h}}}$, such that
$\varepsilon_C(C)=C$, that commutes with $\sigma$ and with
$w_{(C,\bar C)}$, such that\,:\par\smallskip
\centerline{
$\sigma=\varepsilon_C\circ w_{(C,\bar C)}$\,.}
\endroster
The Weyl chamber $C$ is uniquely determined
modulo the analytic Weyl group $\bold{W}_{\frak h}$.
\qed}\edef\teoremaba{\number\q.\number\x}
A Weyl chamber $C$ that satisfies conditions ($i$),
($ii$) and ($iii$) of  Proposition {\teoremaba} will be called
{\it adapted to the conjugation $\sigma$.} \smallskip
We shall denote by 
$\Cal{B}_{\bullet}=\Cal{B}_{\bullet}(C)$ the
set $\Cal{B}(C)\cap\Cal{R}_{\bullet}$
of simple purely imaginary roots of $\Cal{R}^+(C)$ and by
$\Xi=\Xi(C)$ its complement in $\Cal{B}$.
The involution $\varepsilon_C$ transforms 
$\Xi$ into itself.
Moreover, 
from Proposition {\teoremaba}
we obtain the conjugation formula\,:
for every $\alpha\in\Xi$ 
there are integers $n_{\alpha,\beta}\geq 0$ such that
\form{\bar\alpha\, = \, \varepsilon_C(\alpha)\, +  \,
\dsize\sum_{\beta\in\Cal{B}_{\bullet}}{n_{\alpha,\beta}
\beta}\, .}
\edef\formulaej{\number\q.\number\t}
We associate to $C$ the {\it Satake diagram} of $\frak g$.
It is obtained from the Dynkin diagram 
of $\Hat{\frak g}$ whose nodes correspond
to the roots in $\Cal B(C)$ by painting black those 
corresponding to imaginary roots and joining by a curved arrow
those corresponding to distinct roots
$\alpha_1,\alpha_2\in\Xi$
with $\varepsilon_C(\alpha_1)=\alpha_2$.  
\par
We can associate to any
automorphism $\eta\in\bold{A}_{\frak{h}}$ the automorphism
$\tilde\eta=\eta\circ w_{(\eta^{-1}(C),C)}\in \bold{A}_{\frak{h}}$,
that leaves $C$ fixed. We observe that $\tilde\eta(\Cal{B}(C))=
\Cal{B}(C)$, and therefore $\tilde\eta$ defines a permutation of
the nodes of the Satake diagram $\pmb{\Cal{S}}$. Thus the
quotient group $\bold{A}_{\frak{h}}/\bold{W}_{\frak{h}}$ can
be considered as the group $\bold{Aut}(\pmb{\Cal S})$ 
of {\it automorphisms of the Satake diagram $\pmb{\Cal{S}}$}. \par
Note that a Satake diagram is completely determined by the data
of: ($i$) the underlying Dynkin diagram $\Delta$, ($ii$) the color
(white or black) of the nodes, ($iii$) the involution $\epsilon_C$
on the nodes of $\Delta$.\par
The correspondence between real semisimple Lie algebras and
their Satake diagrams is one to one.\par 
We list in the appendix all the connected Satake diagrams of 
the non compact forms. We use the labels ($\roman{A\, I}$, $\hdots$,
$\roman{G}$) devised by Cartan in his classification of symmetric
spaces. We shall also consistently employ the indices attached
to the simple roots in these diagrams throughout the paper.
\medskip
We proved in Proposition {\propaf} that, for a parabolic minimal
$CR$ algebra $(\frak g,\frak q)$, the isotropy subalgebra
$\frak g_+$ contains a
Cartan subalgebra
$\frak h$ of $\frak g$ 
with maximal vector part.
First we prove\,:
\prop{If $(\frak g,\frak q)$ is an effective parabolic minimal
$CR$ algebra and $\frak h$ is a Cartan subalgebra 
of $\frak g$ with maximal vector
part 
and contained in $\frak q$, then there exists a $\sigma$-adapted 
Weyl chamber $C$ for $(\hat\frak g,\hat\frak h)$ such that
$\Cal R^+(C)\subset\Cal Q$.}\edef\teorema27-1{\number\q.\number\x}
\dimo
Modulo an inner automorphism of $\hat\frak g$, we can assume that
any given parabolic subalgebra 
$\frak q$ of $\hat\frak g$
contains a Borel subalgebra $\frak b$ of the form
\par\centerline{$\frak b\, = \,
\hat\frak h\, + \, \dsize\sum_{\alpha\in\Cal R^+(C)}{\frak g^{\alpha}}$}
\par\smallskip\noindent
for a Weyl chamber $C\in\frak C(\Cal R)$ adapted to the conjugation
$\sigma$ defined by the real form $\frak g$.
Then $\frak b\cap\frak g$ 
is contained in $\frak g_+$ and contains a maximal triangular 
subalgebra $\frak t$ of $\frak g$ (see for instance \cite{OV, 4.4, 4.5} or
\cite{Wa 1.1.3, 1.1.4}). The statement follows from 
the uniqueness stated in Theorem \teoremaag .
\enddimo
\medskip
Let $\pmb{\Cal S}$ be the Satake diagram of the semisimple real Lie algebra
$\frak g$. The nodes of $\pmb{\Cal S}$ correspond to the simple roots
$\Cal B(C)$ of a Weyl chamber $C\in\frak C(\Cal R)$ adapted to the
conjugation $\sigma$ defined by $\frak g$. Fix a subset
$\Phi$ of $\Cal B(C)$ and consider the diagram $(\pmb{\Cal S},\Phi)$
obtained from $\pmb{\Cal S}$ by adding a cross-mark on each node of
$\pmb{\Cal S}$ corresponding to a root in $\Phi$.\par
We associate to the pair $(\pmb{\Cal S},\Phi)$
the $CR$ algebra $(\frak g,\frak q_{\Phi})$ with
$\frak q_{\Phi}$ defined by (\formparab).\smallskip
Two cross-marked Satake diagrams $(\pmb{\Cal{S}},\Phi)$ and
$(\pmb{\Cal{S}},\Psi)$ are said to be {\it equivalent}
if there exists an $\varepsilon\in\bold{Aut}(\pmb{\Cal S})$ such that
$\Psi=\varepsilon(\Phi)$. 

Proposition {\teorema27-1} and Theorem {\teoremaag} yield\, :
\thm{The correspondence \par\smallskip
\centerline{$(\pmb{\Cal S},\Phi)\longleftrightarrow (\frak g,\frak q_\Phi)$}
\par\smallskip\noindent
is bijective between 
cross-marked Satake diagrams (modulo automorphisms of cross-marked
Satake diagrams) 
and minimal
effective parabolic $CR$ algebras (modulo $CR$ isomorphisms).}

\par\medskip 
\noindent
\exam{%
The diagram
$$
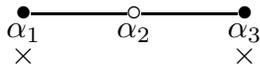
\xymatrix@M=0pt @R=1pt{
\bullet\ar @{-}[r]&\circ\ar @{-}[r]&\bullet\\
\alpha_1&\alpha_2&\alpha_3\\
\times&&\times
}
$$
corresponds to\par 
\qquad\qquad$\frak g=\frak{sl}(2,\Bbb H)\subset\frak{sl}(4,\Bbb C)$,
\par
\smallskip
\qquad\qquad
$\frak q=\{Z\in\frak{sl}(4,\Bbb C)\, | \, Z(\langle e_1\rangle)
\subset\langle e_1\rangle,\;
Z(\langle e_1,e_2,e_3\rangle)
\subset\langle e_1,e_2,e_3\rangle\,\}$,
\par\smallskip\noindent
where $e_1,e_2,e_3,e_4$ is the canonical basis of $\Bbb C^4$
with $e_1\Bbb H=\langle e_1,e_2\rangle$ and
$e_3\Bbb H=\langle e_3,e_4\rangle$.
\par
The associated minimal orbit is the $CR$ manifold $M=M^{3,2}$ 
whose points are the pairs
$(\ell_1,\ell_3)$ 
consisting of a complex line $\ell_1$ and
a complex  
$3$-plane $\ell_3$ of $\Bbb C^4$ with $\ell_1\cdot\Bbb H\subset\ell_3$.
It is strictly nondegenerate,
of $CR$ dimension $3$
and $CR$ codimension $2$;
all its nonzero Levi forms have 
one positive, one
negative and one zero eigenvalues (see for instance \cite{HN2}).
}
\exam{The diagram:
\bigskip
$$
\xymatrix@M=0pt @R=2pt{
\circ\ar @{-} [rr] \ar @/^20pt/@{<->}[rrrr]&&\bullet\ar@{-}[rr]&&\circ\\
\alpha_1&&\alpha_2&&\alpha_3\\
&&\times}
$$
\noindent
corresponds to\par
\smallskip
\qquad\qquad $\frak g=\frak{su}(1,3)\subset\hat\frak g=\frak{sl}(4,\Bbb C)$
\par\smallskip
 \qquad\qquad $\frak q=\{Z\in\frak{sl}(4,\Bbb C)\, | \,
Z(\langle e_1,e_2\rangle)\subset\langle e_1,e_2\rangle\,\}$
\par\smallskip\noindent
where $e_1,e_2,e_3,e_4$ is a basis of $\Bbb C^4$ such that\par\medskip
\centerline{$\frak{su}(1,3)=\left\{
Z\in\frak{sl}(4,\Bbb C)\, \left| \, \left(\smallmatrix
0&0&0&1\\
0&1&0&0\\
0&0&1&0\\
1&0&0&0\endsmallmatrix
\right)\, Z \, +  \, Z^* \left(\smallmatrix
0&0&0&1\\
0&1&0&0\\
0&0&1&0\\
1&0&0&0\endsmallmatrix
\right)\, =\, 0\,\right\}\right.$}
\medskip
The associated minimal orbit is a
$CR$ manifold $M=M^{3,1}$, of hypersurface type,
with a Levi form having one positive, one negative and one zero eigenvalues,
and is
weakly nondegenerate but not strictly nondegenerate.
}
\exam{The diagram\,:
\vskip1cm
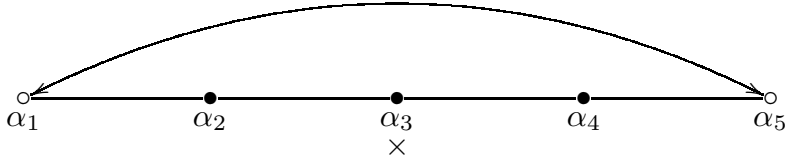
$$
\xymatrix@M=0pt @R=2pt{
\circ\ar @{-}[rr]\ar@/^30pt/ @{<->}[rrrrrrrr]&&\bullet
\ar @{-}[rr]&&\bullet\ar @{-}[rr]&&\bullet\ar @{-}[rr]&&\circ\\
\alpha_1&&\alpha_2&&\alpha_3&&\alpha_4&&\alpha_5\\
&& && \times}
$$
\noindent
corresponds to\par
\qquad\qquad $\frak g=\frak{su}(1,5)\subset\hat\frak g=\frak{sl}(6,\Bbb C)$
\par\smallskip
 \qquad\qquad $\frak q=\{Z\in\frak{sl}(6,\Bbb C)\, | \,
Z(\langle e_1,e_2,e_3\rangle)\subset\langle e_1,e_2,e_3\rangle\,\}$
\par\smallskip\noindent
where $e_1,e_2,e_3,e_4,e_5,e_6$ is a basis of $\Bbb C^4$ such that
\par\medskip
\centerline{$\frak{su}(1,5)=\left\{
Z\in\frak{sl}(4,\Bbb C)\, \left| \, \left(\smallmatrix
&&1\\
&I_3\\
1
\endsmallmatrix
\right)\, Z \, +  \, Z^* \left(\smallmatrix
&&1\\
&I_3\\
1
\endsmallmatrix
\right)\, =\, 0\,\right\}\right.$}
\medskip
The associated minimal orbit is the
$CR$ manifold $M=M^{8,1}$, of hypersurface type,
with a Levi form having two positive, two negative 
and four zero eigenvalues,
and is
weakly nondegenerate but not strictly nondegenerate.
}
\exam{The two diagrams\,:
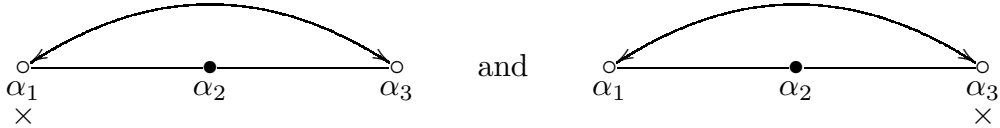
\bigskip
$$
\xymatrix@M=0pt @R=2pt{
\circ\ar @{-} [rr] \ar @/^20pt/@{<->}[rrrr]&&\bullet\ar@{-}[rr]&&\circ\\
\alpha_1&&\alpha_2&&\alpha_3\\
\times}\qquad\text{and}\qquad
\xymatrix@M=0pt @R=2pt{
\circ\ar @{-} [rr] \ar @/^20pt/@{<->}[rrrr]&&\bullet\ar@{-}[rr]&&\circ\\
\alpha_1&&\alpha_2&&\alpha_3\\
&&&&\times}
$$
are isomorphic. Indeed the map $\digamma(\alpha_i)=\alpha_{4-i}$ for
$i=1,2,3$ defines an isomorphism of cross-marked Satake diagrams.
The corresponding effective parabolic minimal $CR$-algebras
correspond to $\frak{g}=\frak{su}(1,3)$ and 
$\frak{q}=\frak{q}_{\{\alpha_1\}}$,
$\frak{q}=\frak{q}_{\{\alpha_3\}}$, respectively.
Let 
$$K=\left(\smallmatrix
0\,&0\,&0\,&1\\
0\,&1\,&0\,&0\\
0\,&0\,&1\,&0\\
1\,&0\,&0\,&0\endsmallmatrix\right)$$
and identify $\frak g$ with the Lie algebra of $4\times 4$ complex
matrices with trace zero that satisfy $X^*K+KX=0$.
The $CR$ isomorphism $(\frak g,\frak q_{\alpha_1})@>>>
(\frak g,\frak q_{\alpha_3})$
is given by the map
$\frak{su}(1,3)\ni X @>>>\, -\,{^tX}\in\frak{su}(1,3) $.
}

%% file: 07min.tex
\se{$\frak g$-equivariant fibrations}   
In this section we discuss morphisms of $CR$ algebras of the
special form $(\frak{g},\frak{q}) @>>>(\frak{g},\frak{q}')$,
for $\frak{q}\subset\frak{q}'$. They have been called in \cite{MeN5}
{\it $\frak{g}$-equivariant fibrations} and describe at the level
of $CR$ algebras the corresponding 
$\bold{G}$-equivariant smooth fibrations
$M(\frak{g},\frak{q}) @>>>M(\frak{g},\frak{q}')$.
In this section we focus on the $CR$ algebra aspects, preparing for
the applications of the next sections.
\smallskip
We keep the notation of the previous sections.
In particular, $\frak g$ is a semisimple real Lie algebra, 
$\frak h$ a Cartan subalgebra of $\frak g$ with maximal vector part,
$\Cal R=\Cal R(\Hat{\frak{g}},\Hat{\frak{h}})$, $C$ a Weyl chamber
adapted to the conjugation $\sigma$ in $\frak{h}^*_{\Bbb{R}}$
induced by the real form $\frak{g}$ of $\Hat{\frak{g}}$,
$\Cal{B}=\Cal{B}(C)$ is the set of simple roots in
$\Cal{R}^+=\Cal{R}^+(C)$.\par
Let $\Psi\subset\Phi\subset\Cal{B}$.
Then
$\frak q_\Phi\subset\frak q_\Psi$
and the identity on $\frak g$ defines
a natural $\frak g$-equivariant
morphism of $CR$ algebras (see \cite{MeN5})\,:
\form{
\pi: (\frak g,\frak q_\Phi) @>>> (\frak g,\frak q_\Psi)\, .
}\edef\formulafiba{\number\q.\number\t}
Its fiber (see \cite{MeN5}) is\,:
\form{(\frak g',\frak q'),\quad\text{where}\quad
\cases
\frak{g}'=\frak{g}\cap\frak{q}_{\Psi}=\frak{g}\cap\bar{\frak{q}}_{\Psi}\\
\Hat{\frak{g}}'=\frak{q}_{\Psi}\cap\bar{\frak{q}}_{\Psi}\\
\frak{q}'=\Hat{\frak{q}}_{\Phi}\cap\Hat{\frak{g}}'=
\Hat{\frak{q}}_{\Phi}\cap\frak{q}_{\Psi}\cap\bar{\frak{q}}_{\Psi}=
\frak{q}_{\Phi}\cap\bar{\frak{q}}_{\Psi}\, .
\endcases}
Denote by $\Cal R'$ and $\Cal Q'$ the sets of roots $\alpha\in\Cal R$ 
for which  $Hat{\frak{g}}^{\alpha}$ is contained in
$\hat\frak g'$ and $\hat\frak q'$, respectively\,:
\form{\cases
\Cal R'=\Cal Q_\Psi\cap\bar\Cal Q_\Psi\\
\Cal Q'=\Cal Q_\Phi\cap\bar\Cal Q_\Psi\, ,
\endcases}
define\,:
\form{\cases
\Cal R''=\Cal R'\cap(-\Cal R')\,=\,\Cal{Q}^r_{\Psi}\cap
\bar{\Cal{Q}}^r_{\Psi}\\
\Cal Q''=\Cal Q'\cap\Cal R'' \\
\Cal A=\Cal R'\setminus\Cal R''\,=\,\left(\Cal{Q}^n_{\Psi}
\cap\bar{\Cal{Q}}_{\Psi}\right)\cup\left(\bar{\Cal{Q}}^n_{\Psi}
\cap{\Cal{Q}}_{\Psi}\right)
\endcases}
and set\,:
\form{\cases
\hat\frak g''=\hat\frak h\oplus\bigoplus_{\alpha\in\Cal R''}\hat\frak
g^\alpha\\
\frak q''=\frak q'\cap \hat\frak g''\\
\hat\frak a=\bigoplus_{\alpha\in\Cal A}\hat\frak g^\alpha\, .
\endcases}
Then $\Cal R''$ is $\sigma$-invariant, $\hat\frak g''=\frak
q_\Psi^r\cap\bar\frak q_\Psi^r$ is reductive, $\frak
q''$ is parabolic in $\hat\frak g''$ and 
$\hat\frak a=(\frak q_\Psi^n\cap\bar\frak q_\Psi)+(\frak q_\Psi\cap\bar\frak q_\Psi^n)$ is an
ideal in $\hat\frak g'$, which is invariant with respect to the
conjugation defined by the real form $\frak{g}$.

\lem{$\hat\frak a\subset\frak q_\Phi$.}
\dimo  We first show that 
$\Cal Q_\Psi^n\cap\bar\Cal Q_\Psi\cap\Cal R_\bullet=\emptyset$. 
Assume by contradiction that there is
$\alpha\in \Cal Q_\Psi^n\cap\bar\Cal Q_\Psi\cap\Cal R_\bullet$.
From $\alpha\in\Cal Q_\Psi^n$ we obtain that
$\bar\alpha=-\alpha\in\bar\Cal Q_\Psi^n$, that is $\alpha\not\in\bar\Cal
Q_\Psi$, which gives a contradiction. \par
Since $\Cal{Q}^n_{\Psi}$ is contained in $\Cal R^+$ and
$\Cal{Q}^n_{\Psi}\cap\bar{\Cal{Q}}_{\Psi}$ does not contain imaginary
roots, also its conjugate 
$\overline{\Cal{Q}^n_{\Psi}\cap\bar{\Cal{Q}}_{\Psi}}=
\bar{\Cal{Q}}^n_{\Psi}\cap{\Cal{Q}}_{\Psi}$ is contained in
$\Cal R^+$. Hence 
$\Cal{A}\subset\Cal{R}^+\subset\Cal{Q}_{\Phi}$.
\enddimo
\smallskip
\lem{$\Cal B''=\Cal B\cap\Cal R''$ is a basis of $\Cal R''$\,.}
\dimo
Indeed, assume that $\alpha\in\Cal{R}''$ is the sum of two positive
roots: $\alpha=\beta+\gamma$ with $\beta,\gamma\in\Cal{R}^+$.
Then $\alpha\in\Cal{Q}^r_{\Psi}$ implies that
also $\beta,\gamma\in\Cal{Q}^r_{\Psi}$. 
If $\beta,\gamma\notin\Cal{R}_{\bullet}$, then by the same argument
applied to $\bar\alpha=\bar\beta+\bar\gamma\in\Cal{Q}^r_{\Psi}$
we obtain that $\beta,\gamma$ also belong to $\bar{\Cal{Q}}^r_{\Psi}$
and hence to $\Cal{R}''$. \par
Consider now the case where, for instance, $\beta\in\Cal{R}_{\bullet}$.
Then $\bar{\beta}=-\beta\in\Cal{Q}^r_{\Psi}$ implies that
$\beta\in\Cal{R}''$ and therefore  $\gamma=\alpha-\beta\in\Cal{R}''$,
showing that also in this case $\alpha$ is not simple in
$\left[\Cal{R}''\right]^+=\Cal R^+\cap\Cal{R}''$. This shows that
$\Cal{B}''$ is exactly the set of simple roots in $\left[\Cal{R}''\right]^+$,
and thus a basis of $\Cal{R}''$.\enddimo
\par
We have obtained:
\prop{The $CR$ algebra $(\frak{g}'',\frak{q}'')$ is parabolic minimal.
Its
cross-marked Satake diagram $({\pmb{\Cal S}''},\Phi'')$  
is the subdiagram of $(\pmb{\Cal{S}},\Phi)$ 
consisting of the
simple roots $\alpha$ such that\,:\par\smallskip\centerline{
either \quad ($i$) $\alpha\in\Cal{R}_{\bullet}\setminus\Psi$, \quad or\quad
($ii$) $\alpha\not\in\Cal R_\bullet$ and
$\left(\{\alpha\}\cup\roman{supp}(\bar\alpha)
\right)\cap\Psi=\emptyset$.}\par\smallskip
The cross-marks are left on the nodes corresponding to roots in
$\Phi\cap\Cal B''$.\qed}\edef\propgef{\number\q.\number\x}
\par
We say that a Satake diagram is {\it $\sigma$-connected} if 
either it is connected or
consists of two connected components, joined by curved arrows.
\thm{Let $(\formulafiba)$ be a $\frak g$-equivariant fibration. Then the
effective quotient of its fiber is the parabolic minimal $CR$ algebra 
whose cross-marked Satake diagram consists of the union of all 
$\sigma$-connected components of the diagram ${\pmb{\Cal S}''}$ described
in Proposition \propgef, containing at least one cross-marked
node.\qed}
\edef\thmgef{\number\q.\number\x}
\exam{%
Let $\frak{g}=\frak{su}(1,3)$ and let
$\Phi=\{\alpha_1,\alpha_2\}$, $\Psi=\{\alpha_1\}$. Then
the cross-marked Satake diagrams corresponding to the $CR$ algebra
$(\frak g,\frak{q}_\Phi)$, the basis $(\frak g,\frak{q}_\Psi)$
and the corresponding effective fiber are given by:
\bigskip
$$\CD
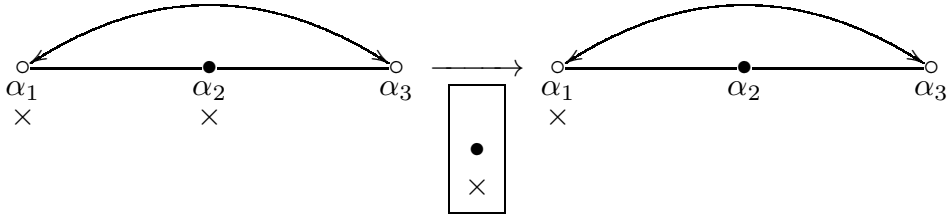
{\xymatrix@M=0pt @R=2pt{
\circ\ar @{-} [rr] \ar @/^20pt/@{<->}[rrrr]&&\bullet\ar@{-}[rr]&&\circ\\
\alpha_1&&\alpha_2&&\alpha_3\\
\times&&\times}}@>>{\boxed{\matrix 
{}\\
\bullet\\
\times
\endmatrix}}>
{\xymatrix@M=0pt @R=2pt{
\circ\ar @{-} [rr] \ar @/^20pt/@{<->}[rrrr]&&\bullet\ar@{-}[rr]&&\circ\\
\alpha_1&&\alpha_2&&\alpha_3\\
\times}}\endCD
$$
In the case $\Psi=\{\alpha_2\}$ we have instead:
\bigskip
$$\CD
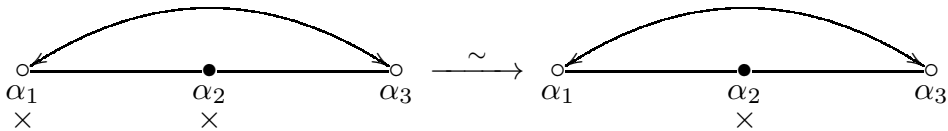
{\xymatrix@M=0pt @R=2pt{
\circ\ar @{-} [rr] \ar @/^20pt/@{<->}[rrrr]&&\bullet\ar@{-}[rr]&&\circ\\
\alpha_1&&\alpha_2&&\alpha_3\\
\times&&\times}}@>{\sim}>>
{\xymatrix@M=0pt @R=2pt{
\circ\ar @{-} [rr] \ar @/^20pt/@{<->}[rrrr]&&\bullet\ar@{-}[rr]&&\circ\\
\alpha_1&&\alpha_2&&\alpha_3\\
&&\times}}\endCD
$$ 
The fiber is trivial and the map is a $CR$ morphism, but not a
$CR$ isomorphism. The corresponding map 
$M(\frak{g},\frak{q}_{\Phi})@>>>M(\frak{g},\frak{q}_{\Psi})$ is
an analytic diffeomorphism and a $CR$ map, but not a $CR$ diffeomorphism. }

\par
A $\frak{g}$-equivariant morphism of $CR$ algebras (\formulafiba)
is a {\it $CR$-fibration} if the quotient map
\form{
\frak{q}_{\Phi}/\left(\frak{q}_{\Phi}\cap\overline{\frak{q}}_{\Phi}\right)
@>>> \frak{q}_{\Psi}/\left(\frak{q}_{\Psi}\cap\overline{\frak{q}}_{\Psi}\right)
}
is onto. Set $M_{\Phi}=M(\frak{g},\frak{q}_{\Phi})$,
$M_{\Psi}=M(\frak{g},\frak{q}_{\Psi})$, and
$F=M(\frak{g}'',\frak{q}'')$. The condition that (\formulafiba)
is a {\it $CR$-fibration} is equivalent to the fact that
every point of $M_{\Phi}$ has an open
neighborhood which is $CR$ diffeomorphic to the product of an open
submanifold of $M_{\Psi}$ and $F$. \par
We have the criterion\,:
\prop{The following conditions are equivalent:
\roster
\item"$(i)$" $(\formulafiba)$ is a $CR$-fibration;
\item"$(ii)$" $\Cal{Q}^r_{\Psi}\setminus\Cal{Q}_{\Phi}\subset\bar{\Cal{Q}}_{\Psi}$; 
\item"$(iii)$" $\Cal{Q}^r_{\Psi}\cap\Cal{Q}^n_{\Phi}\subset\bar{\Cal{Q}}^r_{\Psi}$.
\endroster
}\edef\propcrf{\number\q.\number\x}
\dimo 
First we prove the equivalence $(i)\Leftrightarrow(ii)$.
A necessary and sufficient condition in order that (\formulafiba)
be a $CR$-fibration is that the sum of the $CR$ dimensions of
$(\frak{g},\frak{q}_{\Psi})$ and of the fiber $(\frak{g}',\frak{q}')$
equals the $CR$-dimension of the total space $(\frak{g},\frak{q}_{\Phi})$\,:
$$\matrix\format\r&\l\\
\roman{dim}_{\Bbb{C}}{\frak{q}_{\Phi}}-
\roman{dim}_{\Bbb{C}}{\frak{q}_{\Phi}\cap\bar{\frak{q}}_{\Phi} }&=
\roman{dim}_{\Bbb{C}}{\frak{q}_{\Psi}}-
\roman{dim}_{\Bbb{C}}{\frak{q}_{\Psi}\cap\bar{\frak{q}}_{\Psi} }\\
& + \roman{dim}_{\Bbb{C}}{\frak{q}_{\Phi}\cap\bar{\frak{q}}_{\Psi} }
- \roman{dim}_{\Bbb{C}}{\frak{q}_{\Phi}\cap\bar{\frak{q}}_{\Phi} }
\, .
\endmatrix$$
Since all subspaces considered in this formula contain $\Hat{\frak{h}}$,
this is equivalent to\,:
$$|\Cal{Q}_{\Phi}| \, = \, |\Cal{Q}_{\Psi}| - |\Cal{Q}_{\Psi}\cap 
\bar{\Cal{Q}}_{\Psi}| + |\Cal{Q}_{\Phi}\cap \bar{\Cal{Q}}_{\Psi}|\, =\,
|\Cal{Q}_{\Psi}\setminus 
\bar{\Cal{Q}}_{\Psi}| + |\Cal{Q}_{\Phi}\cap \bar{\Cal{Q}}_{\Psi}|\, ,
\tag $*$ $$
(where we used $|A|$ for the number of elements of the finite set
$A$). Since $\Cal{Q}_{\Phi}\subset\Cal{Q}_{\Psi}$, we always have\,:
$$\Cal{Q}_{\Phi}\subset \left(\Cal{Q}_{\Psi}\setminus 
\bar{\Cal{Q}}_{\Psi}\right)\cup \left(\Cal{Q}_{\Phi}\cap \bar{\Cal{Q}}_{\Psi}
\right)\, .$$
The two sets on the right hand side are disjoint. Hence ($*$)
is equivalent to\,:
$$
\Cal{Q}_{\Psi}\setminus \bar{\Cal{Q}}_{\Psi}\subset\Cal{Q}_{\Phi}.
$$
As $\Cal Q^n_\Psi\subset\Cal R^+\subset\Cal Q_\Phi$, this is equivalent to
$$
\Cal{Q}^r_{\Psi}\setminus \Cal{Q}_{\Phi}\subset\bar{\Cal{Q}}_{\Psi}.
$$
\par\smallskip
Next we prove that $(ii)\Rightarrow(iii)$. We distinguish several cases.
\par
If $\alpha\in\Cal Q^r_\Psi\cap\Cal R_\bullet$, then 
$\bar\alpha=-\alpha\in\Cal Q^r_\Psi$, that is $\alpha\in\bar\Cal Q^r_\Psi$.
\par
If $\alpha\in\Cal Q^r_\Psi\cap\Cal Q^n_\Phi$ and 
$\alpha\not\in\Cal R_\bullet$, 
then $\bar\alpha\succ0$, hence $\alpha\in\bar\Cal Q_\Psi$. On the other hand
$-\alpha\in\Cal Q^r_\Psi\setminus\Cal Q_\Phi$ and, by $(ii)$, 
$-\alpha\in\bar\Cal Q_\Psi$, thus $\alpha\in\bar\Cal Q^r_\Psi$.
\par\smallskip
Finally we prove that $(iii)\Rightarrow(ii)$. Let 
$\alpha\in\Cal Q^r_\Psi\setminus\Cal Q_\Phi$.
Then $-\alpha\in\Cal Q^r_\Psi\cap\Cal Q^n_\Phi$, and $(iii)$ 
implies
that $-\alpha
\in\bar\Cal Q^r_\Psi$, which is
equivalent to  $\alpha\in\bar\Cal Q^r_\Psi$.
\enddimo\edef\proposizionedf{\number\q.\number\x}
\par
In particular, we obtain\,:
\prop{If $\bar{\Cal{Q}}_{\Psi}={\Cal{Q}}_{\Psi}$, then (\formulafiba)
is a $CR$-fibration.}\edef\proposizionedg{\number\q.\number\x}
\dimo Indeed condition ($iii$) of Proposition {\proposizionedf}
is trivially satisfied if $\bar{\Cal{Q}}_{\Psi}={\Cal{Q}}_{\Psi}$.
\enddimo
We recall (see \cite{MeN5}) that a $CR$ algebra $(\frak{g},\frak{q})$
is {\it totally complex} if $\frak{q}+\bar{\frak{q}}=\Hat{\frak{g}}$.
This condition is equivalent to 
$\frak{g}+\frak{q}=\Hat{\frak{g}}$ and to
the fact that every homogeneous
$CR$ manifold $M$ with associated $CR$ algebra $(\frak{g},\frak{q})$
is actually a complex manifold.
\prop{If $\bar{\Cal{Q}}_{\Psi}\cup{\Cal{Q}}_{\Psi}=
\bar{\Cal{Q}}_{\Phi}\cup{\Cal{Q}}_{\Phi}$, then
(\formulafiba)
is a $CR$-fibration
with a totally complex fiber.}\edef\proposizionedh{\number\q.\number\x}
\dimo Indeed we obtain: $\Cal{Q}_{\Psi}\setminus\Cal{Q}_{\Phi}\subset
\bar{\Cal{Q}}_{\Phi}\subset\bar{\Cal{Q}}_{\Psi}$, and hence
($ii$) of Proposition {\proposizionedf} follows because
$\Cal{Q}_{\Psi}^r\supset\Cal{Q}_{\Phi}^r$,
$\Cal{Q}_{\Psi}^n\subset\Cal{Q}_{\Phi}^n$. 
\par
To show that the fiber is totally complex, we need to verify that
$\frak{q}_{\Psi}\cap\bar{\frak{q}}_{\Psi}=
\frak{q}_{\Phi}\cap\bar{\frak{q}}_{\Psi}+
\frak{q}_{\Psi}\cap\bar{\frak{q}}_{\Phi}\, .$ This is obvious
because $\frak{q}_{\Phi}\subset\frak{q}_{\Psi}\subset\frak{q}_{\Phi}
+\bar{\frak{q}}_{\Phi}$. 
\enddimo

Our next aim is to
characterize $\frak g$-equivariant $CR$ fibrations
in terms of cross marked Satake diagrams. For this
we introduce 
some notation.
\par
The {\it component}
$\check\Psi(\alpha)$ of a root $\alpha\in\Cal B(C)$ 
is the set of roots $\beta\in\Cal{B}(C)$ belonging
to 
the connected component of the node corresponding to $\alpha$ 
in the graph obtained from $\pmb{\Cal S}$
by deleting those nodes that correspond to roots in 
$\Psi\setminus\{\alpha\}$ and the
lines and arrows issuing from them.
\par
Given a subset $\Cal E$ of $\Cal{B}(C)$, its
{\it exterior boundary} $\partial_e{\Cal E}$ 
in $\pmb{\Cal S}$ is
the set of roots $\alpha$ in $\Cal{B}(C)\setminus\Cal{E}$
such that, for some $\beta\in\Cal{E}$, $\alpha+\beta\in\Cal{R}$.
\par
It will be convenient in the following to identify the nodes of
$\pmb{\Cal{S}}$ with the corresponding roots in $\Cal{B}(C)$.
In particular,
for a connected subset ${\Cal E}$ of a Satake diagram
$\pmb{\Cal S}$, we set 
$\delta({\Cal E})=\sum_{\alpha\in{\Cal E}}\alpha\in\Cal R$.
\par
We recall the notation $\Xi=\Cal B(C)\setminus\Cal R_\bullet$ for
the set of non imaginary simple roots.
\lem{If $\alpha\in\Cal R\setminus\Cal R_\bullet$, then
$$
\supp(\bar\alpha)\supset\bigl(\partial_e(\supp(\alpha))\cap\Cal R_\bullet\bigr)
\cup\epsilon_C\bigl(\supp(\alpha)\setminus\Cal R_\bullet\bigr).
$$}\edef\lemcrf{\number\q.\number\x}
\dimo
By inspecting the conjugation diagrams in \cite{Ar}, we find that, if $\alpha
\in\Xi$\,:
\form{
\supp(\bar\alpha)=\bigl(\check\Xi(\alpha)\setminus\{\alpha\}\bigr)
\cup\{\epsilon_C(\alpha)\}.
}\xdef\formconjsimpl{\number\q.\number\t}
If $\alpha=\sum k_i\alpha_i\in\Cal R\setminus\Cal R_\bullet$, then
$$
\supp(\bar\alpha)\supset\Bigl(
\bigcup\Sb k_i>0\\\alpha_i\in\Xi\endSb\supp(\bar\alpha_i)\Bigr)
\setminus\bigl( \supp(\alpha)\cap\Cal R_\bullet \bigr),
$$
in particular $\supp(\bar\alpha)$ contains 
$\epsilon_C\bigl(\supp(\alpha)\setminus\Cal R_\bullet\bigr)$.
\par
If $\beta\in\partial_e\big(\supp(\alpha)\big)\cap\Cal R_\bullet$, 
then, since $\supp(\alpha)\not\subset\Cal R_\bullet$, there 
exists $\alpha_i\in\supp(\alpha)\cap\Xi$ such that 
$\beta\in\check\Xi(\alpha_i)$. 
This implies that $\supp(\bar\alpha)\ni\beta$.
\enddimo
\thm{A necessary and sufficient condition for $(\formulafiba)$ to be
a $CR$ 
$\frak{g}$-equivariant
fibration is that for every $\alpha\in\Phi\setminus\Psi$ 
either one
of the following conditions hold:
\roster
\item"$(i)$" $\check\Psi(\alpha)\subset\Cal R_\bullet$\,;
\item"$(ii)$" $\check\Psi(\alpha)\not\subset\Cal R_\bullet$, 
$\epsilon_C(\check\Psi(\alpha)\setminus\Cal R_\bullet)\cap\Psi=\emptyset$, 
and $\partial_e\check\Psi(\alpha)\cap\Cal R_\bullet=\emptyset$\, .
\endroster
}\edef\thmcrfibr{\number\q.\number\x}
\dimo
Condition $(ii)$ in Proposition \propcrf\ is equivalent to the assertion
that, for every root $\beta$:
\form{
\left.\aligned
\supp(\beta)\cap\Psi &= \emptyset\\
\supp(\beta)\cap\Phi &\neq \emptyset
\endaligned\right\}\ \Longrightarrow\ 
\supp(\bar\beta)\cap\Psi = \emptyset.
}\par
Fix $\alpha\in\Phi\setminus\Psi$ and let $\beta=\delta(\check\Psi(\alpha))$.
Then, according to Lemma \lemcrf, either $\beta\in\Cal R_\bullet$ or
$\supp(\bar\beta)\supset
\epsilon_C\big(\check\Psi(\alpha)\setminus\Cal R_\bullet\big)\cup
\big(\partial_e\check\Psi(\alpha)\setminus\Cal R_\bullet\big)$, showing that 
either $(i)$ or $(ii)$ 
must be valid.
\par
Fix again $\alpha\in\Phi\setminus\Psi$ and let 
$\alpha_j\in\check\Psi(\alpha)$. If $\alpha_j\in\Cal R_\bullet$
then $\bar\alpha_j=-\alpha_j$ and $\supp(\bar\alpha_j)\cap\Psi
=\emptyset$. 
If $\alpha_j\not\in\Cal R_\bullet$, formula $(\formconjsimpl)$
implies that either $\supp(\bar\alpha_j)\subset\check\Psi(\alpha)$
or $\bar\alpha_j=\epsilon_C(\alpha_j)$. In both cases
 $\supp(\bar\alpha_j)\cap\Psi=\emptyset$. 
For a generic $\beta\in\Cal R\setminus\Cal R_\bullet$  such
that $\supp(\beta)\subset\check\Psi(\alpha)$ we have that:
$$
\supp(\bar\beta)\subset
\bigcup_{\alpha_j\in\supp(\beta)}\supp(\bar\alpha_j),
$$
hence $\supp(\bar\beta)\cap\Psi=\emptyset$.
\enddimo

%% file: 08min.tex
\se{A rigidity theorem for minimal orbits} 
In this section we will discuss some topological properties
of homogeneous $CR$ manifolds, having an associated $CR$ algebra
which is parabolic minimal. We
generalize here some results 
proved in \cite{MeN1}, \cite{MeN4} in the case
of parabolic minimal $CR$ algebras corresponding to 
{\it semisimple Levi-Tanaka algebras}.
\par
Let $(\frak g,\frak q_{\Phi})$ be an effective parabolic minimal
$CR$ algebra, with associated cross-marked Satake diagram
$(\pmb{\Cal{S}},\Phi)$. We say that  $(\frak g,\frak q_{\Phi})$
has property $(F)$ if $\Phi$ does not contain any real root.
\par
Let (\formulacartan) be a
Cartan decomposition of $\frak{g}$ with (\formuladopocartan).
Denote by $\frak{k}^{(1)}=[\frak{k},\frak{k}]$ the 
maximal semisimple
ideal of $\frak{k}$ and set
$\frak{k}_+=\frak k\cap\frak{g}_+=\frak k\cap\frak q_{\Phi}$. We have\,:
\prop{Assume that $(\frak g,\frak q_{\Phi})$
has property $(F)$. Then\,:
\form{\frak{k}=
\frak{k}^{(1)}+\frak{k}_+\, .}}\edef\formulada{\number\q.\number\t}
\dimo
Let $\Hat{\frak{k}}$ be the complexification of $\frak{k}$. Since
$\Hat{\frak{k}}_+=\Bbb{C}\otimes\frak{k}_+=\Hat{\frak{k}}\cap
\frak{q}_{\Phi}\cap\bar{\frak{q}}_{\Phi}$, our contention 
(\formulada) is equivalent to\,:
\form{\Hat{\frak{k}}=\Hat{\frak{k}}^{(1)}\,
+\,\Hat{\frak{k}}_+\,.}\edef\formuladb{\number\q.\number\t}
For $\alpha\notin\Cal{R}_{\bullet}$, set
$\frak{k}^{\alpha}\, = \, 
\left(\frak{g}^{\alpha}\oplus\frak{g}^{-\bar\alpha}\right)\cap
\Hat{\frak{k}}$.
This is a $1$-dimensional subspace of
$\Hat{\frak{g}}$,  and we obtain the direct sum decomposition\,:
\form{\Hat{\frak{k}}\, = \, \Hat{\frak{h}}^+\oplus
\dsize\sum_{\alpha\in\Cal{R}_{\bullet}}{\Hat{\frak{g}}^{\alpha}}\,
\oplus\,
\dsize\sum_{\smallmatrix
\alpha\succ 0\\
\bar\alpha\succ 0
\endsmallmatrix}{{\frak{k}}^{\alpha}}\, .}
First we note that $ \Hat{\frak{h}}^+\subset\frak{q}_{\Phi}\cap
\bar{\frak{q}}_{\Phi}$\,.
Next we observe that, if $\alpha$ is not real, then
$\alpha(\frak{h}^+)\neq\{0\}$ and hence\,:
$\frak{g}^\alpha=[\Hat{\frak{h}}^+,\frak{g}^\alpha]\subset
\Hat{\frak{k}}^{(1)}$ if $\alpha\in\Cal{R}_{\bullet}$, and
$\frak{k}^{\alpha}=[\Hat{\frak{h}}^+,\frak{k}^\alpha]\subset
\Hat{\frak{k}}^{(1)}$ if $\alpha\succ 0$, $\bar\alpha\neq\pm\alpha$.
\par
By property $(F)$, for a simple real $\alpha$ we have
$\alpha\in\Cal{Q}^r\cap\bar{\Cal{Q}}^r$ and hence
$\frak{k}^{\alpha}\subset\Hat{\frak{k}}_+$. \par
To complete the proof,
we argue by contradiction. If (\formuladb) is not valid, there exists
some real root $\alpha\succ 0$,
minimal with respect to "$\prec$", with
$\frak{k}^{\alpha}\not\subset\Hat{\frak{k}}^{(1)}\cup
\Hat{\frak{k}}_+$. This $\alpha$ is not simple, and hence
$\alpha=\beta+\gamma$, with $\beta,\gamma\succ 0$ and we can assume
that $\beta$ is not imaginary. We consider the different cases.
\par
If $\gamma\in\Cal{R}_{\bullet}$, then 
$\frak{k}^{\alpha}=[\frak{g}^{\gamma},\frak{k}^{\beta}]\subset
\Hat{\frak{k}}^{(1)}$.
\par
We note that, for $\beta,\gamma\notin\Cal{R}_{\bullet}$,
the commutator $[\frak{k}^{\beta},\frak{k}^{\gamma}]$ 
is contained in a sum of eigenspaces $\Hat{\frak{g}}^\eta$ where
$\eta$ is a root in $\{\beta+\gamma, -\bar\beta-\bar\gamma,
\beta-\bar\gamma,\gamma-\bar\beta\}$. Hence\,:\par
If $\gamma\notin\Cal{R}_{\bullet}$ and $\beta-\bar\gamma\notin\Cal{R}$,
then $\frak{k}^{\alpha}=[\frak{k}^{\beta},\frak{k}^{\gamma}]\subset
\Hat{\frak{k}}^{(1)}$.\par
If $\gamma\notin\Cal{R}_{\bullet}$ and $\beta-\bar\gamma$ is a positive
root, we have that $0\prec \beta-\bar\gamma\prec\beta\prec\alpha$.
Thus, by the assumption that $\alpha$ is minimal, 
$\frak{k}^{\beta-\bar\gamma}\subset\Hat{\frak{k}}^{(1)}+\Hat{\frak{k}}_+$
and
$\frak{k}^{\alpha}\subset [\frak{k}^{\beta},\frak{k}^{\gamma}]
+\frak{k}^{\beta-\bar\gamma}\subset\Hat{\frak{k}}^{(1)}+\Hat{\frak{k}}_+$.
\par
Analogously, if  $\gamma\notin\Cal{R}_{\bullet}$ and
$\gamma-\bar\beta$ is a positive root, then
$0\prec \gamma-\bar\beta\prec \gamma\prec \alpha $;
by the assumption that $\alpha$ is minimal,
$\frak{k}^{\gamma-\bar\beta}\subset\Hat{\frak{k}}^{(1)}+\Hat{\frak{k}}_+$
and
$\frak{k}^{\alpha}\subset [\frak{k}^{\beta},\frak{k}^{\gamma}]
+\frak{k}^{\gamma-\bar\beta}\subset\Hat{\frak{k}}^{(1)}+\Hat{\frak{k}}_+$.
\enddimo\edef\teoremada{{\number\q.\number\x}}
As a corollary, we obtain:
\thm{Let $M$ be a connected 
homogeneous $CR$ manifold whose associated
$CR$ algebra is a parabolic minimal $(\frak{g},\frak{q}_{\Phi})$
that satisfies property $(F)$. Then $M$ is compact and has
a finite fundamental group.}
\dimo Let $\bold{G}$ denote the semisimple group with Lie algebra
$\frak{g}$ that acts transitively on $M$. Let (\formulacartan)
be a Cartan decomposition of $\frak{g}$ and let $\bold{K}^{(1)}$
be the analytic subgroup generated by $\frak{k}^{(1)}$. Since
$\frak{k}^{(1)}$ is semisimple and compact, 
the group $\bold{K}^{(1)}$ is 
semisimple and compact. By Proposition \teoremada ,
$\bold{K}^{(1)}$ is transitive on $M$,
because for each $p\in M$ the orbit
$\bold{K}^{(1)}\cdot p$ is open and closed in $M$, and
hence fore coincides with $M$, which, in particular, is compact. 
The universal covering $\Tilde{\bold{K}}^{(1)}$ of
$\bold{K}^{(1)}$ is compact. If $\Tilde{\bold{K}}^{(1)}_+$
is the analytic subgroup of $\Tilde{\bold{K}}^{(1)}$ generated by
$\frak{k}_+^{(1)}=\frak{k}^{(1)}\cap\frak{q}_{\Phi}$, then
$\Tilde{M}=\Tilde{\bold{K}}^{(1)}/\Tilde{\bold{K}}^{(1)}_+$
is simply connected and is the universal covering of $M$.
Therefore, having a compact universal covering, $M$ has a finite
fundamental group.\enddimo\edef\teoremadb{\number\q.\number\x}
\exam{ 
Fix a positive integer $p$, and let $\frak{g}\simeq\frak{su}(p,p)$
be the set of $(2p)\times(2p)$ complex matrices $Z$ with $0$ trace
that satisfy\,:\par\medskip\centerline{
$Z^* K\,+\,K\, Z\, = \, 0\, $, where 
$K=\left({\smallmatrix
I_p\\
&-I_p
\endsmallmatrix}\right)\,.$}\par\medskip
Let $e_1,\hdots,e_{2p}$ be the canonical basis of $\Bbb{C}^{2p}$
and let $\frak{q}\subset\Hat{\frak{g}}\simeq\frak{sl}(2p,\Bbb{C})$
be the set of $(2p)\times (2p)$ matrices in $\frak{sl}(2p,\Bbb{C})$
such that \par\smallskip\centerline{
$Z(\langle e_1+e_{p+1}\, ,\hdots\,,e_{p}+e_{2p}\rangle)\subset
\langle e_1+e_{p+1}\, ,\hdots\,,e_{p}+e_{2p}\rangle\, .$}\par
\smallskip
Then $(\frak{g},\frak{q})$ is parabolic minimal. The corresponding
$CR$ manifold $M=M(\frak{g},\frak{q})$ is the Grassmannian of $p$-planes
$\ell_p$ in $\Bbb{C}^{2p}$ which are totally isotropic for $K$
(i.e. $v^*Kv=0$ for all $v\in\ell_p$). We have\par\smallskip
\centerline{
$M\simeq\left\{\ell_p=\{(v,u(v))\in\Bbb{C}^{2p}\, | \,
v\in\Bbb{C}^p\}\,\left| \, u\in\bold{U}(p)\right\}\right.
\simeq\bold{U}(p)$}\par\smallskip\noindent
where $\bold{U}(p)$ is the group of unitary $p\times p$ matrices,
i.e. $\bold{U}(p)=\{u\in\bold{GL}(p,\Bbb{C})\, | \, u^*u=\pmb{1}\}$.
Then $\pi_1(M)\simeq\Bbb{Z}$ is infinite.
In this case the cross-marked Satake diagram is\,:
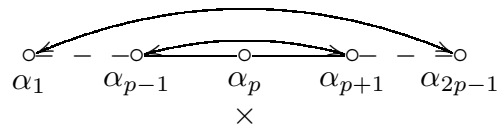
\bigskip
$$ \xymatrix @M=0pt @R=3pt @!C=10pt{
\circ \ar@{--}[r]\ar@/^15pt/@{<->}[rrrr]&\circ\ar@{-}[r]\ar@/^5pt/@{<->}[rr]
&\circ\ar@{-}[r]&\circ\ar@{--}[r]&\circ\\
\alpha_1&\alpha_{p-1}&\alpha_{p}&\alpha_{p+1}&\alpha_{2p-1}\\
&&\times}
$$
and  property ($F$) is not valid, since $\Phi=\{\alpha_p\}$
consists of a real root.}
\prop{Let $(\frak{g},\frak{q}_{\Phi})$ be an effective parabolic minimal
$CR$ algebra. Then there exists a unique minimal (with respect to
inclusion) parabolic subalgebra $\frak q$ of $\Hat{\frak{g}}$
that is contained in $\frak{q}_{\Phi}$ and satisfies
$\frak{q}_{\Phi}\cap\bar\frak{q}_{\Phi}=\frak{q}\cap\bar\frak{q}$.
\par
The parabolic $CR$ algebra $(\frak{g},\frak{q})$ is minimal
and we have\,:
\roster
\item"($i$)" $\frak{q}=\frak{q}_{\Psi}$ for a
$\Psi\supset\Phi$\,;
\item"($ii$)" $\frak{q}_{\Psi}+\bar{\frak{q}}_{\Psi}$ is
a Lie subalgebra of $\Hat{\frak{g}}$\,.
\endroster
Moreover, if property $(F)$
is valid for $(\frak{g},\frak{q}_{\Phi})$, then
it is also valid 
for $(\frak{g},\frak{q}_{\Psi})$.}\edef\teoremadc{{\number\q.\number\x}}
\dimo Let $A\in\frak{h}_{\Bbb{R}}$ be such that\,:\par\smallskip
\centerline{$\Cal{Q}_{\Phi}\, = \,
\{\alpha\in\Cal{R}\, | \, \alpha(A)\geq 0\}$\, }\par\smallskip\noindent
(for this characterization of the parabolic set of roots see
for instance \cite{Wa}). Set $A=A_-\,+\,i\,A_+$, with $A_-\in\frak{h}^-$
and $A_+\in\frak{h}^+$. Fix a real 
positive $\epsilon$ sufficiently small, so
that $\left|\alpha(i\,A_+)\right|<\epsilon^{-1}\left|\alpha(A_-)\right|$
whenever $\alpha(A_-)\neq 0$, and set $B=A_-\,+\,i\,\epsilon\, A_+$.
Then $B\in\frak{h}_{\Bbb{R}}$. We observe that
\par\smallskip
\centerline{$\Cal{Q}\, = \,
\{\alpha\in\Cal{R}\, | \, \alpha(B)\geq 0\}$\, }\par\smallskip\noindent
is the the set of roots of a parabolic minimal corresponding to
some $\Psi\subset\Cal{B}$. Indeed
$\alpha(A_-)>0$ implies that $\alpha(B)>0$ and
$\alpha(B)=\epsilon\alpha(A)$ when $\alpha(A_-)=0$.
This shows in particular that $\Cal{R}^+\subset\Cal{Q}_{\Psi}$,
and hence $\Cal{Q}=\Cal{Q}_{\Psi}$ for some subset $\Psi$
of simple roots of $\Cal{R}^+$. This observation also yields
($i$), while $\frak{q}_{\Psi}+\bar{\frak{q}}_{\Psi}$ is the parabolic
subalgebra corresponding to the set \par\smallskip
\centerline{$\Cal{Q}'\, = \,
\{\alpha\in\Cal{R}\, | \, \alpha(A_-)\geq 0\}$\, .}\par\smallskip\noindent
For a real root $\alpha$, we have
$\alpha(A)=\alpha(B)$, and hence the two parabolic sets
$\Cal{Q}_{\Phi}$ and $\Cal{Q}_{\Psi}$ contain the same
real roots. This implies that they either both have or both do not
have property $(F)$.\par
Let us prove that $\Phi\subset\Psi$.
Let $\alpha\in\Cal{Q}_{\Phi}^n\subset\Cal{R}^+$. If $\alpha(A_-)=0$,
then $\alpha(B)=\epsilon\alpha(A)>0$ and $\alpha\in\Cal{Q}_{\Psi}^n$.
If $\alpha(A_-)\neq 0$, then $\alpha\in\Cal{R}^+\setminus\Cal{R}_\bullet$,
and hence $\alpha,\bar\alpha\in\Cal{R}^+\subset\Cal{Q}_{\Phi}$ implies
that $\alpha(A_-)\geq |\alpha(iA_+)|$. Thus $\alpha(B)>0$ and
again $\alpha\in\Cal{Q}_{\Psi}^n$. Since $\Cal{Q}^n_{\Phi}\subset
\Cal{Q}^n_{\Psi}$, we have $\Phi\subset\Psi$.
\smallskip
Finally we show that $\frak{q}=\frak{q}_{\Psi}$ satisfies the minimality
condition. 
To this aim, we will show that every
parabolic Lie subalgebra $\frak{q}'$
of $\Hat{\frak{g}}$ with $\frak{q}'\subset\frak{q}_{\Phi}$ and
$\frak{q}'\cap\bar\frak{q}'=\frak{q}_{\Phi}\cap
\bar\frak{q}_{\Phi}$ contains $\frak{q}_{\Psi}$.
Since $\frak q'$ contains $\Hat{\frak{h}}$, we have
$\frak{q}'=\Hat{\frak{h}}\oplus\sum_{\alpha
\in\Cal{Q}'}{\Hat{\frak{g}}^\alpha}$ for a parabolic subset
$\Cal{Q}'\subset\Cal{R}$. \par
We claim that $\Cal{R}^+\subset\Cal{Q}'$. Indeed, for
$\alpha\in\Cal{R}^+\setminus\Cal{R}_\bullet$, we have
$\alpha\in\Cal{Q}_{\Phi}\cap\bar\Cal{Q}_{\Phi}=\Cal{Q}'\cap\bar\Cal{Q}'
\subset\Cal{Q}'$. If $\alpha\in\Cal{R}^+\cap\Cal{R}_\bullet$, then
either $\alpha\in\Cal{Q}^r_{\Phi}$ and again
$\alpha\in\Cal{Q}_{\Phi}\cap\bar\Cal{Q}_{\Phi}=\Cal{Q}'\cap\bar\Cal{Q}'
\subset\Cal{Q}'$,
or $\alpha\in\Cal{Q}^n_{\Phi}$ and $-\alpha\notin\Cal{Q}_\Phi\supset
\Cal{Q}'$ implies $\alpha\in\Cal{Q}'$, because $\Cal{Q}'$ is parabolic.
This implies that $(\frak{g},\frak{q}')$ is parabolic minimal and
$\frak{q}'=\frak{q}_{\Psi'}$ for some $\Psi'\supset\Phi$.\par  
To prove minimality, we need to show that $\Psi'\subset\Psi$, i.e.
that $\alpha(B)>0$ for all $\alpha\in\Psi'$. First we observe that
$\Psi'\cap\Cal{R}_\bullet=\Psi\cap\Cal{R}_\bullet=\Phi\cap\Cal{R}_\bullet$,
as we showed above that ${\Cal{Q}'}^n\cap\Cal{R}_\bullet=
\Cal{Q}^n_{\Phi}\cap\Cal{R}_\bullet$ for all parabolic sets
$\Cal{Q}'$ with 
$\Cal{Q}'\cap\bar\Cal{Q}'=\Cal{Q}_{\Phi}\cap\bar\Cal{Q}_{\Phi}$. 
\par
For $\alpha\in\Psi'\setminus\Cal{R}_\bullet$, either 
$\alpha\in\Phi\subset\Psi$, or $\alpha(A)=0$ and $\bar\alpha(A)>0$.
This implies that $\alpha(A_-)>0$ and hence $\alpha(B)>0$.
\par
This completes the proof.\enddimo
\cor{Let $\Phi\subset\Psi\subset\Cal{B}$ be as in the statement of
the previous proposition and let $\bold{G}$ be any connected Lie
group with Lie algebra $\frak{g}$. Denote by $\widetilde{\roman{Ad}}$
the complexification of the adjoint action of $\bold{G}$ in $\frak{g}$
and by \par\smallskip
\centerline{$\bold{N}_{\bold{G}}(\frak{a})=\{
g\in\bold{G}\, | \, \widetilde{\roman{Ad}}(g)(\frak a)\subset\frak a\}$}
\par\smallskip
\noindent the normalizer of a subspace $\frak{a}$ of
$\Hat{\frak{g}}$ in $\bold{G}$. Then
\form{\bold{N}_{\bold{G}}(\frak{q}_{\Psi})=
\bold{N}_{\bold{G}}
(\frak{q}_{\Phi})\,.}}\edef\corollariohe{\number\q.\number\x}
\dimo 
Indeed, if $g\in\bold{N}_{\bold{G}}(\frak{q}_{\Phi})$, then
$\widetilde{\roman{Ad}}(g)(\frak{q}_{\Psi})$ is still a parabolic
subalgebra of $\Hat{\frak{g}}$, minimal among those that are
contained in $\frak{q}_{\Phi}$ and satisfy
$\frak{q}\cap\bar\frak{q}=\frak{q}_{\Phi}\cap\bar\frak{q}_{\Phi}$.
Hence, by the uniqueness stated in Proposition {\teoremadc}, it
coincides with $\frak{q}_{\Psi}$ and therefore
$g\in \bold{N}_{\bold{G}}(\frak{q}_{\Psi})$.
Thus we proved the inclusion
$\bold{N}_{\bold{G}}(\frak{q}_{\Phi})\subset
\bold{N}_{\bold{G}}(\frak{q}_{\Psi})$.
\par
Note that $\widetilde{\roman{Ad}}$ can be considered as a homomorphism
$\bold{G}@>>>\bold{Int}(\Hat{\frak{g}})$, and
$\bold{N}_{\bold{G}}(\frak{a})=\widetilde{\roman{Ad}}^{-1}(
\bold{N}_{\bold{Int}(\Hat{\frak{g}})}(\frak{a})$ for every
subspace $\frak a$ of $\Hat{\frak{g}}$. Then the opposite inclusion
$\bold{N}_{\bold{G}}(\frak{q}_{\Psi})\subset
\bold{N}_{\bold{G}}(\frak{q}_{\Phi})$
follows from the fact that any parabolic subgroup of a connected
complex Lie group is the normalizer of its Lie algebra.
\enddimo
\medskip
\thm{Let $(\frak{g},\frak{q}_{\Phi})$ be an effective parabolic
minimal $CR$ algebra, satisfying condition $(F)$. Let $\Hat{\bold{G}}$
be any connected algebraic complex semisimple Lie group whose Lie algebra
is the complexification $\Hat{\frak g}$ of
$\frak{g}$ and let $\bold{Q}$ be its parabolic subgroup with
Lie algebra $\frak{q}$. Denote by $\bold{G}$ the analytic subgroup
of $\Hat{\bold{G}}$ with Lie algebra $\frak{g}$. Then\,:
\roster
\item"($i$)" $\bold{G}_+=\bold{G}\cap\bold{Q}$ is connected;
\item"($ii$)" $M=\bold{G}/\bold{G}_+$ is compact and simply connected.
\endroster}
\dimo
First we consider the case where $(\frak{g},\frak{q}_{\Phi})$ is
a totally real parabolic minimal $CR$ algebra satisfying condition
$(F)$. In this case $\frak{g}_+$ is a parabolic real Lie subalgebra
of $\frak{g}$. Denote by $\bold{N}$ the normalizer of
$\frak{g}_+$ in $\bold{G}$. Denote by $\Sigma$ the real root system
of $\frak{g}$ with respect to $\frak{h}^-$ and let
$\pi:\Cal{R} @>>>\Sigma$ denote the natural projection.
By \cite{Wi} we know that the fundamental
group of $\bold{G}/\bold{N}$ is a quotient group of the free
group generated by the elements of $\pi(\Phi)$ having multiplicity
$1$. But condition $(F)$ implies that no root in
$\pi(\Phi)$ has multiplicity $1$, and hence 
$\pi_1(\bold{G}/\bold{N})=\pmb{1}$. Thus $\bold{G}/\bold{N}$ is
simply connected. Since
$\bold{N}\supset\bold{G}_+$, we deduce that
$\bold{G}_+=\bold{N}$ is connected.
By Theorem \teoremadb, $M$ is also compact. Thus the proof is complete
in this case.\smallskip
Consider now a general  effective parabolic
minimal $(\frak{g},\frak{q}_{\Phi})$. Let $\Psi$ be the subset of
$\Cal{B}$ found in Proposition \teoremadc. Then we consider the parabolic
$\frak{q}_{\Pi}=\frak{q}_{\Psi}+\bar{\frak{q}}_{\Psi}$,
with $\Pi\subset\Psi$. 
\par
We apply Proposition {\proposizionedh}
to the $\frak{g}$-equivariant fibration 
$(\frak{g},\frak{q}_{\Psi})@>>>(\frak{g},\frak{q}_{\Pi})$.
The fiber is a parabolic minimal totally complex $CR$ algebra
$(\frak{g}',\frak{q}')$. Thus, with $\bold{G}_\Psi$ equal to
$\bold{N}_{\bold{G}}(\frak{q}_{\Psi})$, and
$\bold{G}_\Pi$ equal to 
$\bold{N}_{\bold{G}}(\frak{q}_{\Pi})=\bold{Q}_{\Pi}\cap\bold{G}$,
we obtain a $\bold{G}$-equivariant fibration
$\tilde{M}=\bold{G}/\bold{G}_\Psi@>>>\bold{G}/\bold{G}_\Pi$ whose
fiber is a complex flag manifold. Since both the base space
$\bold{G}/\bold{G}_\Pi$, by the first part of this proof,
and the fiber are simply connected, also
the total space $\tilde{M}$ is simply connected.
\par
By Corollary {\corollariohe}, we have 
$\bold{N}_{\bold{G}}(\frak{q}_{\Psi})=\bold{N}_{\bold{G}}(\frak{q}_{\Phi}) 
=\bold{G}\cap\bold{Q}_{\Phi}$, and therefore $M$ is diffeomorphic
to $\Tilde{M}$ 
and
therefore simply connected. 
By Theorem \teoremadb, $M$ 
is also compact.\enddimo\edef\teoremahf{\number\q.\number\x}
Theorem {\teoremahf} and Theorem {\teoremadb}
yield a rigidity theorem for
$CR$ manifolds\,:
\cor{Let $\bold{G}$ be a semisimple real Lie group and
$M$ a connected
$\bold{G}$-homogeneous $CR$ manifold. If the associated
$CR$ algebra $(\frak{g},\frak{q})$ is parabolic minimal and
has property $(F)$, then $M$ is simply connected and
$CR$-diffeomorphic to $M(\frak{g},\frak{q})$.\qed}
\edef\corollariohg{\number\q.\number\x}


%% file: 09min.tex
\se{Fundamental $CR$ algebras}
We give a criterion to read off the property of being fundamental
from the cross-marked Satake diagram\,:
\thm{An effective parabolic 
minimal $CR$ algebra $(\frak g,\frak q_{\Phi})$
is fundamental if and only if its corresponding 
cross-marked Satake diagram $(\pmb{\Cal S},\Phi)$
has the property:
\form{\alpha\in\Phi\setminus\Cal R_{\bullet}
 \Longrightarrow \epsilon_C(\alpha)\notin \Phi\, .}}
Here $\epsilon_C$ is the involution in $\Cal B(C)$ defined in 
Proposition \teoremaba. 
\smallskip\noindent
\dimo
Assume that $\alpha_1$ and $\alpha_2=\epsilon_C(\alpha_1)$ both
belong to $\Phi$, and let $\Psi=\{\alpha_1,\alpha_2\}$. 
Then $\Psi\subset\Phi$ and hence $\frak q_{\Phi}\subset\frak q_{\Psi}$.
To show that $(\frak g,\frak q_{\Phi})$ is not fundamental, it is
sufficient to check that $\frak q_{\Psi}=\overline{\frak q}_{\Psi}$.
To this aim it suffices to verify that 
$\Cal Q^n_{\Psi}=\overline{\Cal Q}^{\, n}_{\Psi}$.
Let $\Cal B(C)=\{\alpha_1,\alpha_2,\alpha_3,\hdots,\alpha_\ell\}$.
Every root $\alpha\in\Cal Q^n_{\Psi}$ can be written in the form
$\alpha=\sum_{i=1}^\ell{k_i\alpha_i}$ with $k_1+k_2>0$.
Since $C$ is adapted to the conjugation $\sigma$, 
using (\formulaej) we obtain\,:\par\smallskip
\centerline{
$\bar\alpha=\sum_{i=1}^\ell{k_i\epsilon_C(\alpha_i)}+
\sum_{\beta\in\Cal B_\bullet(C)}{k_{\alpha,\beta}\beta}
=\sum_{i=1}^\ell{k'_i\alpha_i}$,}\par\smallskip\noindent
with $k'_1+k'_2=k_2+k_1>0$, showing that also $\bar\alpha\in\Cal Q^n_{\Psi}$.
This shows that the condition is necessary.
\smallskip
Assume vice versa that there exists a 
proper parabolic subalgebra $\frak q'$ of
$\hat\frak g$ with $\frak q_{\Phi}\subset\frak q'=\overline{\frak q}\, '$.
Then $\frak q'=\frak q_{\Psi}$ for some $\Psi\subset\Phi$, 
$\Psi\neq\emptyset$. 
Since $\overline{\Cal Q}^{\, n}_{\Psi}=\Cal Q^n_{\Psi}\subset\Cal R^+(C)$,
we have $\Psi\cap\Cal R_{\bullet}=\emptyset$. Hence, again by
(\formulaej), we obtain that $\epsilon_C(\alpha)\in\Psi$ for
all $\alpha\in\Psi$. 
\enddimo\edef\teoremaea{{\number\q.\number\x}}
\cor{Fundamental effective parabolic minimal $CR$ algebras have 
the $(F)$ property.\qed}
\par
From Theorems {\teoremaea}, {\thmgef} and 
Proposition {\proposizionedg}
we obtain\,:
\thm{Let $(\frak g,\frak q_{\Phi})$ be an effective parabolic minimal
$CR$
algebra and let $(\pmb{\Cal S},\Phi)$ be its corresponding 
cross-marked Satake diagram.
Let \par\smallskip
\centerline{
$\Psi=\{\alpha\in\Phi\setminus\Cal{R}_{\bullet}\, | \, 
\epsilon_C(\alpha)\in\Phi\}$.}\par\smallskip
Then 
\roster
\item"($i$)" The diagram $\pmb{\Cal S}'$ obtained from $\pmb{\Cal S}$
by erasing all the nodes corresponding to the roots in $\Psi$ and
the lines and arrows issued from them is still a Satake diagram,
corresponding to a semisimple real Lie algebra $\frak g'$.
\item"($ii$)" $(\frak g,\frak q_{\Psi})$ is a totally real
effective parabolic minimal
$CR$ algebra.
\item"($iii$)" The natural map
$(\frak g,\frak q_{\Phi}) @>>> (\frak g,\frak q_{\Psi})$, defined by the
inclusion $\frak q_{\Phi}\subset\frak q_{\Psi}$, is a $\frak g$-equivariant
$CR$ fibration. The effective quotient of its fiber is the fundamental
parabolic minimal $CR$ algebra
$(\frak g'',\frak q_{\Phi'})$, associated to the cross-marked
Satake diagram $(\pmb{\Cal{S}}'',\Phi')$, 
where $\Phi'=\Phi\setminus\Psi$ and $\pmb{\Cal{S}}''$ is the union of
the $\sigma$-connected components of $\pmb{\Cal{S}}'$ that contain some 
root of $\Phi'$.\qed
\endroster}\edef\teofundamentalb{\number\q.\number\x}
\par
We call the map in $(iii)$ the {\it fundamental reduction} of
$(\frak g,\frak q_{\Phi})$ and the totally real
$CR$ algebra $(\frak g,\frak q_{\Psi})$ its {\it basis}.
\exam{Let $\frak{g}\simeq\frak{su}(2,2)$ and let
$\Phi=\{\alpha_2,\alpha_3\}$ (we refer to the diagram below). 
We have $\varepsilon_C(\alpha_i)=\alpha_{4-i}$ for $i=1,2,3$
and hence 
$\Psi=\{\alpha\in\Phi\, | \, \varepsilon_C(\alpha)\in\Phi\}=\{\alpha_2\}$. 
In particular $(\frak{g},\frak{q}_{\{\alpha_2,\alpha_3\}})$
is not fundamental. We obtain by Theorem~{\teofundamentalb}
a $\frak{g}$-equivariant
$CR$ fibration 
$(\frak{g},\frak{q}_{\{\alpha_2,\alpha_3\}})@>>>
(\frak{g},\frak{q}_{\{\alpha_2\}})$ with
fundamental fiber 
$(\frak{g}',\frak{q}'_{\{\alpha_3\}})$, with
$\frak{g}'\simeq\frak{sl}(2,\Bbb{C})$.
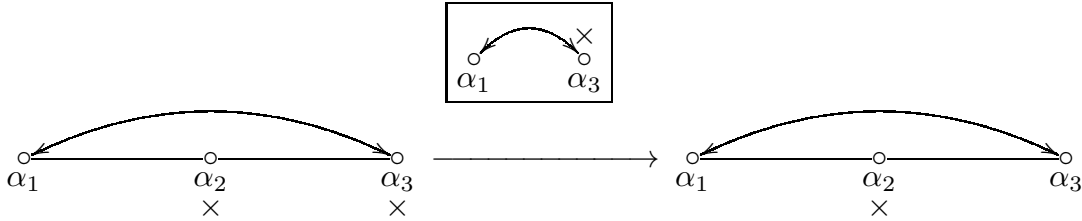
$$\CD
{\xymatrix@M=0pt @R=2pt{
\circ\ar @{-}[rr]\ar @/^15pt/@{<->}[rrrr]&&\circ\ar @{-}[rr]&&\circ\\
\alpha_1&&\alpha_2&&\alpha_3\\
&&\times&&\times}}@>{\matrix
\qquad\\
\boxed{\xymatrix@M=0pt @R=2pt{
\quad\\
\quad \\
&\times\\
\circ\ar @/^10pt/@{<->}[r]&\circ\\
\alpha_1&\alpha_3
}}
\quad\\
\\
\\
\endmatrix}>>
 {\xymatrix@M=0pt @R=2pt{
\circ\ar @{-}[rr]\ar @/^15pt/@{<->}[rrrr]&&\circ\ar @{-}[rr]&&\circ\\
\alpha_1&&\alpha_2&&\alpha_3\\
&&\times}}\endCD
$$
}

\cor{Let $\bold{G}$ be a semisimple Lie group and
$M$ a 
$\bold{G}$-homogeneous $CR$ manifold. Assume that
the $CR$ algebra
$(\frak{g},\frak{q}_{\Phi})$ associated to $M$ is
parabolic minimal. 
Let $(\frak{g},\frak{q}_{\Psi})$ 
be the basis of its fundamental reduction. Then there exists
a (totally real) 
$\bold{G}$-homogeneous $CR$ manifold $N$, with
associated $CR$ algebra $(\frak{g},\frak{q}_{\Psi})$, and a 
$\bold{G}$-equivariant submersion 
$\omega:M@>>>N$ such that the induced map
$\omega_*:\pi_1(M) @>>>\pi_1(N)$ is an isomorphism.}
\dimo
Let $o$ be a point of $M$ and let $\bold{G}_+$ the stabilizer
of $o$ in $\bold{G}$. Let $\bold{H}$ be the analytic subgroup
of $\bold{G}$ generated by $\frak{g}\cap\frak{q}_{\Psi}$. Then
$\bold{H}$ contains $\bold{G}_+^{\circ}$. We claim that
$\bold{H}\cdot\bold{G_+}=\bold{G}'_+$ is a Lie subgroup of
$\bold{G}$. Indeed, for all $g\in\bold{G}_+$, we have
$\roman{Ad}(g)(\frak{q}_{\Phi})=\frak{q}_{\Phi}$. 
Since $g$ is real, we also have
$\roman{Ad}(g)(\bar\frak{q}_{\Phi})=\bar\frak{q}_{\Phi}$ and therefore
$\roman{Ad}(g)(\frak{q}_{\Psi})=\frak{q}_{\Psi}$
because $\frak{q}_{\Psi}$ is generated by 
$\frak{q}_{\Phi}+\bar\frak{q}_{\Phi}$.
This implies that $\roman{ad}(g)(\bold{H})=\bold{H}$ for all
$g\in\bold{G}_+$, and hence $\bold{G}'_+$ is a subgroup of
$\bold{G}$. It is a Lie subgroup because its Lie algebra is
real parabolic. Then $N=\bold{G}/\bold{G}'_+$ is a 
$\bold{G}$-homogeneous manifold. By the inclusion
$\bold{G}_+\subset\bold{G}_+'$ we obtain a 
$\bold{G}$-equivariant submersion
$\omega:M@>>>N$. By construction the fiber is 
connected. It has a natural structure of 
$CR$ manifold, associated to a fundamental $CR$ algebra
$(\frak{g}'',\frak{q}_{\Phi'})$,
as in Theorem {\teofundamentalb},  which is parabolic and minimal. 
By Corollary {\corollariohg} the fiber is simply connected. Hence 
$\omega_*:\pi_1(M)@>>>\pi_1(N)$ is an isomorphism.
\enddimo

%% file: 10min.tex
\se{Totally real and totally complex CR algebras}
From the discussion in the previous section we obtain the criterion\,:
\thm{A simple effective parabolic 
minimal $CR$ algebra $(\frak g,\frak q_{\Phi})$,
with corresponding cross-marked Satake diagram
$(\pmb{\Cal{S}},\Phi)$,  is 
totally real if and only if the following conditions hold true:
\roster
\item "($i$)" $\Phi\cap\Cal R_\bullet=\emptyset$;
\item "($ii$)"
$\epsilon_C(\Phi)=\Phi$.\qed
\endroster}
\par
\thm{A simple effective parabolic 
minimal $CR$ algebra $(\frak g,\frak q_{\Phi})$ 
with associated cross-marked Satake diagram $(\pmb{\Cal{S}},\Phi)$
is 
totally complex if and only if either:
\roster
\item "($i$)"
$\frak g$ is compact, or 
\item "($ii$)"
$\frak g$ is of the complex type and all 
cross-marked nodes are in the 
same connected component of $\pmb{\Cal{S}}$, or
\item "($iii$)"
$(\pmb{\Cal{S}},\Phi)$ is
one of the
following:
$$\align
&\cases
\Phi=\{\alpha_1\}\\
\Phi=\{\alpha_{\ell}\}
\endcases
\tag{$\roman{A\,II}$} \\\vspace{2\jot}
&\cases
\Phi=\{\alpha_{\ell}\}\\
\Phi=\{\alpha_{\ell-1}\}
\endcases
\tag{$\roman{D\,II}$}
\endalign
$$
\endroster
}\edef\thmtotcplx{\number\q.\number\x}
\dimo
The $CR$ algebra $(\frak g, \frak q)$ is totally complex if and only if
$\frak g+\frak q=\hat\frak g$. This is equivalent to the fact that the
standard $CR$ manifold $\bold G\cdot o$ is open in the complex flag
manifold $\hat\bold G/\bold Q$. Since it is also closed, it follows that
$\bold G$ is transitive on $\hat\bold G/\bold Q$. The result then follows
from \cite{W2, Corollary 1.7}.
\enddimo
This yields also a characterization of ideal nondegenerate parabolic
minimal $CR$ algebras. Indeed, for general parabolic
$CR$ algebras, we have\,:
\thm{Let $(\frak g,\frak q)$ be a simple effective parabolic $CR$ algebra.
Then $(\frak g,\frak q)$ is either totally
complex or ideal nondegenerate.
}
\dimo
We recall (see \cite{MeN5}) that an effective $CR$ algebra
$(\frak g,\frak q)$ is ideal nondegenerate if 
$\Cal H_+=\left(\frak q + \bar\frak q\right)\cap\frak g$ does not
contain a non zero ideal of $\frak g$. When $\frak g$ is simple,
this is equivalent to the fact that $\Cal H_+\neq\frak g$, i.e.
that $(\frak g,\frak q)$ is not totally complex.
\enddimo

%% file: 11min.tex
\se{Weak nondegeneracy}
In this section we characterize those parabolic minimal
$(\frak{g},\frak{q}_{\Phi})$ that are weakly nondegenerate.
We recall from Proposition {\proposizioneda} that this means
that there is no nontrivial complex $CR$ fibration 
$M(\frak{g},\frak{q}_{\Phi})@>>>N$ with totally complex fibers.
In turns this is equivalent to the fact that 
$M(\frak{g},\frak{q}_{\Phi})$ is not, locally, $CR$ equivalent
to the product of a $CR$ manifold with the same $CR$ codimension
and of a complex manifold of positive dimension.
\smallskip
From  Proposition {\proposizionedh} we obtain\,:
\lem{A fundamental effective parabolic
minimal $CR$ algebra $(\frak g, \frak q_\Phi)$ is weakly
degenerate if and only if there is  $\Psi\subsetneq\Phi$ such that
the $\frak g$-equivariant
fibration $(\frak g, \frak q_\Phi)@>>>(\frak g,\frak q_\Psi)$
is a $CR$ fibration with totally complex fiber. \qed
}
\lem{Let $(\frak g,\frak q_\Phi)$ be a minimal fundamental
effective parabolic $CR$ algebra. 
A necessary and sufficient 
condition in order that
$(\frak g,\frak q_\Phi)$ be weakly degenerate is that there exists
$\Psi\subset\Phi$ satisfying conditions in Theorem {\thmcrfibr}
and such
that $\frak q_\Psi\subset\frak q_\Phi+\bar\frak q_\Phi$.\qed}
\par\edef\lemhb{\number\q.\number\x}
We now give a characterization of the pairs $(\Phi,\Psi)$ for
which (\formulafiba) is a $CR$ fibration with totally complex
fiber in terms of properties of the roots $\alpha$ in
$\Phi\setminus\Psi$. 
\lem{Let $(\frak g,\frak q_\Phi)$ be a minimal fundamental
effective parabolic $CR$ algebra, with $\frak g$ of the real type
(i.e.\ $\hat\frak g$ is also simple).
Let $\emptyset\ne\Psi\subset\Phi$ and assume that $(\formulafiba)$
is a $CR$ fibration.
Then for each $\alpha\in\Phi\setminus\Psi$ we have the following
possibilities:
\roster
\item"$(i)$" $\check\Psi(\alpha)\subset\Cal R_\bullet$;
\item"$(ii)$" 
$(\roman a)$ $\check\Psi(\alpha)\cap\Cal R_\bullet=\emptyset$ and
$(\roman b)$ $\big(\check\Psi(\alpha)\big)\cap\varepsilon_C
\big(\check\Psi(\alpha)\big)=\emptyset$;
\item"$(iii)$" 
$(\roman a)$ $\emptyset\neq\check\Psi(\alpha)
\cap\Cal R_\bullet\neq\check\Psi(\alpha)$ and
$(\roman a)$ $\varepsilon_C\big(\check\Psi(\alpha)\setminus\Cal 
R_\bullet\big)=\check\Psi(\alpha)\setminus
\Cal R_\bullet$.
\endroster}\edef\lemwnd{\number\q.\number\x}
\dimo
Fix $\alpha\in\Phi\setminus\Psi$ with
$\check\Psi(\alpha)\not\subset\Cal R_\bullet$ and let 
$\delta=\delta\big(\check\Psi(\alpha)\big)$.
\par
If $\beta\in\check\Psi(\alpha)\setminus\Cal R_\bullet$ and 
$\varepsilon_C(\beta)\in\check\Psi(\alpha)$, then $\supp(\bar\delta)
\ni\varepsilon_C(\beta)$. Since it is connected and does not
meet $\Psi$, we obtain $\supp(\bar\delta)\subset\check\Psi(\alpha)$.
This implies that $\check\Psi(\alpha)\setminus\Cal 
R_\bullet=\varepsilon_C\big(\check\Psi(\alpha)\setminus
\Cal R_\bullet\big)$. 
In this way we have shown that either 
$\varepsilon_C\big(\check\Psi(\alpha)\setminus\Cal 
R_\bullet\big)\cap\check\Psi(\alpha)\setminus\Cal R_\bullet=
\emptyset$ 
or $\varepsilon_C\big(\check\Psi(\alpha)\setminus\Cal 
R_\bullet\big)=\check\Psi(\alpha)\setminus\Cal R_\bullet$ 
\par
If $\check\Psi(\alpha)\cap\Cal R_\bullet$ is not empty, then 
there exists $\beta\in\check\Psi(\alpha)\setminus\Cal R_\bullet$
such that $\partial_e\{\beta\}\cap\check\Psi(\alpha)
\cap\Cal R_\bullet\neq\emptyset$. Hence $\supp(\bar\beta)\cap
\check\Psi(\alpha)\neq\emptyset$ and, by the same argument as above,
$\varepsilon_C(\beta)\in\check\Psi(\alpha)$ and 
we get $(iii.\roman b)$.
\par
Finally we consider the case where $\check\Psi(\alpha)
\cap\Cal R_\bullet=\emptyset$. 
The boundary 
$\partial_e\big(\check\Psi(\alpha)\big)$ is not empty, thus it
contains a root $\beta\in\Psi$ and $\beta \not\in\Cal R_\bullet$ because
of Theorem {\thmcrfibr}. 
The fact that $\frak g,\frak q_\Psi$ is fundamental
implies that $\varepsilon_C(\beta)\not\in\Psi$. 
In particular $\varepsilon_C(\beta)\not\in
\partial_e\big(\check\Psi(\alpha)\big)$. 
Applying again Theorem {\thmcrfibr} we have 
$\varepsilon_C\big(\check\Psi(\alpha)\setminus\Cal R_\bullet\big)
\cap\Psi=\emptyset$, hence $\varepsilon_C(\beta)\not\in\check\Psi(\alpha)$.
Since $\varepsilon_C(\beta)\in\partial_e\big(\supp(\bar\delta)\big)$ and
$\supp(\bar\delta)\cap\Cal R_\bullet\not\subset
\check\Psi(\alpha)$, it follows that 
$\supp(\bar\delta)\cap\check\Psi(\alpha)=\emptyset$, thus
$\check\Psi(\alpha)\cap\varepsilon_C\big(\check\Psi(\alpha)\big)=\emptyset$.
\enddimo
\par
\lem{With the same hypotheses of  Lemma {\lemwnd}, 
the effective quotient of the fiber of the $\frak g$-equivariant 
$CR$ fibration $(\frak g, \frak q_\Phi)@>>>
(\frak g,\frak q_\Psi)$ has cross-marked Satake diagram
$$
\pmb{\Cal S}'=\bigcup_{\alpha\in\Phi\setminus\Psi}\check\Psi(\alpha)
\cup\varepsilon_C\big(\check\Psi(\alpha)\setminus\Cal R_\bullet\big)
$$
and $\Phi'=\Phi\cap\pmb{\Cal S}'$.
\par
In particular it is totally complex if and only if 
for each $\alpha\in\Phi\setminus\Psi$, condition $(i)$ 
or $(ii)$ of Lemma {\lemwnd} holds.}
\dimo
The effective quotient  is  described 
in \cite{MeN5} and at the end of
\S\secta .
From Theorem {\thmgef} we know that
$\pmb{\Cal S}'\subset
\bigcup_{\alpha\in\Phi\setminus\Psi}
\check\Psi(\alpha)\cup\varepsilon_C
\big(\check\Psi(\alpha)\setminus\Cal R_\bullet\big)$.
Equality then follows from the
observation that if $\beta\in\check\Psi(\alpha)
\setminus\Cal R_\bullet$ then $\supp(\bar\beta)\cap\Psi=\emptyset$.
\par
To prove the second statement, we can assume that there exists 
exactly one root $\alpha\in\Phi\setminus\Psi$.
In cases $(i)$ and $(ii)$
of Lemma {\lemwnd} the cross-marked Satake diagram of the fiber
is of the types described in Theorem {\thmtotcplx} $(i),(ii)$ and
is totally complex.
If we are in case $(iii)$ of Lemma {\lemwnd}, then 
$\check\Psi(\alpha)\cap\Cal R_\bullet\ne\emptyset$, and 
the fiber is totally complex if and only if $\big(\check\Psi(\alpha),
\Phi\cap\check\Psi(\alpha)\big)$ is one of the diagrams in 
Theorem {\thmtotcplx} $(iii)$.
\par
Since $\partial_e\big(\check\Psi(\alpha)\big)
\cap\Cal R_\bullet=\emptyset$ and 
$\varepsilon_C\big(\partial_e\check\Psi(\alpha)\big)\cap
\big(\check\Psi(\alpha)\cup\partial_e\check\Psi(\alpha)\big)=\emptyset$,
We have that $\varepsilon_C$ is not the identity, hence $\pmb{\Cal S}$ must
be of type $\roman{A\,III}$, 
$\roman{A\,IV}$, $\roman{D\,Ib}$, 
$\roman{D\,IIIb}$, $\roman{E\,II}$ or $\roman{E\,III}$.
We exclude types $\roman{A\,III}$, $\roman{A\,IV}$, $\roman{D\,Ib}$ and
$\roman{E\,II}$
because they do not contain 
subdiagrams of type $\roman{A\,II}$ or 
$\roman{D\,II}$, so we are left with types 
$\roman{D\,IIIb}$ and $\roman{E\,III}$.
\par
Type $\roman{D\,IIIb}$
must be excluded because in this case we have $\alpha=\alpha_1$
or $\alpha_{\ell-2}$, $\check\Psi(\alpha)=\{\alpha_1,\dots,\alpha_{\ell-2}\}$
and $\partial_e\big(\check\Psi(\alpha)\big)=\{\alpha_{\ell-1},
\alpha_{\ell}\}=\varepsilon_C\big(\partial_e\check\Psi(\alpha)\big)$.
\par
Similarly type $\roman{E\,III}$
must be excluded because we have $\alpha=\alpha_3$
or $\alpha_{5}$,
$\check\Psi(\alpha)=\{\alpha_2,\alpha_3,\alpha_4,\alpha_5\}$ and 
$\partial_e\big(\check\Psi(\alpha)\big)=\{\alpha_1,
\alpha_6\}=\varepsilon_C\big(\partial_e\check\Psi(\alpha)\big)$.
\enddimo\edef\lemhd{\number\q.\number\x}
\par

\thm{Let $(\frak g,\frak q_\Phi)$ be a simple fundamental
effective parabolic minimal $CR$ algebra and assume that it is not
totally complex.
Let $\Pi$ be the set of simple roots $\alpha$ in $\Phi$ that
satisfy either one of\,:
\roster
\item"$(i)$" $\check\Phi(\alpha)\subset\Cal R_\bullet$;
\item"$(ii)$" $\big(\check\Phi(\alpha)\cup
\partial_e\check\Phi(\alpha)\big)\cap\Cal R_\bullet=\emptyset$ 
and $\varepsilon_C\big(\check\Phi(\alpha)\big)\cap\Phi=\emptyset$.
\endroster
Then $(\frak g,\frak q_\Phi)$ is weakly nondegenerate if 
and only if $\Pi=\emptyset$.
\par
Set $\Psi=\Phi\setminus\Pi$.
Then $(\frak g,\frak q_\Phi)@>>>(\frak g,\frak q_\Psi)$ is a 
$\frak g$-equivariant $CR$ fibration with totally complex fiber
and fundamental weakly nondegenerate base.}
\dimo
Fix $\alpha\in\Phi\setminus\Psi$. Then conditions $(i)$ and $(ii)$ are
necessary and sufficient for $(\frak g,\frak q_\Phi)@>>>
(\frak g,\frak q_{\Phi\setminus\{\alpha\}})$ to be a 
$\frak g$-equivariant $CR$ fibration with totally complex fiber. This
observation, Lemma {\lemhd} and Lemma {\lemhb}
yield our first statement. 
\par
To prove the last part of the Theorem, we make the following
\par\smallskip
\noindent{\it Claim}. {\sl Let $\alpha,\beta\in\Phi$ with $\alpha\in\Pi$.
Then $\beta$ satisfies either $(i)$ 
or $(ii)$ for $\Phi$ 
if and only if $\beta$ satisfies $(i)$ or $(ii)$ for
$\Phi'=\Phi\setminus\{\alpha\}$.}
\par\smallskip
Assuming that this claim is true, we conclude as follows.
If $\Pi=\{\beta_1,\dots,\beta_k\}$, we have $\frak g$-equivariant 
$CR$ fibrations with totally complex fibers:
$$
(\frak g,\frak q_\Phi)@>>>
(\frak g,\frak q_{\Phi\setminus\{\beta_1\}})@>>>
(\frak g,\frak q_{\Phi\setminus\{\beta_1,\beta_2\}})@>>>\cdots@>>>
(\frak g,\frak q_{\Phi\setminus\Pi}).
$$
Their composition is still a $\frak g$-equivariant 
$CR$ fibration with totally complex fiber, and the base
$(\frak g,\frak q_\Psi)$ is weakly nondegenerate.
\par
Now we prove the claim.
If $\beta\not\in\partial_e\big(\check\Phi(\alpha)\big)\cup
\partial_e\varepsilon_C\big(\check\Phi(\alpha)\big)$,
then $\check\Phi(\beta)=\check\Phi'(\beta)$, 
$\varepsilon_C\big(\check\Phi(\beta)\big)
=\varepsilon_C\big(\check\Phi'(\beta)\big)$, and there is 
nothing to prove.
\par
Assume $\beta\in\partial_e\big(\check\Phi(\alpha)\big)$; 
then $\check\Phi'(\beta)=\check\Phi(\beta)\cup\check\Phi(\alpha)$.
If $\check\Phi(\alpha)\subset\Cal R_\bullet$, then 
$\check\Phi(\beta)\subset\Cal R_\bullet$ if and only if
$\check\Phi'(\beta)\subset\Cal R_\bullet$.
\par
If $\check\Phi'(\beta)\cap\Cal R_\bullet=\emptyset$, we need to 
prove that, if $\beta$ satisfies $(i)$ or $(ii)$, then 
$\check\Phi(\alpha)\cap\varepsilon_C\big(\check\Phi(\beta)\big)
=\emptyset$. This is true because otherwise $\varepsilon_C(\beta)
\in\partial_e\big(\check\Phi(\alpha)\big)$, and this yields
a contradiction because we assumed that
$(\frak{g},\frak{q}_{\Phi})$ is fundamental.
\par
Finally if $\beta\in\partial_e\varepsilon_C
\big(\check\Phi(\alpha)\big)$ then 
$\beta\not\in \Cal R_\bullet$ and $\varepsilon_C(\beta)
\in\partial_e\big(\check\Phi(\alpha)\big)$, again
contradicting the assumption that 
$(\frak{g},\frak{q}_{\Phi})$ is fundamental.
\enddimo

%% file: 12min.tex
\se{Strict nondegeneracy} 
In this section we give necessary and sufficient conditions for a weakly
nondegenerate $CR$ algebra to be strictly nondegenerate. 
We recall from the introduction that 
the $CR$ geometry of strict nondegenerate 
homogeneous $CR$ manifolds can be related to the {\it standard}
models and investigated by using the {\it Levi-Tanaka algebras}
(cf. \cite{MeN1}, \cite{T1}, \cite{T2}). Therefore, by classifying
the weakly degenerate minimal orbits that 
do not have the strict nondegeneracy property,
we single out a class of homogeneous $CR$ manifolds with a
highly non trivial $CR$ structure that
cannot be discussed by using the
standard Levi-Tanaka models. This also explains the need to introduce
$CR$ algebras, as a generalization of the Levi-Tanaka algebras,
in \cite{MeN5}.
\smallskip
First we  
reformulate weak and strict  nondegeneracy in terms of the root system:
\lem{A fundamental effective parabolic minimal CR algebra $(\frak g,\frak
q)$ is weakly nondegenerate if and only if for every root
$\alpha\in\bar\Cal Q\setminus\Cal Q$ there exist 
a sequence  
$\left(\beta_i\in\Cal Q\right)_{1\leq i\leq n}$ such
that 
\form{\alpha_j=\alpha+\sum_{i\leq j}{\beta_i}\in\Cal R\quad
\forall j=1,\hdots,n\,,\qquad \alpha_n\notin \Cal{Q}\cup\bar{\Cal{Q}}
\, .}}
\edef\lemaltwnd{\number\q.\number\x}
\dimo The statement \edef\formuladeibeta{\number\q.\number\t}
is an easy consequence of
\cite{MeN5, Theorem 6.2}. 
\enddimo \hyphenation{str-ong-ly}
Likewise we have\,:
\lem{A fundamental effective parabolic minimal $CR$ algebra 
$(\frak g,\frak q)$ is strictly nondegenerate if and only if for every root
$\alpha\in\bar\Cal Q\setminus\Cal Q$ there exists 
$\beta\in\Cal{Q}$ such that
$\alpha+\beta\in\Cal R$ and $\alpha+\beta\not\in\Cal
Q\cup\bar\Cal Q$.\qed}
Next we prove that it suffices to check this condition on
purely imaginary roots\,:
\prop{A necessary and sufficient condition for a fundamental effective weakly
nondegenerate parabolic minimal $CR$ algebra $(\frak g,\frak q)$ 
to be strictly nondegenerate is that for every root
$\alpha\in\Cal R_\bullet\cap\bar\Cal Q\setminus\Cal Q$ there exists $\beta\in
\Cal{Q}$ such that 
$\alpha+\beta\in\Cal R$ and $\alpha+\beta\not\in\Cal
Q\cup\bar\Cal Q$.}\edef\propsnd{\number\q.\number\x}
\dimo
The condition is obviously necessary.
\par
To prove sufficiency, consider a root $\alpha\in\bar\Cal Q\setminus\Cal
Q$, $\alpha\not\in\Cal{R}_\bullet$; since $\alpha\prec 0$,
we have $-\alpha\in\Cal{R}^+\setminus\Cal{R}_\bullet$. This implies
that $-\bar\alpha\in\Cal{R}^+\subset\Cal{Q}$. Then $-\alpha\in\bar{\Cal{Q}}$
and $\alpha\in\bar\Cal Q^r$. 
By the assumption that
$(\frak g,\frak q)$ is weakly
nondegenerate, using Lemma {\lemaltwnd} we can find
a sequence of roots
$(\beta_i)$ 
satisfying (\formuladeibeta). 
Take the sequence $(\beta_i)_{1\leq i\leq n}$  
of minimal length; we claim that
for every permutation $\tau$ of the indices, 
the sequence $(\beta_{\tau(i)})_{1\leq i\leq n}$ still
satisfies (\formuladeibeta). 
\par
Indeed, fix a Chevalley basis
$\{X_{\alpha}\}_{\alpha\in\Cal R}$. Then, for every transposition $(i,i+1)$:
$$
\align
\frak{q}+\bar{\frak{q}}&\not\ni
[X_{\beta_n},\hdots, X_{\beta_{i+1}}, X_{\beta_{i}},\hdots,
X_{\beta_1},X_{\alpha}] =\\
&= [X_{\beta_n},\hdots, X_{\beta_{i}}, X_{\beta_{i+1}},\hdots,
X_{\beta_1},X_{\alpha}] +
[X_{\beta_n},\hdots, [X_{\beta_{i+1}},X_{\beta_{i}}],\hdots,
X_{\beta_1},X_{\alpha}].
\endalign
$$
The last addendum in the right hand side belongs to 
$\frak{q}+\bar{\frak{q}}$ by our assumption that $(\beta_i)_{1\leq i\leq n}$
has minimal length. Thus
 $$ 
 [X_{\beta_n},\hdots\,,X_{\beta_{i}}, X_{\beta_{i+1}},\hdots\,,
X_{\beta_1},X_{\alpha}]\in \Hat{\frak{g}}\setminus\left(
\frak{q}+\bar{\frak{q}}\right).
$$
In particular $\alpha+\beta_i\in\Cal R$ for
every $i$.
At least one of the $\beta_i$'s, say $\beta_{i_0}$, does not belong to
$\bar\Cal Q$, so $\alpha+\beta_{i_0}\not\in\bar\Cal Q$. 
Indeed, since $\alpha\in\bar{\Cal{Q}}^r$, 
if $\alpha+\beta_i\in\bar{\Cal{Q}}$, then also
$\beta_i=(\alpha+\beta_i)+(-\alpha)\in \bar{\Cal{Q}}$. 
By a permutation, we can take $\beta_{i_0}=\beta_n$. 
Then we claim that $\alpha+\beta_n\notin\Cal{Q}\cup\bar{\Cal{Q}}$.
Indeed we already choose $\beta_n$ so that $\alpha+\beta_n\notin
\bar{\Cal{Q}}$. If $\alpha+\beta_n\in\Cal{Q}$, we have
$[X_{\beta_1},\hdots,X_{\beta_{n-1}},X_{\beta_n},X_\alpha]=
[X_{\beta_1},\hdots,X_{\beta_{n-1}},[X_{\beta_n},X_\alpha]]\in
\frak{q}$, because $X_{\beta_i}\in\frak{q}$ for every $i=1,\hdots,n$,
and hence $\alpha_n\in\Cal{Q}$, contradicting
(\formuladeibeta).
\enddimo
\thm{Let $(\frak g,\frak q_{\Phi})$ be an effective parabolic
minimal $CR$ algebra, with $\frak{g}$ simple. If 
$(\frak g,\frak q_{\Phi})$ is weakly nondegenerate, but 
is not strictly nondegenerate,
then $\Phi$ 
is contained in a
connected component of $\Cal{B}\cap\Cal R_\bullet$.
\par
The strictly nondegenerate $(\frak g,\frak q_{\Phi})$ with
$\frak{g}$ simple and $\Phi\subset\Cal{R}_{\bullet}$ are those
listed below\,:
$$\align\TagsAsText
\Phi=&\{\alpha_{p+1}\} \tag\text{$\roman{B\,Ib\ /\ B\,II}$}\\\vspace{2\jot}
\Phi=&\{\alpha_{2i-1}\}\,,\; 1\leq i\leq p \tag{$\roman{C\,IIa\ /\ 
IIb}$}\\\vspace{2\jot}
\Phi=&\{\alpha_{p+1}\} \tag{$\roman{D\,Ia}$}\\\vspace{2\jot}
\Phi=&\{\alpha_{2}\} \tag{$\roman{D\,II}$}
\endalign
$$
$$\align
\Phi=&\cases
\{\alpha_4\}\\
\{\alpha_3,\alpha_4\}\\
\{\alpha_4,\alpha_5\}\\
\{\alpha_3,\alpha_4,\alpha_5\}
\endcases
\tag{$\roman{E\,III}$}\\\vspace{2\jot}
\Phi=&\cases
\{\alpha_3\}\\
\{\alpha_5\}\\
\{\alpha_3,\alpha_4\}\\
\{\alpha_3,\alpha_5\}\\
\{\alpha_4,\alpha_5\}\\
\{\alpha_2,\alpha_3\}\\
\{\alpha_2,\alpha_5\}
\endcases
\tag"{$\left(\matrix
\roman{E\,IV}\\
\roman{E\,VII}\\
\roman{E\,IX}
\endmatrix\right)$}"\\\vspace{2\jot}
\Phi=&\cases
\{\alpha_2\}\\
\{\alpha_4\}
\endcases
\tag{$\roman{F\,II}$}
\endalign
$$}
\dimo 
We prove the first statement.
The proof of the second will be omitted, as it requires a straightforward 
case by case analysis, chasing over the different Satake diagrams.\par
Suppose that $(\frak g, \frak q_{\Phi})$ is 
weakly, but not strictly, nondegenerate. Then there is some root
$\alpha\in\bar\Cal{Q}_{\Phi}\setminus\Cal{Q}_{\Phi}$, $\alpha\prec 0$, 
such that
$\alpha+\beta\in\Cal{Q}_{\Phi} \cup\bar\Cal{Q}_{\Phi}$ for all 
$\beta\in\Cal{Q}_{\Phi}$ for which $\alpha+\beta\in\Cal{R}$. 
By
Proposition {\propsnd} we can take
$\alpha\in\Cal R_\bullet$. 
Let $\Cal B'$ be the connected component of $\supp (\alpha)$ in
$\Cal{B}\cap\Cal R_\bullet$. Since
$\alpha\not\in\Cal{Q}_{\Phi}$, we have $\Cal{B}'\cap\Phi\neq\emptyset$.\par
Since we assumed that $(\frak{g},\frak{q}_{\Phi})$ is weakly
nondegenerate, for each $\gamma\in\Phi$ the set
$\check{\Phi}(\gamma)$ is not contained in $\Cal{R}_{\bullet}$.
As $\supp(\alpha)\cap\Phi\neq\emptyset$, this implies that
there is some $\beta\in\Cal{Q}_{\Phi}$, with $\beta\prec 0$,
such that $\beta\notin\Cal{R}_{\bullet}$ and 
$\alpha+\beta\in\Cal{R}$. Since
$\beta\in\Cal{Q}_{\Phi}^r$ and $-\alpha\in\Cal{Q}_{\Phi}^n$,
we obtain that
$\alpha+\beta\notin\Cal{Q}_{\Phi}$. If $\Cal{B}'\cap\Phi$
contains some $\alpha_i$ which does not belong to $\supp(\alpha)$,
this $\alpha_i$ would belong to $\supp(\overline{\alpha+\beta})$.
Indeed $\alpha+\beta\notin\Cal{R}_{\bullet}$, hence
$\supp(\overline{\alpha+\beta})$ contains all simple 
imaginary roots
$\gamma$ that are not in $\supp(\alpha+\beta)$ and such that
$\partial_e\check{\Xi}(\gamma)\cap\supp(\alpha+\beta)\neq\emptyset$.
This shows that $\Cal{B}'\cap\Phi=\supp(\alpha)\cap\Phi$. 
\par
Let $\Cal A=\left(\Cal{R}_{\bullet}\cap\left[\Phi\setminus\Cal{B}'
\right]\right)
\cup \varepsilon_C(\Phi\setminus\Cal{R}_\bullet)$.
We want to show that $\Cal{A}=\emptyset$.\par
Assume by contradiction that
$\Cal A$ is not empty. Then there exists a
segment $S$ in $\Cal B\setminus\Phi$ joining
$\Cal{A}$ to $\supp(\alpha)$, i.e. 
such that $\partial_eS\cap\Cal
A\neq\emptyset$, $\partial_eS\cap\supp(\alpha)\neq\emptyset$.
By taking $S$ of minimal lenght,
we can also assume that $S\cap (\Cal A\cup\supp(\alpha))=\emptyset$\,.
\par
Let $\beta=-\delta(S)$. Then $\beta\prec 0$, $\beta\in\Cal{Q}_{\Phi}^r$ and
$\beta\not\in\Cal R_\bullet$, so that $\alpha+\beta
\in\Cal{R}\setminus\Cal{Q}_{\Phi}$\,.
\par
If there is some $\alpha_i$ in  
 $\partial_eS\cap\Cal A\cap\Cal R_\bullet\neq\emptyset$, then
$\alpha+\beta\in\Cal{R}$, 
$\supp(\overline{\alpha+\beta})\ni\alpha_i$, and $\alpha+\beta\not\in\bar
\Cal{Q}_{\Phi}$, contradicting our assumption.
\par
If $\partial_eS\cap\Cal A\cap\Cal R_\bullet=\emptyset$, there
is $\alpha_i$ in $\Phi\setminus\Cal{R}_{\bullet}$ with
$\varepsilon_C(\alpha_i)\in\partial_e S\cap\Cal{A}$. Set
$\beta'=\beta-\epsilon_C(\alpha_i)$. Then $\beta'\in\Cal{Q}_{\Phi}$,
and
$\alpha+\beta'\in\Cal{R}\setminus\left(\Cal{Q}_{\Phi}\cup\bar\Cal{Q}_{\Phi}
\right)$, 
yielding a contradiction; this shows that $\Cal A$ is empty,
completing the proof of our first claim.
\enddimo
\medskip

%% file: 13min.tex
\se{Essential pseudoconcavity for minimal orbits}
Let $(M,HM,J)$ be a $CR$ manifold of finite kind (cf. \S\sectionb).
We say that $(M,HM,J)$ is 
{\it essentially pseudoconcave} (see \cite{HN2}) if it is possible
to define a Hermitian symmetric smooth scalar product 
$h$ on the
fibers of $HM$ such that for each $\xi\in H^0M$ the Levi form
$\Cal{L}_{\xi}$ has zero trace with respect to $h$.
For a homogeneous $CR$ manifold, this last condition is equivalent to
the fact that for each $\xi\in H^0M$ the Levi form
$\Cal{L}_{\xi}$ is either $0$ or has at least one positive and one
negative eigenvalue. \par
The $CR$ functions defined on 
essentially pseudoconcave $CR$ manifolds enjoy 
some nice properties, like local smoothness and the
local maximum modulus principle;
$CR$ sections of $CR$ complex line bundles
have the weak unique continuation property (see \cite{HN2}, \cite{HN3}).
When $M$ is compact and essentially pseudoconcave,
global $CR$ functions are constant
and $CR$-meromorphic functions form a field of finite transcendence degree
(see \cite{HN4}).
\par
In this section we classify the essentially pseudoconcave
minimal orbits of complex flag manifolds.
\par
We keep the notation of the previous sections. In particular,
$(\frak{g},\frak{q}_{\Phi})$ is an effective parabolic minimal
$CR$ algebra, with associated cross-marked Satake diagram
$(\pmb{\Cal{S}},\Phi)$. Moreover, we introduce a Chevalley system
for $(\Hat{\frak{g}},\Hat{\frak{h}})$, i.e. a family
$(Z_{\alpha})_{\alpha\in\Cal{R}}$ with the properties
(\cite{B2, Ch.VIII,\S{2}})\,: 
\roster
\item"($i$)" $Z_{\alpha}\in\Hat{\frak{g}}^{\alpha}$ for all
$\alpha\in\Cal{R}$;
\item"($ii$)" $[Z_{\alpha},Z_{-\alpha}]=-H_{\alpha}$, 
where $H_\alpha$ is the unique element of 
$[\Hat{\frak{g}}^{\alpha},\Hat{\frak{g}}^{-\alpha}]$ such that
$\alpha(H_{\alpha})=2$;
\item"($iii$)" the $\Bbb{C}$-linear map that transforms
each $H\in\Hat{\frak{h}}$ into $-H$ and $Z_{\alpha}$ into
$Z_{-\alpha}$ for every $\alpha\in\Cal{R}$ is an automorphism
of the complex Lie algebra $\Hat{\frak{g}}$.
\endroster
In particular,
$(Z_{\alpha})_{\alpha\in\Cal{R}}\cup (H_{\alpha})_{\alpha\in\Cal{B}}$
is a basis of $\Hat{\frak{g}}$ as a $\Bbb{C}$-linear space. We denote
by $(\xi^{\alpha})_{\alpha\in\Cal{R}}
\cup(\omega^{\alpha})_{\alpha\in\Cal{B}}$
the corresponding dual basis in $\Hat{\frak{g}}^*$.
\par
Let $\frak{M}$ be the complex flag manifold $\Hat{\bold{G}}/\bold{Q}$
and $M$ the minimal orbit $\bold{G}/\bold{G}_+$
of $\bold{G}$ in $\frak{M}$. As usual, 
$o\simeq e\cdot\bold{G}_+\simeq e\cdot\bold{Q}$ is the base point.
We note that \par\smallskip\centerline{
$T^{1,0}_o\frak{M}\simeq \Hat{\frak{g}}/\frak{q}\simeq
\left\langle Z_{\alpha}\, | \, -\alpha\in\Cal{Q}^n\right\rangle_{\Bbb{C}} \,.$}
\par\smallskip
Therefore a Hermitian metric in $\frak{M}$ is expressed at the point
$o$ by\,:
\par\smallskip
\centerline{
$\tilde{h}_o\,=\, \dsize\sum_{\alpha,\beta\in\Cal{Q}^n}{c_{\alpha,\bar\beta}
\xi^{-\alpha}\otimes\bar{\xi}^{-\beta}}\,.$}\par\smallskip
\noindent where $(c_{\alpha,\bar\beta})$ is Hermitian symmetric
and positive definite.
For the minimal orbit we have\,:\par\smallskip\centerline{
$T^{1,0}_o{M}\simeq \bar{\frak{q}}/(\frak{q}\cap\bar{\frak{q}})\simeq
\left\langle Z_{\alpha}\, | \, -\alpha\in\Cal{Q}^n
\,,\; \alpha\in\bar{\Cal{Q}}\,\right\rangle_{\Bbb{C}} \,.$}
\par\smallskip
Thus a Hermitian metric $h$ in $T^{1,0}M$ can be represented at $o$
by\,:
\par\smallskip
\centerline{
${h}_o\,=\, \dsize\sum_{\alpha,\beta\in\Cal{Q}^n
\setminus{\bar{\Cal{Q}}}^n}{c_{\alpha,\bar\beta}
\xi^{-\alpha}\otimes\bar{\xi}^{-\beta}}\,.$}\par\smallskip
\noindent where $(c_{\alpha,\bar\beta})$ is 
again Hermitian symmetric
and positive definite. 
\par
The subspace 
$\Hat{\frak{t}}=
\sum_{\alpha\in\Cal{Q}^n\cap\bar{\Cal{Q}}^n}{\Hat{\frak{g}}^{-\alpha}}$
is a nilpotent Lie subalgebra  of $\Hat{\frak{g}}$, which is the
complexification of a real subalgebra $\frak{t}=\Hat{\frak{t}}\cap\frak{g}$
of $\frak{g}$. It can be identified to the quotient
$T_oM/H_oM$ and hence its dual space $\frak t^*$ to the stalk $H_o^0M$
of the characteristic bundle of $M$ at $o$. 
\par
From this discussion we obtain the criterion\,:
\prop{A necessary and sufficient condition for $M$ to be
essentially pseudoconcave is that there exists a positive definite
Hermitian symmetric matrix $(c_{\alpha,\bar\beta})_{\alpha,\beta\in
\Cal{Q}^n\setminus\bar{\Cal{Q}}^n}$ such that
\form{\dsize\sum_{\smallmatrix
\alpha,\beta\in
\Cal{Q}^n\setminus\bar{\Cal{Q}}^n\\
\alpha+\bar\beta=\gamma
\endsmallmatrix}{c_{\alpha,\bar\beta}[Z_{\alpha},\bar Z_{\beta}]}
=0\qquad\forall \gamma\in\Cal{Q}^n\cap\bar{\Cal{Q}}^n\, .}
\xdef\formtracc{\number\q.\number\t}}
\dimo Indeed (\formtracc) ie equivalent to the formula we obtain by
changing $\alpha,\beta,\gamma$ into $-\alpha,-\beta,-\gamma$.
\enddimo
Denote by $\check\frak{T}^{1,0}$ the $\Bbb{C}$-linear
subspace of $\Hat{\frak{g}}$ with basis
$(Z_\alpha)_{\alpha\in\Cal{Q}^n\setminus\bar\Cal{Q}^n}$.
To each $\gamma\in\Cal{Q}^n\cap\bar{\Cal{Q}}^n$ we associate a
complex-valued form of type $(1,1)$ in $\check\frak{T}^{1,0}$\,:
\form{\bold{L}_{\gamma}:
\check\frak{T}^{1,0}\times\check\frak{T}^{1,0}\ni (Z,W)
@>>>\bold{L}_{\gamma}(Z,W)=(1/i)\kappa_{\Hat{\frak{g}}}(Z_{-\gamma},
[Z,\bar W]) \in\Bbb{C}\, ,}
where $\kappa_{\Hat{\frak{g}}}$ is the Killing form in
$\Hat{\frak{g}}$. When $\gamma=\bar\gamma$ is real, 
we take $Z_{-\gamma}$ in $\frak g$,
to obtain a Hermitian symmetric $\bold L_\gamma$. 
\par
We obviously have\,:
\lem{The following are equivalent\,:
\roster
\item"($i$)" $M=M(\frak g,\frak q)$ is essentially pseudoconcave\,;
\item"($ii$)" There exists a Hermitian symmetric positive definite
form $\bold{h}$ in $\check\frak{T}^{1,0}$ such that all
$\bold{L}_{\gamma}$, for $\gamma\in\Cal{Q}^n\cap\bar\Cal{Q}^n$
have zero trace with respect to $\bold{h}$\,;
\item"($iii$)" For each $\gamma\in\Cal{Q}^n\cap\bar\Cal{Q}^n$
the Hermitian quadratic forms in $\check\frak{T}^{1,0}$\,:
\form{\check\frak{T}^{1,0}\ni Z @>>> \Re\bold{L}_{\gamma}(Z,\bar Z)\in
\Bbb{R}
\quad\text{and}\quad 
\check\frak{T}^{1,0}\ni Z @>>> \Im\bold{L}_{\gamma}(Z,\bar Z)\in\Bbb{R}}
are either $0$ or have at least one positive and one negative eigenvalue.
\endroster}
\dimo The equivalence was proved in \cite{HN2}.
\enddimo
\prop{Let $(\frak g,\frak q)$ be an effective parabolic minimal 
fundamental $CR$ algebra.
A necessary and sufficient condition for $M=M(\frak g,\frak q)$
to be essentially pseudoconcave is that for all real roots
$\gamma\in\Cal{Q}^n\cap\bar{\Cal{Q}}^n$ the Hermitian symmetric
form $\bold L_\gamma$ is either zero or has at least one positive
and one negative eigenvalue.}
\dimo
The condition is obviously necessary. We prove
sufficiency. Let $\Gamma$ be a subset of $\Cal{Q}^n\cap\bar\Cal{Q}^n$
and let $\Cal{H}(\Gamma)$ the $\Bbb{R}$-linear space
consisting of the Hermitian symmetric parts of all linear combinations
$\sum_{\gamma\in\Gamma}{a_\gamma\bold{L}_{\gamma}}$ with
$a_\gamma\in\Bbb{C}$. 
When $\gamma\in\Cal{Q}^n\cap\bar\Cal{Q}^n$ is not real, the 
Hermitian symmetric part $h$ of $a\bold{L}_\gamma$, for $a\in\Bbb{C}$,
satisfies $h(Z_\alpha,\bar Z_\alpha)=0$ for all 
$\alpha\in\Cal{Q}^n\setminus\bar\Cal{Q}^n$.
More generally, 
if $\Gamma_0$ is the set of all $\gamma\in
\Cal{Q}^n\cap\bar\Cal{Q}^n$ for which
$\sum_{\alpha\in\Cal{Q}^n\setminus\bar\Cal{Q}^n}{\bold{L}_{\gamma}
(Z_\alpha,\bar Z_\alpha)}=0$, 
then the matrices
$(h(Z_\alpha,\bar Z_\beta))_{\alpha,\beta\in\Cal{Q}^n\setminus\bar\Cal{Q}^n}$
corresponding to $h\in\Cal{H}(\Gamma_0)$ have 
zero trace and thus every  $h\in\Cal{H}(\Gamma_0)$ that is
$\neq 0$ has at least one positive and one negative eigenvalue.
\par
Choose $\Gamma$ as a maximal subset of $\Cal{Q}^n\cap\bar\Cal{Q}^n$
that contains $\Gamma_0$ and has the property that all non zero
$h\in\Cal{H}(\Gamma)$ have at least one positive and one negative
eigenvalue.
\par
If $\Gamma=\Cal{Q}^n\cap\bar\Cal{Q}^n$, then $M(\frak g,\frak q)$
is essentially pseudoconcave. Assume by contradiction that 
there is $\gamma\in\Cal{Q}^n\cap\bar\Cal{Q}^n\setminus\Gamma$. 
\par
Then
$\gamma$ is real, $\bold L_\gamma$ is Hermitian symmetric and
$\Cal{H}(\Gamma\cup\{\gamma\})=\Cal{H}(\Gamma)+\Bbb{R}\cdot
\bold{L}_\gamma$\,. Moreover,
there is at least one root 
$\alpha_0\in\Cal{Q}^n\setminus\bar\Cal{Q}^n$ such that
$\gamma=\alpha_0+\bar\alpha_0$. 
Assume that there is another root $\alpha_1\in
\Cal{Q}^n\setminus\bar\Cal{Q}^n$ with 
$\alpha_1+\bar\alpha_1=\gamma$ and 
$\bold{L}_{\gamma}(Z_{\alpha_0},\bar Z_{\alpha_0})\cdot
\bold{L}_{\gamma}(Z_{\alpha_1},\bar Z_{\alpha_1})<0$.
If $h\in\Cal{H}(\Gamma)$, then  
$h(Z_{\alpha_0},\bar Z_{\alpha_0})=h(Z_{\alpha_1},\bar Z_{\alpha_1})=0$.
Then the matrix  associated 
in the basis $(Z_{\alpha})$ to
a linear combinations $h+c\bold L_{\gamma}$ with
$c\in\Bbb{R}$, $c\neq 0$, 
has two entries of opposite sign on the main diagonal
and therefore at least one negative and one positive eigenvalue.
This would contradict the maximality of $\Gamma$. Hence we must
assume that all terms $\bold{L}_\gamma(Z_\alpha,\bar Z_\alpha)$
have the same sign. \par
By the assumption that $\bold L_\gamma$ has at least one positive
and one negative eigenvalue, we deduce that there are roots
$\beta_1,\beta_2\in\Cal{Q}^n\setminus\bar\Cal{Q}^n$ such that
$\beta_2\neq\bar\beta_1$ and $\beta_1+\bar\beta_2=\bar\beta_1+\beta_2=
\gamma$, so that $\bold{L}_\gamma(Z_{\beta_1},\bar Z_{\beta_2})\neq 0$.
If $h(Z_{\beta_2},\bar Z_{\beta_2})=0$ for all $h\in\Cal{H}(\Gamma)$,
then the matrix corresponding to $h+c\bold{L}_{\gamma}$, for
$h\in\Cal{H}(\Gamma)$, $c\in\Bbb{R}$, $c\neq 0$ in the basis
$(Z_\alpha)$ contains a principal $2\times 2$ minor
matrix, corresponding
to $\beta_1,\beta_2$, of the form
$$\pmatrix
a&\lambda\\
\bar\lambda&0
\endpmatrix \quad\text{with}\quad 
a\in\Bbb{R} \quad\text{and}\quad 
\lambda\in\Bbb{C},\;\lambda\neq 0\,.$$
Thus it would have at least one positive and one negative eigenvalue,
contradicting the choice of $\Gamma$. \par
Therefore, if 
$\Gamma\neq\Cal{Q}_{\Phi}^n\cap\bar\Cal{Q}_{\Phi}^n$,
we have\,:
\par
($i$) there exists $\alpha_0\in\Cal{Q}^n\setminus\bar\Cal{Q}^n$
such that $\alpha_0+\bar\alpha_0=\gamma\in\Cal{Q}^n\cap\bar\Cal{Q}^n$\,;
\par
($ii$) there exists $\alpha_1,\alpha_2\in\Cal{Q}^n\setminus\bar\Cal{Q}^n$
with $\alpha_2\neq\alpha_1$,
$\alpha_2\neq\bar\alpha_1$ and $\alpha_1+\bar\alpha_2=\gamma$\,;
\par
($iii$) for all $\alpha,\beta\in\Cal{Q}^n\setminus\bar\Cal{Q}^n$
with $\alpha\neq\beta$, 
$\beta\neq\bar\alpha$ and $\alpha+\bar\beta=\gamma$, we have
$\alpha+\bar\alpha\in\Cal{Q}^n\cap\bar\Cal{Q}^n$
and $\beta+\bar\beta\in\Cal{Q}^n\cap\bar\Cal{Q}^n$.
\smallskip
The roots $\alpha_0,\bar\alpha_0,\alpha_1,\bar\alpha_1,\alpha_2,\bar\alpha_2$
generate a root system $\Cal{R}'$ 
in their span in
$\frak{h}_{\Bbb{R}}^*$, that is closed under conjugation. 
Since we have the relations $\alpha_0+\bar\alpha_0=\alpha_1+\bar\alpha_2
=\alpha_2+\bar\alpha_0$, the span of $\Cal{R}'$ has dimension 
$\leq 4$. 
Moreover, $\alpha_0+\bar\alpha_0$, $\alpha_1+\bar\alpha_1$ and
$\alpha_2+\bar\alpha_2$ must be three distinct roots in $\Cal{R}'$.
Indeed, set $\alpha_1+\bar\alpha_1=\gamma_1$,
$\alpha_2+\bar\alpha_2=\gamma_2$. By assumption
$\gamma_1\neq\gamma\neq\gamma_2$. 
Moreover we obtain 
$\alpha_1-\alpha_2=\gamma_1-\gamma=\gamma-\gamma_2$,
i.e. $\gamma_1+\gamma_2=2\gamma$, which implies that
$\gamma_1\neq\gamma_2$ when $\gamma_1\neq\gamma\neq\gamma_2$. 
\par
Thus the dimension of the span of $\Cal{R}'$ is
$\leq 4$. An inspection of the Satake diagrams 
corresponding to bases of at most $4$ simple roots shows that
no such root system 
contains $3$ distinct positive real roots that are sum
of a root and its conjugate. Denote by $\Omega$ the set of 
positive real roots $\gamma$ that are of the form $\gamma=
\alpha+\bar\alpha$ with $\alpha\in\Cal{R}$.
To verify our claim, we only need to consider
the diagrams with $\ell=3,4$ and $\Omega\neq\emptyset$\,:
$$\matrix\format\c &\qquad\l\\
\frak{su}(1,3):&\Omega=\{\alpha_1+\alpha_2+\alpha_3\}\\
\frak{su}(2,2):&\Omega=\{\alpha_1+\alpha_2+\alpha_3\}\\
\frak{su}(1,4):&\Omega=\{\alpha_1+\alpha_2+\alpha_3+\alpha_4\}\\
\frak{su}(2,3):&
\Omega=\{\alpha_2+\alpha_3,\alpha_1+\alpha_2+\alpha_3+\alpha_4\}
\endmatrix
\tag $\roman{A\,IIIa,IIIb}$ $$
$$\matrix\format\c &\qquad\l\\
\frak{sp}(1,2):&\Omega=\{\alpha_1+2\alpha_2+\alpha_3\}\\
\frak{sp}(1,3):&\Omega=\{\alpha_1+2\alpha_2+2\alpha_3+\alpha_4\}
\endmatrix
\tag $\roman{C\,IIa}$
$$
Thus we obtained a contradiction,
proving our statement.
\enddimo

\thm{Let $(\frak g,\frak q_{\Phi})$ be a simple effective
and fundamental parabolic minimal $CR$ algebra. Then 
$M(\frak g,\frak q_{\Phi})$ is always essentially pseudoconcave if
$\frak g$ is either of the complex type, or compact, or of real 
type $\roman{A\,II}$, $\roman{A\,IIIb}$,
$\roman{B}$, $\roman{C\,IIb}$, $\roman{D\,I}$, $\roman{D\,II}$, 
$\roman{D\,IIIa}$,
$\roman{E\,II}$, $\roman{E\,IV}$, $\roman{E\,VI}$, 
$\roman{E\,VII}$, $\roman{E\,IX}$. 
In the remaining cases $M(\frak g,\frak q_{\Phi})$ is
essentially pseudoconcave if and only if 
we have one of the following\,:
$$\cases
\Phi\subset\Cal{R}_{\bullet}\\
\Phi\subset\{\alpha_i\,| i<p\}\cup\{\alpha_i\,| i>q\}
\endcases\;
\tag{A\,IIIa-IV} $$
$$\cases
\Phi\subset\{\alpha_{2h-1}\,|\, 1\leq h\leq p\}\\
\Phi\subset\{\alpha_i\,| i>2p\}
\endcases\;
\tag{C\,IIa}$$
$$\Phi\cap\{\alpha_{\ell-1},\alpha_\ell\}=\emptyset
\tag{D\,IIIb}
$$
$$
\cases
\{\alpha_4\}\subset\Phi\subset\Cal{R}_{\bullet}\\
\Phi=\{\alpha_3,\alpha_5\}
\endcases\;
\tag{E\,III}
$$
$$
\Phi\subset\{\alpha_1,\alpha_2\}
\tag{F\,II} $$
}\edef\teoremalc{\number\q.\number\x}
[See the table of Satake diagrams for the 
types and the references to the roots in the statement.]
\dimo
We exclude in the statement the split forms, because in these
cases $(\frak g,\frak q)$ is not fundamental.
When $\frak{g}$ is compact, $(\frak g,\frak q)$ is totally
complex and thus essentially pseudoconcave, since
the condition on the Levi form is trivially fulfilled.
\par
For $\frak{g}$ of the complex types or of the real
types $\roman{A\,II}$, $\roman{A\,IIIb}$,
$\roman{B}$, $\roman{C\,IIb}$, $\roman{D\,I}$, $\roman{D\,II}$, 
$\roman{D\,IIIa}$, $\roman{E\,II}$, $\roman{E\,IV}$,
$\roman{E\,VI}$, 
$\roman{E\,VII}$, $\roman{E\,IX}$ 
the statement follows from the fact that $\Cal{Q}^n\cap\bar{\Cal{Q}}^n$
cannot possibly contain a root of the form $\alpha+\bar\alpha$ with
$\alpha\in\Cal{Q}^n\setminus\bar{\Cal{Q}}^n$.
\smallskip
To discuss the remaining cases, we shall use the following\,:
\lem{Let $\frak g$ be a semisimple real Lie algebra, with a Cartan
decomposition, $\frak{g}=\frak{k}\oplus\frak{p}$,
and
$\frak{h}$ a Cartan subalgebra which is invariant with respect to
the corresponding Cartan involution $\vartheta$ and with maximal vector
part. Denote by $\sigma$ the conjugation of $\Hat{\frak{g}}$ with
respect to the real form $\frak g$ and let $\tau=\sigma\circ\vartheta$
the conjugation with respect the compact form $\frak{k}\oplus\,i\,\frak{p}$
of $\Hat{\frak{g}}$. Set
$\Cal{R}=\Cal{R}(\Hat{\frak{g}},\hat{\frak{h}})$. 
Then there exists a Chevalley system 
$\{X_\alpha\}_{\alpha\in\Cal{R}}$ with $X_\alpha\in\Hat{\frak{g}}^{\alpha}$
such that\,:
$$\cases
[X_{\alpha},X_{-\alpha}]=-H_{\alpha}\quad\forall\alpha\in\Cal{R}\\
[H_\alpha,X_{\beta}]=\beta(H_\alpha) X_\beta\\
[X_{\alpha},X_{\beta}]=N_{\alpha,\beta}X_{\alpha+\beta}\\
\tau(X_{\alpha})=\sigma(X_{\alpha})=\bar{X}_\alpha=X_{-\alpha}\quad
\forall \alpha\in\Cal{R}_{\bullet}
\endcases$$ 
where the $H_\alpha$ and the coefficients $N_{\alpha,\beta}$ satisfy\,:
$$\cases
\beta(H_\alpha)=q-p\\
N_{\alpha,\beta}=\pm(q+1)\\
N_{\alpha,\beta}\cdot N_{-\alpha,\alpha+\beta}=-p(q+1)\\
\roman{if}\; \beta-q\alpha,\hdots,\beta+p\alpha\;\text{is the
$\alpha$-string through $\beta$.}
\endcases
$$
}\edef\lemmale{\number\q.\number\x}\demo{Proof of Lemma {\lemmale}}
For the proof of this lemma we refer the reader to \cite{B2, Ch.VIII}, or
\cite{He, Ch.III}. \enddimo
\lem{With the notation of  Lemma {\lemmale}\,: let $\alpha,\beta\in\Cal{R}$,
with $\alpha\in\Cal{R}_{\bullet}$, and
$\alpha+\beta\in\Cal{R}$, $\alpha-\beta\notin\Cal{R}$,
$\beta+\bar\beta\in\Cal{R}$. Let
\par\smallskip\centerline{
$\beta, \,\hdots,\,\beta+p\alpha$ \qquad and\qquad
$\beta+\bar\beta-q'\alpha,\,\hdots,\,\beta+\bar\beta+p'\alpha$}\par
\smallskip\noindent
be the $\alpha$-strings through $\beta$ and $\beta+\bar\beta$, respectively.
Then we have\,:
\form{\left[X_{\alpha+\beta},\bar{X}_{\alpha+\beta}\right]\, = \, 
\left[[X_{\alpha},X_{\beta}],\overline{[X_{\alpha},X_{\beta}]}\right]
=\left(p-p'(1+q')\right)\,[X_{\beta},\bar{X}_\beta]\, 
.}\xdef\formulald{\number\q.\number\t}}\edef\lemmalf{\number\q.\number\x}
\demo{Proof of Lemma {\lemmalf}}
We observe that $[X_{\alpha},X_{\beta}]=\pm X_{\alpha+\beta}$,
because $\beta-\alpha\notin\Cal{R}$.
We have\,:
$$\matrix\format\r&\l\\
\left[X_{\alpha+\beta},\bar X_{\alpha+\beta}\right]&=
\left[[X_{\alpha},X_{\beta}],\overline{[X_{\alpha},X_{\beta}]}\right]
=\left[[X_{\alpha},X_{\beta}],[X_{-\alpha},\bar X_{\beta}]\right]\\
\vspace{2\jot}
&=[[[X_{\alpha},X_{\beta}],X_{-\alpha}],\bar X_{\beta}]+
[X_{-\alpha},[[X_\alpha,X_\beta],\bar X_\beta]]\\\vspace{2\jot}
&=[[[X_{\alpha},X_{-\alpha}],X_{\beta}],\bar X_{\beta}]
+[X_{-\alpha},[X_{\alpha},[X_\beta,\bar X_{\beta}]]]\\\vspace{2\jot}
&= \left(-\beta(H_{\alpha})+
N_{\alpha,\beta+\bar\beta}N_{-\alpha,\beta+\bar\beta
+\alpha}\right)\,[X_\beta,\bar X_\beta]\,, 
\endmatrix
$$
which, by Lemma {\lemmale}, yields (\formulald).\enddimo
\smallskip
\noindent
{\it Continuation of the Proof of Theorem {\teoremalc}.}\par
We proceed by a case by case analysis of the 
simple
real Lie algebras containing real roots $\gamma$ of
the form $\gamma=\alpha+\bar\alpha$.
\smallskip
\noindent 
$\boxed{\roman{A\,IIIa-IV}}$ \quad The positive real 
roots that are of the
form $\alpha+\bar\alpha$ for some $\alpha\in\Cal{R}$ are\,:
\smallskip\centerline{$
\gamma_h=\sum_{j=h}^{p+q-h}{\alpha_j}$ \quad
for \; $h=1,\hdots,p$.}
\par \medskip
$(i)$\quad  Assume that 
$\Phi\subset\Cal{R}_{\bullet}$. All $\gamma_h$'s belong
to 
$\Cal{Q}_{\Phi}^n\cap\bar\Cal{Q}_{\Phi}^n$ and are sums 
$\alpha+\bar\alpha$ with $\alpha\in\Cal{Q}_{\Phi}^n\setminus
\bar\Cal{Q}_{\Phi}^n$. 
To prove that $\bold{L}_{\gamma_h}$ has at least one positive and
one negative eigenvalue, we consider the roots
$\beta=\sum_{j=h}^{q-1}{\alpha_j}$ and
$\delta=\sum_{j=p+1}^{p+q-h}$. 
They both belong to 
$\Cal{Q}_{\Phi}^n\setminus
\bar\Cal{Q}_{\Phi}^n$ and
$\beta+\bar\beta=\delta+\bar\delta=\gamma_h$.
We have $\delta=\bar\beta+\eta$ with $\eta\in\Cal{R}_{\bullet}$
and $\bar\beta-\eta\notin\Cal{R}$. Since 
$\gamma_h\pm\eta\notin\Cal{R}$,
by Lemma {\lemmalf} we obtain\,:
$$[X_\delta,\bar X_{\delta}]= [[X_{\bar\beta},X_{\eta}],
\overline{[X_{\bar\beta},X_{\eta}]}]= -\, [X_\beta,\bar X_\beta]\, .$$
\smallskip
$(ii)$\quad Assume that $\Phi\cap\left(\Cal{R}_{\bullet}
\cup\{\alpha_p,\alpha_q\}\right)=\emptyset$. Let
$\Phi=\{\alpha_{j_1},\hdots,\alpha_{j_r},\alpha_{h_1},\hdots,\alpha_{h_s}\}$
with $1\leq j_1<\cdots<j_r<p<q<h_1<\cdots<h_s\leq \ell=p+q-1$.
We can assume that $r\geq 1$ and, if $s\geq 1$, that
$p-j_r<h_1-q$. Let $h_1'=p+q-h_1$
if $s\geq 1$, and $h'_1=0$ otherwise. 
The real roots in $\Cal{Q}_{\Phi}^n\cap\bar\Cal{Q}_{\Phi}^n$ are the
$\gamma_k$'s with $1\leq k\leq j_r$.
All $\bold{L}_{\gamma_k}$'s with $k\leq h'_1$ are $0$.
To show that the $\bold{L}_{\gamma_k}$'s with $h'_1<k\leq j_r$ have
at least one positive and one negative eigenvalues, we consider
$\alpha=\sum_{i=k}^{j_r}{\alpha_i}$ and
$\beta=\sum_{i=k}^{\ell-j_r}{\alpha_i}$. 
They both belong to $\Cal{Q}_{\Phi}^n\setminus
\bar\Cal{Q}_{\Phi}^n$, are distinct, and
$\alpha+\bar\beta=\gamma_k$.
\smallskip
$(iii)$\quad When $\Phi\cap\{\alpha_p,\alpha_q\}\neq\emptyset$,
we can 
assume, modulo a $CR$ isomorphism,
that $\alpha_p\in\Phi$. Then $\gamma_p
\in\Cal{Q}_{\Phi}^n\cap\bar\Cal{Q}_{\Phi}^n$ and all pairs
$(\alpha,\beta)$ of roots in $\Cal{Q}_{\Phi}^n\setminus
\bar\Cal{Q}_{\Phi}^n$ with $\alpha+\bar\beta=\gamma_p$ are of the form
$(\beta_k,\beta_k)$ with $\beta_k=\sum_{i=p}^k{\alpha_i}$ for some
$p\leq k<q$. By Lemma {\lemmalf}, we have
$$[X_{\beta_k},\bar X_{\beta_k}]=[X_{\alpha_p},\bar{X}_{\alpha_p}]\,,
$$
and hence the corresponding 
$\bold{L}_{\gamma_p}$ is $\neq 0$ and semi-definite.
\smallskip
$(iv)$\quad Assume that
$\Phi\cap\Cal{R}_{\bullet}\neq\emptyset$ and
$\Phi\not\subset\Cal{R}_{\bullet}$. We can assume,
modulo a $CR$ isomorphism, that there is
$\alpha_j\in\Phi$ with $j\leq p$ and that $\alpha_i\notin\Phi$ if
either $j<i\leq p$, or $q\leq i\leq p+q-j$. Let $r$ be the largest
integer $<q$ such that $\alpha_r\in\Phi$. We observe that
$\gamma_j
\in\Cal{Q}_{\Phi}^n\cap\bar\Cal{Q}_{\Phi}^n$
and that all pairs
$(\alpha,\beta)$ of roots in $\Cal{Q}_{\Phi}^n\setminus
\bar\Cal{Q}_{\Phi}^n$ with $\alpha+\bar\beta=\gamma_j$ are of the form
$(\beta_k,\beta_k)$ with $\beta_k=\sum_{i=j}^k{\alpha_i}$ for some
$r\leq k<q$. As in the previous case, for all $p\leq k<q$\,: 
$$[X_{\beta_k},\bar X_{\beta_k}]=[X_{\beta_p},\bar{X}_{\beta_p}]\,,
$$
and hence
$\bold{L}_{\gamma_j}$ is $\neq 0$ and semi-definite.
\medskip
\noindent
$\boxed{\roman{C\,IIa}}$\qquad The positive real roots that can be
written a sum $\alpha+\bar\alpha$ with $\alpha\in\Cal{R}$ are\,:
\par\smallskip\centerline{
$\gamma_h=
\alpha_{2h-1}+\alpha_{\ell}+2\sum_{i=2h}^{\ell-1}{\alpha_i}$
\quad for\quad $h=1,\hdots,p$.}
\par\medskip
$(i)$\quad Assume that $\Phi=\{\alpha_{2h_1-1},\hdots,\alpha_{2h_r-1}\}$
with $1\leq h_1<\cdots<h_r\leq p$. The roots in 
$\Cal{Q}_{\Phi}^n\cap\bar\Cal{Q}_{\Phi}^n$ that are of the form
$\alpha+\bar\alpha$ are the $\gamma_h$ with $1\leq h\leq h_r$.
The root $\gamma_{h_r}$ is the only one that can be written as
$\alpha+\bar\alpha$ 
with $\alpha\in\Cal{Q}_{\Phi}^n\setminus\bar\Cal{Q}_{\Phi}^n$.
But this root can also be written as $\alpha+\bar\beta$ with
$\alpha=\alpha_{2h_r-1}+\gamma_{h_r}$ and $\beta=\alpha_{2h_r-1}$,
and therefore $\bold{L}_{\gamma_{h_r}}$ has at least one positive
and one negative eigenvalue.
\par\smallskip
$(ii)$\quad Assume that $\Phi=\{\alpha_{k_1},\hdots,\alpha_{k_r}\}$
with $2p< k_1<\cdots<k_r\leq \ell$. Then all $\gamma_h$ belong
to $\Cal{Q}_{\Phi}^n\cap\bar\Cal{Q}_{\Phi}^n$. Fix $1\leq h\leq p$,
and consider the roots $\beta=\sum_{i=2h}^{2p}{\alpha_i}+\alpha_{\ell}+
2\sum_{i=2p+1}^{\ell-1}{\alpha_i}$ and $\alpha=\alpha_{2h-1}$.
Then $\beta,\alpha+\beta\in
\Cal{Q}_{\Phi}^n\setminus\bar\Cal{Q}_{\Phi}^n$ and
$\beta+\bar\beta=(\alpha+\beta)+\overline{(\alpha+\beta)}=\gamma_h$.
By Lemma {\lemmalf} we have\,:
\par\smallskip\centerline{
$[X_{\alpha+\beta},\bar X_{\alpha+\beta}]=
[[X_{\alpha},X_{\beta}],[X_{-\alpha},\bar X_{\beta}]]=
-[X_{\beta},\bar X_{\beta}] $\,,}\par\smallskip\noindent
showing that $\bold{L}_{\gamma_h}$ has at least one positive and
one negative eigenvalue.
\par\smallskip
$(iii)$\quad 
Assume that $\Phi\supset\{\alpha_{2h-1},\alpha_k\}$ with
$1\leq h\leq p$ and $k>2p$. We can take
$h$ to be the largest integer $\leq p$
with $\alpha_{2h-1}\in\Phi$ and $k$ to be the smallest integer $>2p$
with $\alpha_k\in\Phi$. Then 
$\gamma_h\in\Cal{Q}_{\Phi}^n\cap\bar\Cal{Q}_{\Phi}^n$.
The set of pairs $(\alpha,\beta)$ of elements of
$\Cal{Q}_{\Phi}^n\setminus\bar\Cal{Q}_{\Phi}^n$ with
$\alpha+\bar\beta=\gamma_h$ consists of the pairs $(\beta_r,\beta_r)$,
where\,: \par\smallskip\centerline{
$\beta_r=\sum_{i=2h-1}^{r}{\alpha_i}+\alpha_{\ell}+
2\sum_{i=r}^{\ell-1}{\alpha_i}$}\par\smallskip\noindent
for $r=2p+1,\hdots,k$. We observe that 
$\beta_{r}=\beta_{r+1}+\alpha_r$, and that $\gamma_h\pm\alpha_r\notin\Cal{R}$.
Hence by Lemma {\lemmalf} we have\,:\par\smallskip\centerline{
$[X_{\beta_{r}},\bar X_{\beta_{r}}]=
[[X_{\alpha_r},X_{\beta_{r+1}}],[X_{-\alpha_r},\bar X_{\beta_{r+1}}]]=
[X_{\beta_{r+1}},\bar X_{\beta_{r+1}}] $\,,}\par\smallskip\noindent
for all $r=2p+1,\hdots,k-1$. Hence
$\bold{L}_{\gamma_h}$ is $\neq 0$ and semi-definite.
\medskip
\noindent
$\boxed{\roman{D\,IIIb}}$ \quad The positive real roots 
that can be written as $\alpha+\bar\alpha$ with
$\alpha\in\Cal{R}$ are\,:\par\smallskip\centerline{
$\gamma_h=\alpha_{2h-1}+\alpha_{\ell-1}+\alpha_{\ell}+2\sum_{i=1}^{\ell-2}
{\alpha_i}$,\quad
$h=1,\hdots,p$,\quad for\quad $p=\frac{\ell-1}{2}$.} \par
\medskip
$(i)$\quad Assume that $\Phi\cap\{\alpha_{\ell-1},\alpha_{\ell}\}\neq
\emptyset$. Then $\gamma_{p}\in\Cal{Q}_{\Phi}^n\cap\bar\Cal{Q}_{\Phi}^n$,
and the same discussion of case $\roman{A\,IV}$ shows that
$\bold{L}_{\gamma_{p}}$ is $\neq 0$ and semi-definite.
\smallskip
$(ii)$\quad Assume that 
$\Phi=\{\alpha_{2h_1-1},\hdots,\alpha_{2h_r-1}\}$ with
$1\leq h_1<\cdots<h_r\leq p$. Then 
$\gamma_1,\hdots,\gamma_{h_r}\in\Cal{Q}_{\Phi}^n\cap\bar\Cal{Q}_{\Phi}^n$,
but only $\gamma_{h_r}$ can be represented as a sum
$\alpha+\bar\alpha$ with 
$\alpha\in\Cal{Q}_{\Phi}^n\setminus\bar\Cal{Q}_{\Phi}^n$.
If $h_r=p$, we reduce to the case of $\roman{A\,IV}$.
Assume that $h_r<p$. Then we consider the 
two distinct roots\,:
\par\smallskip\centerline{
$\beta=\alpha_{2h_r-1}+\alpha_{2h_r}$\qquad and\qquad 
$\delta=\beta+\gamma_{h_{r+1}}$\,. }\par\smallskip\noindent
They
both belong to $\Cal{Q}_{\Phi}^n\setminus\bar\Cal{Q}_{\Phi}^n$
and $\beta+\bar\delta=\gamma_{h_r}$, showing that 
$\bold{L}_{\gamma_{h_r}}$ has at least one positive and one negative
eigenvalue.
\medskip
\noindent
$\boxed{\roman{E\,III}}$ \quad
Set $\gamma_1=\alpha_1+\alpha_3+\alpha_4+\alpha_5+\alpha_6$,\quad
$\gamma_2=\alpha_1+2\alpha_2+2\alpha_3+3\alpha_4+2\alpha_5+\alpha_6$.
These are the real 
positive roots in $\Cal{R}$ that can be written as
a sum $\alpha+\bar\alpha$ for a root $\alpha\in\Cal{R}$.
Note that $\gamma_1,\gamma_2$ both  belong
to $\Cal{Q}_{\Phi}^n\cap\bar\Cal{Q}_{\Phi}^n$ for
every choice of $\Phi$. 
The discussion of the signature of $\bold{L}_{\gamma_1}$ 
reduces to the one we did for $\roman{A\,IV}$.
\medskip
$(i)$\quad
Assume that $\Phi\cap\{\alpha_1,\alpha_6\}\neq\emptyset$. In this case
the discussion for $\roman{A\,IV}$ shows that 
$\bold{L}_{\gamma_1}$ is $\neq 0$ and semi-definite.
\smallskip
$(ii)$\quad Assume that $\Phi=\{\alpha_3\}$ (the case $\Phi=\{\alpha_5\}$
is analogous). Then the set of pairs $(\alpha,\beta)$ of roots of
$\Cal{Q}_{\Phi}^n\setminus\bar\Cal{Q}_{\Phi}^n$ such that
$\alpha+\bar\beta=\gamma_2$ contains only the pair $(\alpha,\alpha)$
with $\alpha=\alpha_1+\alpha_2+2\alpha_3+2\alpha_4+\alpha_5$.
Hence $\bold{L}_{\gamma_2}$ has rank $1$ and is $\neq 0$ and
semi-definite.
\smallskip
$(iii)$\quad
Assume that either
$\alpha_4\in\Phi\subset\Cal{R}_{\bullet}$, or
$\Phi=\{\alpha_3,\alpha_5\}$. 
Then the set of pairs $(\alpha,\beta)$ of roots of
$\Cal{Q}_{\Phi}^n\setminus\bar\Cal{Q}_{\Phi}^n$ such that
$\alpha+\bar\beta=\gamma_2$ is empty, so that
$\bold{L}_{\gamma_2}=0$. The discussion 
for $\roman{A\,IV}$ shows that in this case
$\bold{L}_{\gamma_1}$ has one positive and one negative eigenvalue.
\medskip
\noindent
$\boxed{\roman{F\,II}}$ \quad
The real root $\gamma=\alpha_1+2\alpha_2+3\alpha_3+2\alpha_4$
is the only positive root which can be written in the form 
$\alpha+\bar\alpha$ for some $\alpha\in\Cal{R}$. It
belongs to $\Cal{Q}_{\Phi}^n\cap\bar\Cal{Q}_{\Phi}^n$
for every choice of $\Phi$.
\smallskip
$(i)$\quad
Assume that $\alpha_3\in\Phi$. 
Then $(\alpha_4,\alpha_4)$ is the only pair $(\alpha,\beta)$ of roots 
in $\Cal{Q}_{\Phi}^n\setminus\bar\Cal{Q}_{\Phi}^n$ with
$\alpha+\bar\beta=\gamma$. Thus $\bold{L}_{\gamma}$ has
rank $1$ and hencefore is $\neq 0$ and semi-definite.
\smallskip
$(ii)$\quad Assume that $\Phi\subset\{\alpha_1,\alpha_2\}$. 
Set $\beta=\alpha_1+2\alpha_2+2\alpha_3+\alpha_4=\bar\alpha_4-\alpha_3$
and $\alpha=\alpha_3$. Then
$\beta$ and $\beta+\alpha$ both belong to 
$\Cal{Q}_{\Phi}^n\setminus\bar\Cal{Q}_{\Phi}^n$. With the notation
of Lemma {\lemmalf}, we have $p=1$, $p'=1$, $q'=1$. Thus\,:
\par\smallskip\centerline{$
[X_{\alpha+\beta},\bar X_{\alpha+\beta}]
=[[ X_{\alpha},X_{\beta}],[X_{-\alpha},\bar X_{\beta}]]
=-[X_{\beta},\bar X_{\beta}]$\,,}\par\smallskip\noindent
showing that $\bold{L}_{\gamma}$ has at least one positive and
one negative eigenvalue.
\enddimo

%% file: 14min.tex
\se{$CR$ functions on minimal orbits}
Let $M$ be a $CR$ manifold and let $\Cal O_M(M)$ be the space of 
smooth $CR$ functions on $M$
(see \S{\sectpreliminaries}).
We say that $M$ is {\it locally $CR$ separable} if the functions of 
$\Cal O_M(M)$ locally
separate points, and
{\it $CR$ separable} if the functions in $\Cal O_M(M)$ separate 
points.
\par
In this section we discuss $CR$ separability for
the minimal orbit $M=M(\frak g,\frak q)$ 
associated to
the parabolic minimal $CR$ algebra $(\frak g,\frak q)$.\par
When $\frak{g}=\frak{g}_1\oplus\frak{g}_2$ is the product
of two semisimple ideals, we set
$\frak{q}_i=\frak{q}\cap\Hat{\frak{g}}_i$, for $i=1,2$.
By Proposition {\teoremaeb} each
$\frak{q}_i$ is parabolic in
$\Hat{\frak{g}}_i$, and the
$(\frak{g}_i,\frak{q}_i)$'s are parabolic minimal.
By the remarks following Proposition {\teoremaeb}, we have
$M\simeq M_1\times M_2$
and therefore $\Cal{O}_{M_1}(M_1)\otimes\Cal{O}_{M_2}(M_2)$ is
dense in 
$\Cal O_M(M)$.
\par
When $M$ is totally real, all smooth functions in $M$ are $CR$ and
the $CR$ separability is trivial. Thus
the $CR$ separability of
a general $M(\frak g,\frak q)$ 
reduces to that of the fibers of its
fundamental reduction (cf. Theorem {\teofundamentalb} and its Corollary).
\par
Indeed, let $M@>{\rho}>>N$ be  the fundamental reduction. We have
$\Cal{O}_N(N)=\Cal{C}^{\infty}(N,\Bbb{C})$. Thus the $CR$ functions
certainly separate points on distinct fibers. Furthermore, 
let 
$f$ 
be a $CR$ function defined on a fiber $\rho^{-1}(x_0)$ ($x_0\in N$).
We choose a $CR$ trivialization 
$U\times \rho^{-1}(x_0)\ni
(x,y)@>>>\phi(x,y)\in\rho^{-1}(U)$, where
$U$ is an open neighborhood of $x_0$ in $N$
and $\phi$ is a $CR$ diffeomorphism. Then $f$ extends to a
$CR$ function $F$ in $\rho^{-1}(U)$, by $F(\phi(x,y))=f(y)$.
Take a cut-off function $\chi\in\Cal{C}^\infty(U,\Bbb{C})$,
with compact support in $U$ and 
equal to $1$ in $x_0$. Then $\tilde{f}(z)=\chi(\rho(z))\cdot F(z)$
for $z\in\rho^{-1}(U)$,
extended by $\tilde{f}=0$ outside $\rho^{-1}(U)$, is a $CR$ function in $M$
that extends $f$. This shows that $M$ is (locally) $CR$  separable
if and only if the fiber $\rho^{-1}(x_0)$ is 
(locally) $CR$ separable.
\par
In this way we can restrain our discussion of $CR$ separability
to the case where $(\frak g,\frak q)$ 
is simple, effective and fundamental.
\par
First we indicate how 
$CR$ separability can be read off the cross-marked Satake diagram.
We have\,:
\thm{Let $(\frak g, \frak q_\Phi)$ be a simple fundamental effective parabolic 
minimal $CR$ algebra and $M=M(\frak g,\frak q_\Phi)$ its associated $CR$ 
manifold.
Then $M$ is locally $CR$ separable if and only if its 
cross-marked Satake diagram
is one of the following\,:
$$\matrix\format\c\quad&\l\\
(\roman{A\,IIIa-IV}) &\Phi\subset\{\alpha_p,\alpha_q\}\\
(\roman{D\,IIIb}) &\Phi\subset\{\alpha_{\ell-1},\alpha_\ell\}\\
(\roman{E\,III}) &\Phi\subset\{\alpha_1,\alpha_6\}\, .
\endmatrix
$$
%
In these cases $M$ is also $CR$ separable by real analytic $CR$ functions.}
\edef\teoCRsep{\number\q.\number\x}
\dimo
We prove that if 
there is a simple root $\alpha$ satisfying either\,:
\roster
\item"$(i)$"  $\alpha\in\Phi\cap\Cal R_\bullet$, or
\item"$(ii)$" $\alpha\in\Phi\setminus\Cal R_\bullet$,
$\bar\alpha\in\Cal B$ and $\alpha+\bar\alpha\not\in\Cal R$,
\endroster
then $M$ is not $CR$ separable. Inspection of the 
Satake diagrams then shows that
the only possibilities left are those listed above. Finally in the 
examples below we show that in those cases $M$ has a $CR$ embedding
into an affine complex space $E$, hence is separable by analytic 
$CR$ functions
(the restrictions to $M$ of the global
holomorphic functions in $E$).
\par
Let $\alpha$ be a
simple root, satisfying either $(i)$ or $(ii)$. 
Then $\alpha\pm\bar\alpha\not\in\Cal R$.
Denote by $\hat\frak a$ the subalgebra generated by $\hat\frak g^{\alpha}+ 
\hat\frak g^{-\alpha}+\hat\frak g^{\bar\alpha}
+\hat\frak g^{-\bar\alpha}$. It is semisimple, $\frak q\cap\hat\frak a$ is 
parabolic in $\hat\frak a$ and $(\frak a,\frak q\cap\hat\frak a)$ is 
totally complex. 
\par
Let 
$\hat\frak b=\hat\frak a\cap\hat\frak k$ (in case $(i)$ 
we have $\hat\frak b=\hat\frak a$). 
Then $\frak b$ is compact semisimple and $(\frak b,\frak k\cap\hat\frak b)$ 
is totally complex. If $\bold B\subset\bold G$ is the analytic subgroup with 
Lie algebra $\frak b$ then $\bold B\cdot o$ is a compact complex submanifold 
of $M$ of
positive dimension. All smooth $CR$ functions on $M$ 
restrict to holomorphic functions in
$\bold B\cdot o$, that are constant in
$\bold B\cdot o$ by Liouville's Theorem, and therefore
$M$ is not $CR$ separable.
\enddimo
\exam{Fix positive integers $p<q$ and let $n=p+q$.
We identify the simple real Lie algebra $\frak{g}\simeq\frak{su}(p,q)$
with the set of $(n\times n)$ complex matrices $Z$ with zero trace
that satisfy\,:\par\smallskip\centerline{
$Z^* K\,+\,K\, Z\, = \, 0\, $, where 
$K=\left({\smallmatrix
I_p\\
&-I_q
\endsmallmatrix}\right)\,.$}\par\smallskip
Let $e_1,\hdots,e_{n}$ be the canonical basis of $\Bbb{C}^{n}$
and let $\frak{q_{\alpha_p}}\subset\Hat{\frak{g}}\simeq\frak{sl}
(n,\Bbb{C})$ be the set of $(n\times n)$ matrices in $\frak{sl}
(n,\Bbb{C})$ such that \par\smallskip\centerline{
$Z(\langle e_1+e_{p+1}\, ,\hdots\,,e_{p}+e_{2p}\rangle)\subset
\langle e_1+e_{p+1}\, ,\hdots\,,e_{p}+e_{2p}\rangle\, .$}\par
\smallskip\noindent
Then $(\frak{g},\frak{q_{\alpha_p}})$ is parabolic minimal. 
\par
The corresponding $CR$ manifold $M=M(\frak{g},\frak{q_{\alpha_p}})$ is the 
Grassmannian of $p$-planes $\ell_p$ in $\Bbb{C}^{n}$ which are totally 
isotropic for $K$ (i.e. $v^*Kv=0$ for all $v\in\ell_p$). 
We have\par\smallskip
\centerline{
$M\simeq\big\{\ell_p=\{(v,u(v))\in\Bbb{C}^{n}\, | \,
v\in\Bbb{C}^p\}\,\big| \, u\in\bold{U}(\Bbb C^p,\C^q)\big\}
\simeq\bold{U}(\C^p,\C^q)$}\par\smallskip\noindent
where $\bold{U}(\C^p,\C^q)=\{u\in\Cal{M}_{q\times p}(\Bbb{C})\, | \, 
u^*u=I_p\}$  is the set of unitary $q\times p$ matrices.\par
Give $\bold{U}(\C^p,\C^q)$ the $CR$ structure induced by the 
embedding in $\Cal{M}_{q\times p}(\Bbb{C})$.
The compact subgroup $\bold K^{(1)}\simeq\bold{SU}(p)\times\bold{SU}(q)$
of matrices of $\bold{SU}(p,q)$ of the form $\left(\smallmatrix A_p&0\\0&B_q
\endsmallmatrix\right)$ acts transitively by $CR$
automorphisms on $\bold{U}(\C^p,\C^q)$, the action being given by:
$\left(\smallmatrix A&0\\0&B\endsmallmatrix\right)\cdot u=BuA^{-1}$.
\par
The associated $CR$ algebra is $(\frak k^{(1)},\frak q')$ where 
$\frak k^{(1)}\simeq\frak{su}(p)\oplus\frak{su}(q)$ and
$\frak q'$ is the set of matrices in $\frak{sl}(p)\oplus\frak{sl}(q)$
of the form $\left(\smallmatrix A_p&0&0\\0&A_p&D\\0&0&C_{q-p}
\endsmallmatrix\right)$.
\par
The group $\bold K^{(1)}$ acts transitively on $M$ by Theorem {\teoremadb}, 
and the 
associated $CR$ algebra is $(\frak k^{(1)},\Hat{\frak k}^{(1)}\cap
\frak q)=(\frak k^{(1)},\frak q')$. Thus the 
diffeomorphism $M\simeq\bold{U}(\C^p,\C^q)$ is in fact a $CR$
isomorphism.
\par
In this way we obtain the embedding
$M\hookrightarrow\Cal{M}_{q\times p}(\Bbb{C})
\simeq\Bbb{C}^{qp}$. This is a $CR$ 
embedding into a Stein manifold, and therefore $M$ is separable,
since $\Cal{O}_M(M)$ contains the
restrictions to $M$ of the
holomorphic functions in $\Cal{M}_{q\times p}(\Bbb{C})\simeq\Bbb{C}^{qp}$.}
\par\edef\esempiomb{\number\q.\number\x}
\exam{Fix a positive integer $p$ and let $n=2p+1$.
We identify the simple real Lie algebra
$\frak g=\frak{so}^*(2n)$ with the set of $(2n\times 2n)$ complex 
matrices $Z$ with zero trace
that satisfy\,:\par\smallskip\centerline{
$\cases
ZJ=J\bar Z\, ,\\
{^t\! Z} K+KZ = 0\,,
\endcases$}\par\noindent
where:\par\smallskip\centerline{ 
$J=\left({\smallmatrix
0&I_n\\
-I_n&0\endsmallmatrix}\right)\quad$ 
and
$\quad K=\left({\smallmatrix
0&I_n\\
I_n&0\endsmallmatrix}\right)\,. \quad$}\par\smallskip\noindent
Let $\frak q_{\alpha_1}$ be the parabolic subalgebra of matrices in 
$\frak g$ that stabilize the subspace
\par\smallskip\centerline{ 
$V_n=\langle e_1+e_{n+2p},\dots,e_p
+e_{n+p+1},e_{p+1}-e_{n+p},\dots,e_{2p}-e_{n+1},e_{2p+1}\rangle.$}
\par\noindent
Then $(\frak g, \frak q_{\alpha_1})$ is parabolic minimal.
\par
The maximal compact subgroup $\bold K\simeq\bold U(n)$ of $\bold G$ of 
matrices of the form\,:\hfill\linebreak 
$\left(\smallmatrix A_n&0\\0&{^t\! A}^{-1}_n\endsmallmatrix\right)$, 
$A_n\in\bold U(n)$, acts transitively by $CR$ isomorphisms on 
$M(\frak g,\frak q_{\alpha_1})$.
The associated $CR$ algebra is $(\frak k,\frak q')$ where 
$\frak k\simeq\frak u(n)$ and $\frak q'=\Hat{\frak k}\cap\frak q_{\alpha_1}$.
This is the subalgebra of matrices in $\frak{so}(2n,\C)$ of the form 
$\left(\smallmatrix A_n&0\\0&-{^t\! A}_n\endsmallmatrix\right)$ where 
$A_n\in\frak{gl}(n,\C)$ is of the form 
$\left(\smallmatrix B_p&C_p&v_p\\D_p&-{{^t\! B}_p}
&w_p\\0&0&is\endsmallmatrix\right)$ with $B_p={{^t\! B}_p}$, 
$D_p={{^t\! D}_p}$.
\par
We let $\bold K$ act on $\frak{so}(n,\C)$ by\,: 
$k\cdot X=A\, X\, {^t\!\! A}$ if 
$k=\left(\smallmatrix A_n&0\\0&{^t\! A}^{-1}_n\endsmallmatrix\right)$.
Let $N$ be the $\bold K$-orbit of 
$o=\left(\smallmatrix 0&-I_p&0\\I_p&0&0\\0&0&0\endsmallmatrix\right)$.
The associated $CR$ algebra is $(\frak k,\frak q')$ and the isotropy is 
connected and contains a generator of $\pi_1\big(\bold U(n)\big)$.
Thus $M$ is $CR$ isomorphic to $N$.
\par
Since $N$ is an embedded $CR$ submanifold of the Stein manifold
$\frak{so}(n,\C)\simeq\Bbb{C}^{n(n-1)/2}$, it follows that
$M$ is separable by the restrictions to $M$ of  the
holomorphic functions in $\frak{so}(n,\C)\simeq\Bbb{C}^{n(n-1)/2} $.}
\edef\esempiomc{\number\q.\number\x}
\exam{Let $D$ be the exceptional bounded symmetric domain of type V.
Its Shilov boundary $S$ is a real flag manifold 
(see \cite{Fa, Part III,Ch.IV\S 2.8}) for the
group E\,III and is compact, hence it is a minimal orbit $M(\frak g,\frak q)$ 
where $\frak g$ is of type E\,III.
Furthermore it has $CR$ dimension $8$ and $CR$ codimension $8$ 
(see \cite{KZ, p.\ 180}), hence $\frak q=\frak q_{\alpha_1}$ or 
$\frak q=\frak q_{\alpha_6}$.
Thus $M(\frak g,\frak q_{\alpha_1})\simeq S$ is an embedded $CR$ submanifold 
of $\C^{16}$ and it is separable,
since $\Cal{O}_M(M)$ contains the restrictions of the holomorphic functions
in $\C^{16}$.\xdef\esempiomd{\number\q.\number\x}
\par
A similar argument could have been 
applied also to discuss  the two previous examples
{\esempiomb} and {\esempiomc}.
Indeed the three classes of minimal orbits
of examples {\esempiomb}, {\esempiomc} and {\esempiomd}
are exactly the Shilov boundaries
of the  bounded
symmetric domains
that are not of tube type (and
that are not totally real; see also \cite{He, Ch.X, Ex. D.1}
and \cite{Hi}).}
\par
Identify $\Hat{\frak{g}}$ with a complex Lie algebra of
left invariant complex valued vector fields in $\bold{G}$.
Given a $\Bbb{C}$ linear subspace $\frak{a}$ of $\Hat{\frak{g}}$,
denote by 
$\Cal O_{\bold G,\frak a}$ the sheaf of smooth complex 
valued functions $f$ 
on $\bold G$ 
such that $L(f)=0$ for all $L\in\frak a$.
\par
Let $\pi:\bold{G}@>>>M$ be the principal fibration,
$\Cal F_M(\bold{G})
=\pi^*\Cal O_M(M)$ and $\Cal T_M$ the sheaf of local smooth 
complex 
valued vector fields $L$ on $\bold G$ such that $L(f)=0$ for all $f\in
\Cal F_M(\bold{G})$. 
Let $\frak q'=\{X\in\Hat{\frak{g}}|X_e\in(\Cal T_M)_e\}$. We note
that, since $\Cal{F}_M(\bold{G})$ is invariant for the left action of
$\bold{G}$ on functions, 
the elements of $\frak{q}'$ define left invariant
global sections of $\Cal T_M$.
\lem{$\frak q'$ is a parabolic subalgebra of $\hat\frak g$ and $\frak q
\subset\frak q'$.}\edef\lemCRsepprim{\number\q.\number\x}
\dimo
The sheaf $\Cal T_M$ is invariant for the left action of $\bold G$,
hence it is generated at every point by the global left invariant complex 
vector fields that belong to $\frak q'$.
Since $\Cal T_M$ is involutive, $\frak q'$ is  a subalgebra.
Clearly it contains $\frak q$, and thus is parabolic.
\enddimo
\par
\lem{We have: $\Cal{F}_M(\bold{G})=\Cal{O}_{\bold G,\frak q'}(\bold G)
=\Cal{O}_{\bold G,\frak q}(\bold G)$\,.}
\edef\lemCRsepdue{\number\q.\number\x}
\dimo
The inclusions $\Cal{F}_M(\bold{G})\subset\Cal{O}_{\bold G,\frak q'}(\bold G)
\subset\Cal{O}_{\bold G,\frak q}(\bold G)$
follow from the definition of $\frak{q}'$ and 
Lemma {\lemCRsepprim}. To complete the proof we need to show that
$\Cal O_{\bold G,\frak q}(\bold G)\subset\Cal F_M(\bold{G})$.
An $f\in\Cal O_{\bold G,\frak q}(\bold G)$
is constant on the analytic subgroup 
of $\bold{G}$ that has Lie algebra 
$\frak g\cap\frak q$, i.e. on $\bold G_+$, 
(recall that $\bold G_+$ 
is connected). By left invariance, $f$ is constant on the left cosets
$g\bold{G}_+$ and hence is
the pullback of a function 
$\tilde f$ defined in $M$. 
Furthermore $\tilde f$ is $CR$ on $M$ because $T^{0,1}_oM=\pi_*(\frak q)$
(here $\Hat{\frak{g}}$ is identified with the complexification
of $T_e\bold{G}$) and therefore
$L(\Tilde{f})=0$ for all $L\in T^{0,1}M$ by left invariance.
\enddimo
Let $M'=M(\frak g,\frak q')$, $\pi':\bold G@>>>M'$ 
and 
$\rho=\pi'\circ\pi^{-1}:M@>>>M'$ the natural projection.
We have the commutative diagram\,:
$$\xymatrix {
&{\bold{G}}\ar[dl]_{\pi}\ar[dr]^{\pi'}\\
M\ar[rr]_{\rho}&& M'}
$$
\lem{$M'$ is $CR$ separable.}
\dimo
Assume that $M'$ is not locally $CR$ separable.
Then there exists a tangent vector $X\in\frak g$ at $e\in\bold G$ such that 
$\pi'_*(X)\neq 0$ and $\pi'_*(X)(f)=0$ for
every $CR$ function $f$ on $M'$. 
Then $X\in\frak q'$, hence $X\in\frak g_+'$ and $\pi'_*(X)=0$, 
yielding a contradiction.
\par
By Theorem {\teoCRsep} 
the $CR$ manifold $M'$, being locally $CR$ separable,
is also $CR$ separable.
\enddimo
\lem{$\Cal F_M(\bold{G})=\Cal F_{M'}(\bold{G})$. Hence\,: 
$\Cal O_M(M)=\rho^*\big(O_{M'}(M')\big)$.}
\dimo
This follows by applying Lemma {\lemCRsepdue} to $M'$.
\enddimo
\lem{The complex Lie subalgebra 
$\frak q'$ is minimal in the set of (not necessarily proper) parabolic
subalgebras of $\Hat{\frak{g}}$, 
containing $\frak q$, such that $M(\frak g,\frak q)$ is 
$CR$ separable.}
\dimo
Suppose that $\frak{q}''$ is a parabolic subalgebra of 
$\Hat{\frak{g}}$ with 
$\frak q\subset\frak q''\subset\frak q'$ and such that
$M''=M(\frak g,\frak q'')$ is $CR$ separable.
Then $\Cal F_{M''}(\bold G)=\Cal F_{M'}(\bold G)$
by Lemma {\lemCRsepdue}.
This implies that 
$\frak q''\cap\bar\frak q''=\frak q'\cap\bar\frak q'$,
yielding $M''=M'$ and thus $\frak{q}''=\frak{q}'$.
\enddimo
The discussion above leads to the following\,:
\thm{Let $(\frak g,\frak q_\Phi)$ be a simple effective fundamental parabolic 
minimal $CR$ algebra and $M=M(\frak g,\frak q_\Phi)$.
Then there exists $\Psi\subset\Phi$ and a $\bold G$-equivariant fibration 
$\rho:M@>>>M_s=M(\frak g,\frak q_\Psi)$ such that $M_s$ is $CR$ separable and
$\Cal O_M(M)=\rho^*\big(O_{M_s}(M_s)\big)$.
\par
Furthermore $\Psi=\Phi\cap\Sigma$, where $\Sigma$ is defined according
to the type of $\frak g$\,:
\roster\widestnumber\item{"{\rom Type A\,IIIa:}"}
\item"{\rom Type A\,IIIa :}" $\Sigma=\{\alpha_p,\alpha_{q}\}$;
\item"{\rom Type D\,IIIb :}" $\Sigma=\{\alpha_{\ell-1},\alpha_\ell\}$;
\item"{\rom Type E\,III :}" $\Sigma=\{\alpha_1,\alpha_5\}$;
\item"{\rom All other types :}" $\Sigma=\emptyset$.
\endroster
The space $\Cal O_M(M)$ is one dimensional when 
$\Psi=\emptyset$, infinite  
dimensional 
when $\Psi\neq\emptyset$.\qed}

%% file: 15min.tex
\head {} Appendix: Table of noncompact real forms and Satake diagrams\endhead
%
\def\rul{\noalign{\vskip 1mm\hrule height0.5pt\vskip 2mm}}
\def\mat#1{{\xymatrix @M=0pt @R=2pt @!C=6pt{#1}}}
$$\vbox{
\offinterlineskip
\vskip 2mm\hrule height1pt\vskip 2mm
\halign {\strut 
\quad $\roman{\bold{#}}$\hfil \quad&
\quad \hfil $#$ \hfil\quad&\quad \hfil $#$\hfil\quad\cr
\text{Name}&\frak g&\text{Satake diagram}\cr
\noalign{\vskip 2mm\hrule height1pt\vskip 2mm}
A\,I&
\frak{sl}(\ell+1,\Bbb{R})&
\mat{
{\circ}\ar @{-}[r]&{\circ}
\ar @{--}[rr]&&{\circ}\ar @{-}[r]&{\circ}\\
{\alpha_1}&&&&{\alpha_{\ell}}}\cr
\rul
A\,II&
\matrix\frak{sl}(p,\Bbb H)\\2p+1=\ell\endmatrix&
\mat{
{\bullet}\ar @{-}[r]&{\circ}\ar @{-}[r]&{\bullet}
\ar @{--}[r]&{\bullet}\ar @{-}[r]&{\circ}\ar @{-}[r]&{\bullet}\\
{\alpha_1}&&&&&{\alpha_{\ell}}}\cr
\rul
A\,IIIa&
\matrix\frak{su}(p,q)\\p+q=\ell+1\\2\leq p\leq\ell/2\endmatrix&
\vbox{\vskip16pt\mat{
{\circ}\ar @{--}[r]\ar@/^18pt/@{<->}[rrrrr]&
{\circ}\ar @{-}[r]\ar@/^8pt/@{<->}[rrr]&
{\bullet}\ar @{--}[r]&{\bullet}\ar @{-}[r]&{\circ}\ar @{--}[r]&{\circ}\\
{\alpha_1}&{\alpha_p}&&&{\alpha_{\ell-p+1}}&{\alpha_{\ell}}}}
\cr
\rul
A\,IIIb&
\matrix\frak{su}(p,p)\\1\leq p=(\ell+1)/2\endmatrix&
\vbox{\vskip16pt\mat{
{\circ}\ar @{--}[r]\ar@/^18pt/@{<->}[rrrr]&
{\circ}\ar @{-}[r]\ar@/^8pt/@{<->}[rr]&
{\circ}\ar @{-}[r]&{\circ}\ar @{--}[r]&{\circ}\\
{\alpha_1}&&{\alpha_{p}}&&{\alpha_{\ell}}}}
\cr
\rul
A\,IV&
\matrix\frak{su}(1,\ell)\endmatrix&
\vbox{\vskip12pt\mat{
{\circ}\ar @{-}[r]\ar@/^10pt/@{<->}[rrr]&
{\bullet}\ar @{--}[r]&{\bullet}\ar @{-}[r]&{\circ}\\
{\alpha_1}&&&{\alpha_{\ell}}}}
\cr
\rul
B\,I&
\matrix\frak{so}(p,2\ell+1-p)\\2\leq p\leq\ell\endmatrix&
\mat{
{\circ}\ar@{--}[r] & {\circ}\ar@{-}[r] & {\bullet}\ar@{--}[r] &
{\bullet}\ar@{=}[r]|{\SelectTips{cm}{}\dir{>}} & {\bullet} \\ 
{\alpha_1} & {\alpha_p} &&& {\alpha_\ell}}
\cr\rul
B\,II&
\matrix\frak{so}(1,2\ell)\endmatrix&
\mat{
{\circ}\ar@{-}[r] & {\bullet}\ar@{--}[rr] & &
{\bullet}\ar@{=}[r]|{\SelectTips{cm}{}\dir{>}} & {\bullet} \\ 
{\alpha_1} &&&& {\alpha_\ell}}
\cr\rul
C\,I&
\frak{sp}(2\ell,\Bbb{R})&
\mat{
{\circ}\ar@{-}[r] & {\circ}\ar@{--}[rr] & &
{\circ}\ar@{=}[r]|{\SelectTips{cm}{}\dir{<}} & {\circ} \\ 
{\alpha_1} &&&& {\alpha_\ell}}
\cr
\rul
C\,IIa&
\matrix\frak{sp}(p,\ell-p)\\2p<\ell\endmatrix&
\mat{
{\bullet}\ar @{-}[r]&{\circ}\ar @{-}[r]&{\bullet}
\ar @{--}[r]&{\circ}\ar @{-}[r]&{\bullet}\ar @{--}[r]&{\bullet}\ar@{=}[r]|{\SelectTips{cm}{}\dir{<}} & {\bullet} \\
{\alpha_1}&&&{\alpha_{2p}}&&&{\alpha_{\ell}}}
\cr\rul
C\,IIb&
\matrix\frak{sp}(p,p)\\2p=\ell\endmatrix&
\mat{
{\bullet}\ar @{-}[r]&{\circ}\ar @{-}[r]&{\bullet}
\ar @{--}[r]&{\circ}\ar @{-}[r]&{\bullet}\ar@{=}[r]|{\SelectTips{cm}{}\dir{<}} & {\circ} \\
{\alpha_1}&&&&&{\alpha_{\ell}}}
\cr
\rul
D\,Ia&
\matrix\frak{so}(p,2\ell-p)\\2\leq p\leq\ell-2\endmatrix&
\mat{
&&&&{\bullet}\\
&&&&{\alpha_{\ell-1}} \\
{\circ}\ar@{--}[r] & {\circ}\ar@{-}[r] & {\bullet}\ar@{--}[r] &
{\bullet}\ar@{-}[uur]\ar@{-}[ddr] \\ 
{\alpha_1} &  {\alpha_p} &&&{\alpha_\ell}\\
&&&& {\bullet}}
\cr}
\vskip 2mm\hrule height1pt\vskip 2mm}
$$
$$\vbox{
\offinterlineskip
\vskip 2mm\hrule height1pt\vskip 2mm
\halign {\strut 
\qquad $\roman{\bold{#}}$\hfil \quad&
\quad \hfil $#$ \hfil\quad&\quad \hfil $#$\hfil\qquad\cr
\text{Name}&\frak g&\text{Satake diagram}\cr
\noalign{\vskip 2mm\hrule height1pt\vskip 2mm}
D\,Ib&
\frak{so}(\ell-1,\ell+1)&
\mat{
&&&&{\circ}\ar@{<->}@/^13pt/[dddd]\\
&&&&{\alpha_{\ell-\!1}} \\
{\circ}\ar@{-}[r] & {\circ}\ar@{--}[rr] &&
{\circ}\ar@{-}[uur]\ar@{-}[ddr] \\ 
{\alpha_1} &&&&{\alpha_\ell}\\
&&&& {\circ}}
\cr\rul
D\,Ic&
\frak{so}(2\ell,\Bbb R)&
\mat{
&&&&{\circ}\\
&&&&{\alpha_{\ell-1}} \\
{\circ}\ar@{-}[r] & {\circ}\ar@{--}[rr] &&
{\circ}\ar@{-}[uur]\ar@{-}[ddr] \\ 
{\alpha_1} &&&&{\alpha_\ell}\\
&&&& {\circ}}
\cr\rul
D\,II&
\frak{so}(1,2\ell-1)&
\mat{
&&&&{\bullet}\\
&&&&{\alpha_{\ell-1}} \\
{\circ}\ar@{-}[r] & {\bullet}\ar@{--}[rr] &&
{\bullet}\ar@{-}[uur]\ar@{-}[ddr] \\ 
{\alpha_1} &&&&{\alpha_\ell}\\
&&&& {\bullet}}
\cr\rul
D\,IIIa&
\matrix\frak{so*}(2\ell)\\\ell=2p\endmatrix&
\mat{
&&&&{\bullet}\\
&&&&{\alpha_{\ell-1}} \\
{\bullet}\ar@{-}[r] & {\circ}\ar@{-}[r] & {\bullet}\ar@{--}[r] &
{\circ}\ar@{-}[uur]\ar@{-}[ddr] \\ 
{\alpha_1} &&&&{\alpha_\ell}\\
&&&& {\circ}}
\cr
\rul
D\,IIIb&
\matrix\frak{so*}(2\ell)\\\ell=2p+1\endmatrix&
\mat{
&&&&&{\circ}\ar@{<->}@/^13pt/[dddd]\\
&&&&&{\alpha_{\ell-\!1}} \\
{\bullet}\ar@{-}[r] & {\circ}\ar@{-}[r] & {\bullet}\ar@{--}[r] &
{\circ}\ar@{-}[r] &{\bullet}\ar@{-}[uur]\ar@{-}[ddr] \\ 
{\alpha_1} &&&&&{\alpha_\ell}\\
&&&&& {\circ}}
\cr\rul
E\,I&&
\mat{
&&{\alpha_4}\\
{\circ}\ar@{-}[r] & {\circ}\ar@{-}[r] 
& {\circ}\ar@{-}[r] & {\circ}\ar@{-}[r] & {\circ}\\
{\alpha_1}&{\alpha_3}&&{\alpha_5}&{\alpha_6}\\
\\
&&{\circ}\ar@{-}[uuu]\\
&&{\alpha_2}}
\cr\rul
E\,II&&
\vbox{\vskip12pt\mat{
&&{\alpha_4}\\
{\circ}\ar@{-}[r]\ar@/^18pt/@{<->}[rrrr] & {\circ}\ar@{-}[r]\ar@/^11pt/@{<->}[rr] 
& {\circ}\ar@{-}[r] & {\circ}\ar@{-}[r] & {\circ}\\
{\alpha_1}&{\alpha_3}&&{\alpha_5}&{\alpha_6}\\
\\
&&{\circ}\ar@{-}[uuu]\\
&&{\alpha_2}}}
\cr\rul
E\,III&&
\vbox{\vskip10pt\mat{
&&{\alpha_4}\\
{\circ}\ar@{-}[r]\ar@/^16pt/@{<->}[rrrr] & {\bullet}\ar@{-}[r] 
& {\bullet}\ar@{-}[r] & {\bullet}\ar@{-}[r] & {\circ}\\
{\alpha_1}&{\alpha_3}&&{\alpha_5}&{\alpha_6}\\
\\
&&{\circ}\ar@{-}[uuu]\\
&&{\alpha_2}}}
\cr\rul
E\,IV&&
\mat{
&&{\alpha_4}\\
{\circ}\ar@{-}[r]& {\bullet}\ar@{-}[r] 
& {\bullet}\ar@{-}[r] & {\bullet}\ar@{-}[r] & {\circ}\\
{\alpha_1}&{\alpha_3}&&{\alpha_5}&{\alpha_6}\\
\\
&&{\bullet}\ar@{-}[uuu]\\
&&{\alpha_2}}
\cr\rul
E\,V&&
\mat{
&&{\alpha_4}\\
{\circ}\ar@{-}[r] & {\circ}\ar@{-}[r] 
& {\circ}\ar@{-}[r] & {\circ}\ar@{-}[r] & {\circ}\ar@{-}[r] & {\circ}\\
{\alpha_1}&{\alpha_3}&&{\alpha_5}&{\alpha_6}&{\alpha_7}\\
\\
&&{\circ}\ar@{-}[uuu]\\
&&{\alpha_2}}
\cr\rul
E\,VI&&
\mat{
&&{\alpha_4}\\
{\circ}\ar@{-}[r]& {\circ}\ar@{-}[r] 
& {\circ}\ar@{-}[r] & {\bullet}\ar@{-}[r] & {\circ}\ar@{-}[r] & {\bullet}\\
{\alpha_1}&{\alpha_3}&&{\alpha_5}&{\alpha_6}&{\alpha_7}\\
\\
&&{\bullet}\ar@{-}[uuu]\\
&&{\alpha_2}}
\cr}
\vskip 2mm\hrule height1pt\vskip 2mm}
$$
$$\vbox{
\offinterlineskip
\vskip 2mm\hrule height1pt\vskip 2mm
\halign {\strut 
\qquad $\roman{\bold{#}}$\hfil \quad&
\quad \hfil $#$ \hfil\quad&\quad \hfil $#$\hfil\qquad\cr
\text{Name}&\frak g&\text{Satake diagram}\cr
\noalign{\vskip 2mm\hrule height1pt\vskip 2mm}
E\,VII&&
\mat{
&&{\alpha_4}\\
{\circ}\ar@{-}[r]& {\bullet}\ar@{-}[r] 
& {\bullet}\ar@{-}[r] & {\bullet}\ar@{-}[r] & {\circ}\ar@{-}[r] & {\circ}\\
{\alpha_1}&{\alpha_3}&&{\alpha_5}&{\alpha_6}&{\alpha_7}\\
\\
&&{\bullet}\ar@{-}[uuu]\\
&&{\alpha_2}}
\cr\rul
E\,VIII&&
\mat{
&&{\alpha_4}\\
{\circ}\ar@{-}[r] & {\circ}\ar@{-}[r] 
& {\circ}\ar@{-}[r] & {\circ}\ar@{-}[r] & {\circ}\ar@{-}[r] &
{\circ}\ar@{-}[r] & {\circ}\\
{\alpha_1}&{\alpha_3}&&{\alpha_5}&{\alpha_6}&{\alpha_7}&{\alpha_8}\\
\\
&&{\circ}\ar@{-}[uuu]\\
&&{\alpha_2}}
\cr\rul
E\,IX&&
\mat{
&&{\alpha_4}\\
{\circ}\ar@{-}[r]& {\bullet}\ar@{-}[r] 
& {\bullet}\ar@{-}[r] & {\bullet}\ar@{-}[r] & {\circ}\ar@{-}[r]&
{\circ}\ar@{-}[r] & {\circ}\\
{\alpha_1}&{\alpha_3}&&{\alpha_5}&{\alpha_6}&{\alpha_7}&{\alpha_8}\\
\\
&&{\bullet}\ar@{-}[uuu]\\
&&{\alpha_2}}
\cr\rul
F\,I&&
\mat{
{\circ}\ar@{-}[r] &{\circ}\ar@{=}[r]|{\SelectTips{cm}{}\dir{>}}&
{\circ}\ar@{-}[r]&{\circ}\\
{\alpha_1}&{\alpha_2}&{\alpha_3}&{\alpha_4}}
\cr\rul
F\,II&&
\mat{
{\bullet}\ar@{-}[r] &{\bullet}\ar@{=}[r]|{\SelectTips{cm}{}\dir{>}}&
{\bullet}\ar@{-}[r]&{\circ}\\
{\alpha_1}&{\alpha_2}&{\alpha_3}&{\alpha_4}}
\cr\rul
G&&
\mat{{\circ}\ar@{=}[r]|{\SelectTips{cm}{}\dir{>}}\ar@{-}[r]&{\circ}\\
{\alpha_1}&{\alpha_2}}
\cr}
\vskip 2mm\hrule height1pt\vskip 2mm}
$$

%% file: min.bbl
\Refs
\widestnumber\key{AzHuR}
\ref\key AS1
\by D.Alekseevsky, A.Spiro 
\paper Invariant CR structures on compact homogeneous manifolds
\jour Hokkaido Math. J. 
\vol 32 
\yr 2003
\pages 209-276
\endref

\ref\key AS2
\bysame 
\paper Flag manifolds and homogeneous CR structures
\inbook  Recent advances in Lie theory (Vigo, 2000)
\pages 3-44
\jour Res. Exp. Math. 
\vol 25 
\publ Heldermann
\publaddr Lemgo
\yr 2002
\endref
\ref\key AS3
\bysame 
\paper Compact homogeneous CR manifolds
\jour J. Geom. Anal. 
\vol 12 
\yr 2002
\pages 183-201
\endref 

\ref\key AlN
\by A.Altomani, M.Nacinovich
\paper Abelian extensions of semisimple
graded $CR$ algebras
\yr 2004
\pages 433-457
\jour Adv. Geom.
\vol 4
\endref

\ref\key AnF
\by A.Andreotti, G.A.Fredricks
\paper Embeddability of real analytic Cauchy-Riemann manifolds
\jour Ann. Scuola Norm. Sup. Pisa Cl. Sci.
\yr 1979
\vol  6
\pages 285-304
\endref
 
\ref\key AnH
\by A.Andreotti, C.D.Hill
\paper E.E.Levi convexity and the H.Lewy problem I, II
\jour Ann. Scuola Norm. Sup. Pisa
\yr 1972
\vol 26
\pages 352-363, 748-806
\endref

\ref\key Ar
\by S.Araki
\paper On root systems and an infinitesimal classification
of irreducible symmetric spaces
\jour J. Math. Osaka City Univ.
\vol A 13
\yr 1962
\pages 1-34
\endref

\ref\key AzHuR
\by H.Azad, A.Huckleberry, W.Richthofer
\paper Homogeneous CR-manifolds
\jour J. Reine Angew. Math. 
\vol 358 
\yr 1985
\pages 125-154. 
\endref

\ref\key BaERo
\by M.S.Baouendi, P.Ebenfelt, L.P.Rothshild
\book Real submanifolds in complex space and their mappings
\publ Princeton University Press
\yr 1999
\publaddr Princeton
\endref 

\ref\key B1
\by N.Bourbaki
\book Groupes et alg\`ebres de Lie, Ch. 4,5 et 6
\publ Masson
\yr 1981
\publaddr Paris
\endref

\ref\key B2
\by N.Bourbaki
\book Groupes et alg\`ebres de Lie, Ch. 7 et 8
\publ Masson
\yr 1990
\publaddr Paris
\endref

\ref\key CM
\by S.S.Chern, J.K.Moser
\paper  Real hypersurfaces in complex manifolds
\jour Acta Math. 
\vol 133 
\yr 1974 
\pages 219-271
\endref
 

\ref\key Fa
\by J. Faraut\ et al.
\book Analysis and geometry on complex homogeneous domains
\bookinfo Progr.\ Math.\ 185
\publ Birkh\"auser Boston
\publaddr Boston, MA
\yr 2000
\endref

\ref\key GiMa
\by S.Gindikin, T.Matsuki
\paper A remark on Schubert cells and the duality of orbits on flat manifolds
\jour J. Math. Soc. Japan 
\vol 57 
\yr 2005
\pages 157-165
\endref


\ref\key He
\by S. Helgason
\book Differential geometry, Lie groups, and Symmetric Spaces
\publaddr New York
\yr 1978
\publ Academic Press
\endref



\ref\key HN1
\by C.D.Hill, M.Nacinovich
\paper Pseudoconcave $CR$ manifolds
\inbook Complex analysis and geometry
\bookinfo Lect. Notes in Pure and Appl. Math.
\vol 173
\pages 275-297
\publ Marcel and Dekker
\publaddr New York
\yr 1996
\endref

\ref\key HN2
\bysame
\paper A weak pseudoconcavity condition for abstract almost
$CR$ manifolds
\jour Invent. Math.
\vol 142
\yr 2000
\pages 251-283
\endref

\ref\key HN3
\bysame
\paper Weak pseudoconcavity and the maximum modulus principle
\jour Ann. Mat. Pura Appl. (4) 
\vol 182 
\yr 2003 
\pages 103-112.
\endref

\ref\key HN4
\bysame
\paper Fields of CR meromorphic functions
\jour Rend. Sem. Mat. Univ. Padova
\vol 111 
\yr 2004
\pages 179-204
\endref


\ref\key Hi
\by U.Hirzebruch
\paper \"Uber Jordan Algebren und beschr\"ankte symmetrische Gebiete
\jour Math. Z.
\vol 94
\yr 1966
\pages 387-390
\endref

\ref\key Hu
\by A.Huckleberry
\paper On certain domains in cycle spaces of flag manifolds
\jour Math. Ann. 
\vol 323 
\yr 2002
\pages 797-810
\endref

\ref\key HuW
\by A.Huckleberry, J.Wolf
\paper Schubert varieties and cycle spaces
\jour Duke Math. J. 
\vol 120 
\yr 2003
\pages 229-249
\endref
\ref\key K
\by W.Kaup
\paper On the local equivalence of homogeneous CR-manifolds
\jour Arch. Math. (Basel) 
\vol 84 
\yr 2005 
\pages 276-281
\endref

\ref\key KZ
\by W.Kaup, D.Zaitsev
\paper On symmetric Cauchy-Riemann manifolds
\jour Adv. Math. 
\vol 149
\yr 2000
\pages 145--181
\endref

\ref\key LN
\by A.Lotta, M.Nacinovich
\paper On a class of symmetric CR manifolds
\jour Adv. Math.  
\vol 191  
\yr 2005  
\pages 114-146
\endref

\ref\key MeN1
\by C.Medori, M.Nacinovich 
\paper Levi-Tanaka algebras and homogeneous $CR$ manifolds
\jour Compositio Math.
\yr 1997
\vol 109
\pages 195-250
\endref

\ref\key MeN2
\bysame
\paper Classification of semisimple Levi-Tanaka algebras
\jour Ann. Mat. Pura Appl.
\yr 1998
\vol 174
\pages 285-349
\endref

\ref\key MeN3
\bysame
\paper Standard $CR$ manifolds of codimension $2$
\jour Transformation groups
\yr 2001
\vol 6
\pages 53-78
\endref

\ref\key MeN4
\bysame
\paper Complete nondegenerate locally
standard $CR$ manifolds
\yr 2000
\jour Math. Ann.
\pages 509-526
\vol
317
\endref



\ref\key MeN5
\bysame
\paper Algebras of infinitesimal $CR$ automorphisms
\jour Journal of Algebra
\yr 2005
\vol 287 
\pages 234-274
\endref

\ref\key Mo
\by G.D.Mostow
\paper On maximal subgroups of real Lie groups
\jour Ann. Math.
\yr 1961
\vol 74
\pages 503-517
\endref


\ref\key O
\ed A.L. Onishchik
\book Lie groups and Lie algebras I
\publaddr Berlin
\yr 1993
\publ Springer
\bookinfo Encyclopaedia of Mathematical Sciences \vol 20
\endref


\ref\key OV
\eds A.L. Onishchik, E.B. Vinberg
\book Lie groups and Lie algebras III
\publaddr Berlin
\yr 1994
\publ Springer
\bookinfo Encyclopaedia of Mathematical Sciences \vol 43
\endref

\ref\key St
\by N.Stanton
\paper Homogeneous CR submanifolds in $\bold C\sp n$
\jour Math. Z. 228 
\yr 1998
\pages 247-253
\endref


\ref\key T1
\by N.Tanaka
\paper On generalized graded Lie algebras and geometric structures I
\jour J. Math. Soc. Japan
\vol 19
\yr 1967
\pages 215-254
\endref

\ref\key T2
\bysame
\paper On differential systems, graded Lie algebras and pseudogroups
\jour J. Math. Kyoto Univ.
\vol 10
\yr 1970
\pages 1-82
\endref

\ref\key Tr
\by F.Treves
\book Hypo-analytic structures
\yr 1992
\publ Princeton University press
\vol 40
\bookinfo Princeton Mathematical Series
\publaddr Princeton, NJ
\endref

\ref\key V
\by E.B.Vinberg
\paper The Morozov-Borel theorem for real Lie groups
\jour Dokl. Akad Nauk SSSR
\vol 141
\yr 1961\pages 270-273
\transl\nofrills English transl. in \jour Sov. Math. Dokl.
\vol 2 \pages 1416-1419
\endref

\ref\key Wa
\by G.Warner
\book Harmonic analysis on semi-simple Lie groups I
\yr 1972
\publ Springer
\publaddr Berlin
\endref

\ref\key Wi
\by M.Wiggerman
\paper The fundamental group of a real flag manifold
\yr 1998
\jour Indag. Math.
\vol 9
\pages 141-153
\endref


\ref\key W1
\by J.A.Wolf
\paper The action of real semisimple Lie groups on a complex
manifold I: Orbit structure and holomorphic arc components
\jour Bull. Amer. Math. Soc.
\vol 75
\yr 1969
\pages 1121-1137
\endref

\ref\key W2
\bysame
\paper Real groups transitive on complex flag manifolds
\jour Proc. Amer. Math. Soc. 
\vol 129 
\yr 2001
\pages 2483--2487
\endref

\ref\key Zi
\by R.Zierau
\paper Representations in Dolbeault cohomology
\inbook Representation Theory of Lie groups
\bookinfo IAS/PARK CITY Mathematics series vol.8
\publ AMS
\eds J.Adams, D.Vogan
\publaddr Providence, Rhode Island
\pages 89-146
\yr 2000
\endref

\endRefs
